\documentclass[10pt,letterpaper]{amsart} 

\usepackage{color}
\usepackage{multicol}
\usepackage{empheq}
\usepackage[letterpaper, margin=0.8in]{geometry}
\usepackage{subfig}

\usepackage{listings}
\usepackage{amsmath, amsthm, amssymb, wasysym, verbatim, bbm, color, graphics}

\usepackage{url}
\usepackage{mathtools}
\usepackage[colorlinks,linkcolor=blue]{hyperref}
\usepackage{romannum}

\usepackage{xfrac}
\usepackage{float}
\usepackage{amsmath,bm}
\usepackage[title]{appendix}

\newcommand{\R}{\mathbb{R}}

\newcommand{\N}{\mathbb{N}}

\DeclareMathOperator{\tr}{\text{tr}}

\newcommand\innprod[2]{\langle #1,#2 \rangle}
\newcommand\norm[1]{\left\|#1\right\|}
\newcommand{\Grad}{\nabla}
\newcommand{\Div}{\operatorname{div}}
\newcommand{\dom}{\Omega}

\newcommand{\hQ}{\hat{Q}} 
\newcommand{\Us}{U^{\sigma}}
\newcommand{\Qs}{Q^{\sigma}}
\newcommand{\rs}{r^{\sigma}}
\newcommand{\hU}{\hat{U}}
\newcommand{\hr}{\hat{r}}

\newtheorem{lemma}{Lemma}[section]
\newtheorem{proposition}[lemma]{Proposition}
\newtheorem{corollary}[lemma]{Corollary}
\newtheorem{theorem}[lemma]{Theorem}
\newtheorem{remark}[lemma]{Remark}

\newtheorem*{maintheorem*}{Main Theorem}
\theoremstyle{definition}{\newtheorem{definition}[lemma]{Definition}}
\theoremstyle{note}{
\newtheorem*{claim*}{Claim}}

\newtheorem{rem}[lemma]{Remark}

\allowdisplaybreaks

\numberwithin{equation}{section}

\title[Q-tensor zero inertia limit]{The Zero Inertia Limit for the Q-Tensor Model of Liquid Crystals: Analysis and Numerics}
\date{\today}

\author[M. Hirsch]{Max Hirsch}
\address[Max Hirsch]{\newline Department of Mathematics \newline UC Berkeley \newline Berkeley, CA, 94720, USA.}
\email[]{mhirsch@berkeley.edu}
\author[F. Weber]{Franziska Weber}
\address[Franziska Weber]{\newline Department of Mathematics \newline UC Berkeley \newline Berkeley, CA, 94720, USA.}
\email[]{fweber@berkeley.edu}
\author[Y. Yue]{Yukun Yue}
\address[Yukun Yue]{\newline Department of Mathematics \newline University of Wisconsin-Madison \newline 500 Lincoln Dr, Madison, WI 53706, USA.}
\email[]{yyue@math.wisc.edu}

\thanks{F.W. and M.H. were supported in part by NSF grant DMS-2042454 and M.H. and Y.Y. were supported in part by NSF grant DMS-1912854. This material is based upon work supported by the National Science Foundation Graduate Research Fellowship Program under Grant No. DGE 2146752. Any opinions, findings, and conclusions or recommendations expressed in this material are those of the authors and do not necessarily reflect the views of the National Science Foundation.}
\begin{document}
 \pagenumbering{arabic}
\maketitle
\begin{abstract}
The goal of this work is to rigorously study the zero inertia limit for the Q-tensor model of liquid crystals. Though present in the original derivation of the Ericksen-Leslie equations for nematic liquid crystals, the inertia term of the model is often neglected in analysis and applications. We show wellposedness of the model including inertia and then show using the relative entropy method that solutions of the model with inertia converge to solutions of the model without inertia at a rate ${\sigma}$ in $L^\infty(0,T;H^1(\dom))$, where $\sigma$ is the inertial constant. Furthermore, we present an energy stable finite element scheme that is stable and convergent for all $\sigma$ and study the zero inertia limit numerically. We also present error estimates for the fully discrete scheme with respect to the discretization parameters in time and space.
\end{abstract}

\section{Introduction}
Liquid crystals are an important class of materials with unique properties that lie between those of conventional liquids and solid crystals. They exhibit a remarkable combination of fluidity and substantial order, which has led to their widespread use in various technological applications such as display devices, sensors, and optical components~\cite{Stewart2008,Castellano2005}. In order to understand and study their properties better, researchers have derived various mathematical models to describe their dynamics and be able to simulate them. One of the first such models was the director field model by Oseen and Frank~\cite{Oseen1933TheTO} and received extensive attention \cite{adler2015energy,adler2016constrained,kinderlehrer1993elementary}. Then Ericksen and Leslie extended it to a full hydrodynamic theory of liquid crystal flows~\cite{Ericksen1962HydrostaticTO,ERICKSEN1976233}. Later, Landau and de Gennes developed the so-called Q-tensor model~\cite{de1993physics}. The Q-tensor theory replaces the directer field $n\in\mathbb{S}^{d-1}$, which is used in both the Oseen-Frank model and the Ericksen-Leslie model, by a symmetric, trace-free $d\times d$ tensor $Q$ that describes the second moment of the probability distribution of liquid crystal orientations~\cite{borthagaray2020structure,bouck2024hydrodynamical,SonnetVirga2001,swain2024linear,Virga1995}.
The tensor $Q$ is then assumed to minimize the energy functional \begin{equation}\label{eq:LG}
    E_{LG}(Q)=\int_\Omega\mathcal{F}_B(Q)+\mathcal{F}_E(Q)
\end{equation}
in the equilibrium case. Here $\Omega$ denotes a bounded domain in $\R^d$ with smooth or polygonal boundary, $\mathcal{F}_B$ is the bulk potential, and $\mathcal{F}_E$ is the elastic energy density~\cite{Mottram2014IntroductionTQ}: 

\begin{equation}\label{eq:Energy_definition}
	\begin{split}
	\mathcal{F}_B(Q) &= \frac{a}{2}\tr(Q^2)-\frac{b}{3}\tr(Q^3)+\frac{c}{4}\tr^2(Q^2),\\
	\mathcal{F}_E(Q)&= \frac{1}{2}\sum_{i,j,k=1}^d\left[ L_1(\partial_k Q_{ij})^2 + L_2 \partial_j Q_{ij}\,\partial_kQ_{ik} + L_3\partial_jQ_{ik}\,\partial_kQ_{ij} \right], 
	\end{split}\end{equation}
where $L_1, L_2, L_3, a, b, c$ are constants with $c, L_1>0$, $L_2+L_3\geq0$. The non-equilibrium case can be described by taking the variational derivative  of the Landau-de Gennes energy functional $E_{LG}(Q)$ and considering the following gradient flow resulting from the Euler-Lagrange equations~\cite{BallJohn,HuangDing2015}: 
\begin{equation}\label{eq:parabolic}
\begin{split}
\frac{\partial Q_{ij}}{\partial t} = M &\Bigg[L_1 \Delta Q_{ij} + \frac{L_2+L_3}{2} \left(\sum_{k=1}^d\left( \partial_{ik}Q_{jk} + \partial_{jk}Q_{ik}\right) - \frac{2}{d} \sum_{k,s=1}^d \partial_{ks}Q_{ks}\, \delta_{ij}\right) \\
&- \left(a Q_{ij} - b\left((Q^2)_{ij} -\frac{1}{d} \tr(Q^2)\,\delta_{ij}\right)+c \tr(Q^2)Q_{ij}\right) \Bigg].
\end{split}
\end{equation}
Here $M>0$ is a constant that we will assume to be equal to one for the reminder of the paper (the analysis for $M\neq 1$ carries over identically).
This model excludes inertial effects coming from the kinetic energy $\sigma \int |Q_{t}|^2$ of the Q-tensor.
However, while the inertia is often small in liquid crystal applications and therefore disregarded, there are specific cases and applications where the effects of inertia become significant, and the behavior of liquid crystals is affected. For example, the effect of inertia plays an important role when liquid crystals are used as the medium for certain sound waves in acoustic devices~\cite{refId0, Kapustina2004-he, Lisin1995-xc}.  Further experimental evidence that the inertial constant matters was found by Guozhen~\cite{Guozhen1982}.
Ericksen and Leslie's original model also includes the inertia terms.
Leslie~\cite{Leslie1979} showed that the inertial
constant plays a nontrivial role when the anisotropic axis of the liquid crystal molecules is subjected to large accelerations.

 The Ericksen-Leslie equations are based on the director field model, though a similar hydrodynamic system including the effects of inertia based on the Q-tensor model was more recently derived by Qian and Sheng~\cite{Qian1998}. The goal of this paper is therefore to study system~\eqref{eq:parabolic} with an additional inertia term $\sigma Q_{tt}$, as well as the system's rigorous zero inertia limit $\sigma\to 0$. We will also present and analyze a numerical scheme that is stable and convergent in the discretization for any choice of $\sigma$. To be specific, we will study the system

\begin{equation}\label{eq:qtensorflow}
\begin{split}
\frac{\partial Q_{ij}}{\partial t} =&-\sigma \frac{\partial^2 Q_{ij}}{\partial t^2}+  \Bigg[L_1 \Delta Q_{ij} + \frac{L_2+L_3}{2} \left(\sum_{k=1}^d\left( \partial_{ik}Q_{jk} + \partial_{jk}Q_{ik}\right) - \frac{2}{d} \sum_{k,s=1}^d \partial_{ks}Q_{ks}\, \delta_{ij}\right) \\
&- \left(a Q_{ij} - b\left((Q^2)_{ij} -\frac{1}{d} \tr(Q^2)\,\delta_{ij}\right)+c \tr(Q^2)Q_{ij}\right) \Bigg],
\end{split}
\end{equation}
where we denote the bracketed terms on the right-hand side of~\eqref{eq:qtensorflow} 
\begin{equation}
    \label{eq:H}
    H(Q) \coloneqq L_1\Delta Q + \frac{L_2+L_3}{2}\alpha(Q) 
    - f(Q),
\end{equation}
with
\begin{equation}
    \label{eq:alpha}
    \alpha(Q)_{ij} =  \sum_{k=1}^d(\partial_{ik}Q_{jk} + \partial_{jk}Q_{ik}) -\frac{2}{d}\sum_{k,s=1}^d  \partial_{ks}Q_{ks} \,\delta_{ij},
\end{equation}
and
\begin{equation}
    \label{eq:f}
   f(Q) = \frac{\delta \mathcal{F}_B(Q)}{\delta Q} = a Q-b\left[Q^2-\frac{1}{d}\text{tr}(Q^2)I\right]+c\text{tr}(Q^2)Q.
\end{equation}
So~\eqref{eq:qtensorflow} is equivalent to 
\begin{equation}
    \label{eq:Qt}
    Q_t = -\sigma Q_{tt}+ \left(L_1\Delta Q+\frac{L_2+L_3}{2}\alpha(Q)
    -f(Q)\right) :=-\sigma Q_{tt}+  H.
\end{equation}
We complete the problem by adding the initial and boundary conditions
\begin{subequations}
    \label{eq:initial_condition}
\begin{empheq}[left=\empheqlbrace]{align}
 & Q(0,x)=Q_0(x)\quad \text{in }\Omega\\
 & Q_t(0,x)=Q_{t,0}(x)\quad \text{in }\Omega
\end{empheq}
\end{subequations}
and
\begin{equation}
    \label{eq:boundary_condition}
    Q(t,x)=0 \quad x\in\partial\Omega, t\in [0,T].
\end{equation}
We will focus on the Dirichlet boundary conditions case, though the treatment of Neumann boundary conditions is similar.
Smooth solutions of~\eqref{eq:qtensorflow} satisfy the energy dissipation law
\begin{equation}\label{eq:dissipationlaw}
    \frac{d}{dt}E(Q) = -\int_{\dom}|Q_t|^2 dx \leq 0,
\end{equation}
where $E$ is given by 
\begin{equation}
    E(Q) = \frac12\int_{\dom}\left(L_1 |\Grad Q|^2 + \sigma |Q_t|^2 +\mathcal{F}_B(Q) +(L_2 + L_3) |\Div Q|^2 \right) dx. 
\end{equation}
This can be seen by either taking the derivative of the energy, integrating by parts and plugging in the equations, or alternatively, by taking the inner product of the equations~\eqref{eq:qtensorflow} with $Q_t$, then integrating over the domain and integrating by parts where suitable.
Our goal in the following will be to show that~\eqref{eq:qtensorflow} possesses a unique strong solution, as defined later on in Definition~\ref{def:strong_solution}, under the assumption that $L_2+L_3$ are sufficiently small compared to $L_1$ and under assumptions on the regularity of the initial condition in Sobolev spaces. We will then proceed to developing a finite element scheme for~\eqref{eq:qtensorflow} and prove error estimates for the approximations generated by it.

To solve problem~\eqref{eq:Qt}--\eqref{eq:boundary_condition} numerically, we will adapt the Invariant Energy Quadratization (IEQ) method which has been introduced and extensively studied in~\cite{Yang2016, YANG2017691,YangConvergenceIEQ2020,YANG2017104, YangZhaoWangShen2017,  Zhao2017NumericalAF} to deal with gradient-flow type problems. The main idea is to introduce a new variable which is the square root of the nonlinear bulk potential. The resulting system of PDEs, which includes the evolution of the new variable, can be discretized with a second order linearly implicit time stepping scheme.
To reformulate the Q-tensor flow system, we introduce the scalar auxiliary variable $r$:
\begin{equation}\label{eq:defr}
r(Q)=\sqrt{2\left(\frac{a}{2}\tr(Q^2)-\frac{b}{3}\tr(Q^3)+\frac{c}{4}\tr^2(Q^2)+A_0\right)},
\end{equation}
where $A_0$ is a constant ensuring that $r$ is positive. This is well defined since the bulk energy is bounded from below, given $c>0$~\cite{MajumdarZarnescu2010,majumdar_2010}. It follows that
\begin{equation}\label{eq:defP}
\frac{\delta r(Q)}{\delta Q} = \frac{f(Q)}{r(Q)}:= P(Q),
\end{equation}
for symmetric, trace free tensors $Q$.
Then one can formally rewrite the gradient flow~\eqref{eq:Qt} as a system for $(Q,r)$:
\begin{subequations}\label{eq:reformulation}
	\begin{align}
	Q_t &=-\sigma Q_{tt}+ \left(L_1\Delta Q +\frac{L_2+L_3}{2}\alpha(Q)-r P(Q)\right):=-\sigma Q_{tt}+ {H},\label{eq:QZ}\\
	r_t & = P(Q):Q_t.\label{eq:rZ}
	\end{align}
\end{subequations}
Note that there is no need to prescribe boundary conditions to $r$, as it only solves an ordinary differential equation at every point in space. As an initial condition, we take $r(0,x)=r_0(x)=r(Q_0(x))$. The boundary condition for $Q$ is the same as~\eqref{eq:boundary_condition}. 
This system satisfies the energy dissipation law
\begin{equation}\label{eq:dissipationlawr}
    \frac{d}{dt}\frac12\int_{\dom} \left(L_1 |\Grad Q|^2 +(L_2+L_3)|\Div Q|^2+\sigma|Q_t|^2  +r^2 \right) dx = - \int_{\dom} |Q_t|^2 dx, 
\end{equation}
which can be easily obtained by taking the inner product of the first equation~\eqref{eq:QZ} with $Q_t$ and the second equation~\eqref{eq:rZ} with $r$, adding the two, integrating, and integrating by parts in the terms involving spatial derivatives. Due to the simplicity of this energy law, a discrete version can be obtained for a suitable time discretization, as demonstrated in Section~\ref{sec:fem}.
To discretize in space, we will use finite elements with mass lumping for the bulk potential $rP(Q)$. This will enable us to preserve a discrete version of the energy law.

Finally, we will investigate the limit $\sigma\to 0$ theoretically and numerically. We find that the solution of~\eqref{eq:qtensorflow} converges in the $L^\infty(0,T;H^1)$-norm to the solution of~\eqref{eq:parabolic} at a rate ${\sigma}$ as $\sigma\to 0$ and this rate is confirmed by our numerical experiments in Section~\ref{sec:num}.

There are a couple of related mathematical works that investigate liquid crystal models with inertia, as well as the zero inertia limit of such models. Mathematical analysis of Q-tensor models with inertia was conducted in~\cite{DEANNA20181080,Zarnescu2018_hyper,Luo2022}. Zarnescu and de Anna~\cite{DEANNA20181080} prove the global well-posedness of the inertial Qian-Sheng model in $\R^d$, $d=2,3$ for sufficiently small initial data and also prove the global existence of twist wave solutions in $\R^2$. In~\cite{Zarnescu2018_hyper}, the inviscid version of the inertial Qian-Sheng model is studied; here the Navier-Stokes equations in the coupled system reduce to the incompressible Euler equations. They prove the local existence of smooth solutions and the global existence of dissipative solutions, a solution concept related to measure-valued solutions~\cite{DiPerna1985}. They also show a weak-strong uniqueness result for dissipative solutions. Finally, Luo and Ma~\cite{Luo2022} study the zero inertia limit of the Qian-Sheng model on $\R^3$ rigorously for small initial data. They find a rate of $\sqrt{\sigma}$ which is optimal according to~\cite{Chill2004}. Our rate is better, however this may be due to the specific rescaling of the initial data that is used in~\cite{Luo2022}. If we take different initial data for~\eqref{eq:qtensorflow} and~\eqref{eq:parabolic}, then the convergence rate is the minimum of $\sigma$ and the $H^1$-norm of the difference between the initial data (for the precise statement, see Theorem~\ref{thm:sigma_convergence}). These rates are also confirmed numerically in the experiments in Section~\ref{sec:num}.
The zero inertia limit for the damped wave map, which can be viewed as the one-constant approximation of the inertial Oseen-Frank model and a simplified version of the Ericksen-Leslie equations, was studied in~\cite{Weber2018,Zarnescu2019}. In~\cite{Zarnescu2019}, the zero inertia limit is analyzed at the PDE level, whereas~\cite{Weber2018} proposes an asymptotically stable (with respect to $\sigma$) convergent numerical scheme. A one-dimensional version of the Oseen-Frank model with inertia is the nonlinear variational wave equation which has been studied extensively, see e.g.~\cite{Glassey1996,Zhang2003,Bressan2006,Bressan2015,Holden2011,Holden2009,Weber2016,Aursand2015} and references therein. Besides investigating the zero inertia limit theoretically and numerically and providing an improved convergence rate wtih respect to $\sigma$, we also present the first error estimates for a fully discrete finite element scheme for~\eqref{eq:qtensorflow}. Using the relative entropy method introduced by Dafermos~\cite{Dafermos1979} and DiPerna~\cite{Diperna1979} for hyperbolic conservation laws, we obtain a convergence rate of  order 1/2 with respect to the time discretization and order 1 with respect to the spatial discretization in the $H^1$-norm. The proof also implies an error estimate for the parabolic case $\sigma=0$, in which case we get the improved rate of 1 in the time discretiztion. The numerical experiments indicate that the rate of 1/2 in time for the $\sigma\neq 0$ case is not optimal, however, we believe to improve the rate different proof techniques are necessary.

The rest of this paper is organized as follows: In Section \ref{sec:notation}, we introduce frequently used notation and different notions of solutions of~\eqref{eq:qtensorflow}. In Section \ref{sec:strong_Solution}, we prove higher-order energy estimates, from which we derive the existence and uniqueness of a strong solution for~\eqref{eq:qtensorflow}. In Section \ref{sec:fem}, we construct a finite element scheme based on the mass lumping method and also provide error estimates. Finally, in Section~\ref{sec:zero_inertia}, we show that the solution of~\eqref{eq:Qt} converges to its formal zero inertia limit as the inertial term approaches zero, i.e., $\sigma\to0$; and in Section~\ref{sec:num}, we conclude with numerical experiments.

\section{Preliminaries}\label{sec:notation}

\subsection{Notations}

We begin by introducing notations that will be frequently used throughout this paper. For matrix-valued functions $A,B: \R^d \to \R^{d \times d}$, we define the following:

\begin{itemize}
	\item Scalar product of $A$ and $B$: $A:B = \sum_{i,j=1}^d A_{ij}B_{ij}$,
	\item Inner product of $A$ and $B$: $\innprod{A}{B} = \int_{\Omega} A:B\, dx$,
	\item Frobenius norm of $A$: $|A|:=|A|_F = \sqrt{A:A}$,
	\item $L^2$ norm of $A$: $\norm{A}^2_{L^2} = \int_{\Omega} |A|_F^2\, dx$, 
	\item Partial derivative of $A$: $\partial_i A = (\partial_i A_{jk})_{jk}$, where $\partial_i=\partial_{x_i}$,
	\item Gradient of $A$: $\nabla A=(\partial_1 A,\dots,\partial_d A)$,
	\item Norm of $\nabla A$: $|\nabla A|^2 =\sum_{i=1}^d |\partial_i A|_F^2$,
	\item $L^2$ norm of the gradient of $A$: $\norm{\nabla A}_{L^2}^2 = \int_{\Omega}|\nabla A|^2\, dx$,
 \item $H^1$ norm of the gradient of $A$: $\|A\|_{H^1}^2=\|A\|^2_{L^2}+\|\nabla A\|^2_{L^2}$.
\end{itemize}
Let $S^{(d)}$ be the space of symmetric, trace-free $d\times d$ matrices. We shall use standard notations for Sobolev spaces and Bochner spaces. 
For simplicity, we will not distinguish the scalar-valued, vector-valued, or matrix-valued functional spaces if it is clear from the context. If we write $C(0,T;X)$ with $X$ a Banach space, we are referring to a space defined on a closed interval $[0,T]$. Additionally, the notation $A_{ij,t}$ will be used to denote the time derivative of the entry $A_{ij}$ in a matrix $A$. For the constants appearing throughout the paper, we will denote the constant depending on the parameters of the Q-tensor system and time range $T$ to be $C$. As an exception, if the constant depends on $\sigma$, we will denote it as $C(\sigma)$.

\subsection{Definitions of the  weak and strong solutions}

Next, we give the definitions of weak and strong solutions of~\eqref{eq:qtensorflow}. These definitions pertain to the original system as expressed in Equation \eqref{eq:Qt} and the reformulated system detailed in Equation \eqref{eq:reformulation}.

\begin{definition}[Weak solutions of \eqref{eq:qtensorflow}]
\label{def:Q_inertia_weakQ}
		By a weak solution of~\eqref{eq:qtensorflow}, we mean a function $Q:[0,T]\times\Omega\to\R^{d\times d}$ that is trace-free and symmetric for almost every $(t,x)$ and satisfies
		\begin{equation*}
		Q\in L^\infty(0,T;H^1(\dom)),\quad Q_t\in L^2([0,T]\times\Omega)\cap L^\infty(0,T;L^2(\dom)),
		\end{equation*}
		and
		\begin{multline}
		\label{eq:weakform}
		-\int_{\dom} Q(T,x):\varphi(T,:) dx +\int_{\dom} Q_0(x):\varphi(0,x) dx +\int_0^T\int_\Omega Q:\partial_t\varphi\,dx\,dt\\
   = -\sigma\int_{0}^{T}\int_{\Omega} Q_{t}:\varphi_{t}\, dx\, dt + \sigma\int_{\Omega} Q_{t}(T, x):\varphi(T, x)\, dx - \sigma\int_{\Omega} Q_{t,0}(x):\varphi(0, x)\, dx\\ 
		+  \int_0^T\int_\Omega \Bigg(L_1\sum_{i,j=1}^d\Grad Q_{ij}\cdot\Grad\varphi_{ij} \\
  + \frac{L_2+L_3}{2} \sum_{i,j,k=1}^d \left(\partial_{k}Q_{jk}\partial_i \varphi_{ij} + \partial_{k}Q_{ik}\partial_j\varphi_{ij} -\frac{2}{d}  \partial_{i}Q_{ki} \partial_k\varphi_{jj}\right)\Bigg)dxdt \\
		+ \int_0^T\int_\Omega \left(a Q - b\left((Q^2) -\frac{1}{d} \tr(Q^2)I\right)+c \tr(Q^2)Q\right):\varphi \, dxdt,
		\end{multline}
		for all smooth $\varphi = (\varphi_{ij})_{ij=1}^d:[0,T]\times \Omega\to \R^{d\times d}$ that are compactly supported within $\dom$ for every $t\in [0,T]$.
	\end{definition}

\begin{remark}
\label{rem:L2_test_functions}
Given the density of smooth functions with compact support in $W^{1,q}([0,T];L^p(\Omega))$ and $L^q([0,T];W^{1,p}_0(\Omega))$, for all $1\leq p,q<\infty$, along with the imposed regularity constraints on $Q$, we are prompted to revisit the weak formulations stated in \eqref{eq:weakform}. In fact, we are able to weaken the conditions on the test functions, allowing ${\varphi}\in L^2([0,T];H^1(\Omega))\cap W^{1,1}(0,T;L^2(\dom))$. Following a similar line of reasoning, we find that the prerequisites imposed on the test functions in Definition \ref{def:Q_inertia_weakQr} below can likewise be relaxed.
\end{remark}

\begin{definition}[Weak solutions of \eqref{eq:reformulation}]
\label{def:Q_inertia_weakQr}
A weak solution to \eqref{eq:reformulation} is a pair of functions $Q : [0,T]\times\Omega\to\mathbb{R}^{d\times d}$ and $r : [0,T]\times\Omega\to\mathbb{R}$ with $Q(t, x)$ trace-free and symmetric for every $(t, x)$ and satisfying
\begin{equation*}
    Q \in L^{\infty}(0, T; H^{1}_0(\dom)),\quad Q_{t} \in L^{2}([0, T]\times\Omega),\quad r\in L^{\infty}(0, T; L^{2}(\dom))
\end{equation*}
and
\begin{multline}
\label{eq:weakQ_refomulation_inertia}
    -\int_{\dom} Q(T,x):\varphi(T,:) dx +\int_{\dom} Q_0(x):\varphi(0,x) dx +\int_0^T\int_\Omega Q:\partial_t\varphi\,dx\,dt\\
  = -\sigma\int_{0}^{T}\int_{\Omega} Q_{t}:\varphi_{t}\, dx\, dt + \sigma\int_{\Omega} Q_{t}(T, x):\varphi(T, x)\, dx - \sigma\int_{\Omega} Q_{t,0}(x):\varphi(0, x)\, dx\\ 
  + \int_{0}^{T}\int_{\Omega}\Bigg(L_1 \sum_{i,j=1}^{d} \nabla Q_{ij} \cdot \nabla \varphi_{ij} \\
  + \frac{L_2+L_3}{2} \sum_{i,j,k=1}^d \left(\partial_{k}Q_{jk}\partial_i  \varphi_{ij} + \partial_{k}Q_{ik}\partial_j \varphi_{ij} -\frac{2}{d}  \partial_{i}Q_{ki} \partial_k \varphi_{jj}\right)\Bigg)dxdt
	+ \int_{0}^{T}\int_{\Omega} rP(Q):\varphi\, dx\, dt
   \end{multline}
and
\begin{equation}
\label{eq:weakr_reformulation_inertia}    \int_{0}^{T}\int_{\Omega} r \phi_{t}\, dx\, dt - \int_{\Omega} r(T, x) \phi(T, x)\, dx + \int_{\Omega} r_{0}(x)\phi(0, x)\, dx = -\int_{0}^{T}\int_{\Omega} P(Q) : Q_{t}\phi\, dx\, dt
\end{equation}
for all smooth $\varphi = (\varphi_{ij})_{ij=1}^{d} : [0,T]\times\Omega\to\mathbb{R}^{d\times d}$ and $\phi: [0, T]\times\Omega\to\mathbb{R}$ which are compactly supported within $\Omega$ for every $t\in[0,T]$.
\end{definition}

In \cite[Lemma 5.2]{GWY2020}, we have shown the equivalence of Definition \ref{def:Q_inertia_weakQ} and \ref{def:Q_inertia_weakQr} in the case $\sigma=0$. The same proof carries over for $\sigma>0$. Therefore, we will use the two definitions interchangeably in the following, depending on which one is more suitable for our needs. Finally, we define strong solutions to~\eqref{eq:qtensorflow}:

\begin{definition}[Strong solution of \eqref{eq:qtensorflow}]
\label{def:strong_solution}
A strong solution of \eqref{eq:qtensorflow} in $[0,T]\times \Omega$ for any fixed $T>0$ is a function $Q:[0,T]\times \Omega\to\mathbb{R}^{d\times d}$ that is trace-free and symmetric for every $(t,x)$. The solution satisfies Equation \eqref{eq:qtensorflow} in $[0,T]\times\Omega$ almost everywhere and fulfills the following regularity condition:

\begin{equation*}
Q\in L^\infty(0,T;H^2(\dom)),\quad Q_t\in L^\infty(0,T; H^1(\dom)),\quad Q_{tt}\in L^\infty(0,T;L^2(\dom)).
\end{equation*}
\end{definition}

In the upcoming sections, we will focus on establishing the existence and uniqueness of weak and strong solutions of~\eqref{eq:qtensorflow}. Later, we will use these insights to derive error estimates for a numerical scheme to approximate solutions of~\eqref{eq:qtensorflow}.

\section{Global Strong solution}\label{sec:strong_Solution}

Our goal is to establish the existence and uniqueness of a solution for the Q-tensor model with an inertia term, as presented in \eqref{eq:Qt}-\eqref{eq:boundary_condition}, provided that the initial and boundary conditions exhibit appropriate regularity. In particular, we will prove, as the result of a sequence of lemmas:
\begin{theorem}\label{thm:Existence_Strong_solution}
	For fixed $T>0$, the problem ~\eqref{eq:Qt} has a unique strong solution $Q$ which satisfies Definition~\ref{def:strong_solution} given $Q_0\in H^2(\Omega), Q_{0,t}\in H^1(\Omega)$. As a consequence, $Q\in L^\infty\left([0,T];L^\infty(\Omega)\right)$. In addition, this solutions is unique among all weak solutions as in Definition~\ref{def:Q_inertia_weakQ}.
\end{theorem}

In order to prove the existence of a strong solution as in Definition \ref{def:strong_solution}, it is essential to derive  higher-order regularity estimates. We will address this in the first subsection, then leverage these estimates to construct the desired solution.

\subsection{Higher-order Energy Estimates}\label{section:energy_inequality}
The primary objective of this section is to obtain higher-order energy estimates for the problem~\eqref{eq:Qt}-\eqref{eq:boundary_condition}. These estimates are crucial for establishing the existence of a strong solution. For the sake of simplicity, we will assume that $Q$ is smooth and satisfies \eqref{eq:Qt} throughout the calculations in this section. This assumption will be validated by using a Galerkin approximation scheme, which will be discussed in Section~\ref{subsection:galerkin}.

We start with some preliminary technical lemmas.
First, we state the following lemma on the Lipschitz continuity of the function $P(Q)$.
\begin{lemma}\label{lem:PLipschitz}
	There exists  a constant $D>0$ such that for any matrices $Q, \delta Q\in \mathbb{R}^{d\times d}$,
	\begin{equation*}
	|P(Q+\delta Q)-P(Q)|_F\leq D|\delta Q|_F.
	\end{equation*}
\end{lemma}
The proof of this lemma can be found in~\cite[Theorem 4.11]{GWY2020}.
In order to prove the existence of strong solutions, not only is a uniform $H^1$-bound essential, which follows from the energy inequality, but an $H^2$-bound for the tensor $Q$ also comes into play. Given the various types of second-order terms present in the equation, it is reasonable to articulate the following Lemma, which serves to simplify our analysis by focusing on the norm of $\Delta Q$.
\begin{lemma}\cite[Theorem 4.3.1.4]{NonsmoothElliptic}
	\label{lem:Q_H2_Laplacian}
	If $\partial \Omega$ is $C^2$ or polygonal~, then the following estimate holds for $Q:\dom\to \R^{d\times d}$:
	\begin{equation}
	\label{eq:H2_bound_Lap}
	\|Q\|_{H^2}\leq C_\Omega\left(\|Q\|_{L^2}+\|\Delta Q\|_{L^2}\right).
	\end{equation}
\end{lemma}
Nonetheless, it is noteworthy that securing a uniform $H^2$-bound for \eqref{eq:Qt} poses a nontrivial task, particularly when the terms following $L_2, L_3$ enter the equation. Specifically, as clarified in \cite{chen2016existence}, a fundamental maximum principle that can keep $\tr(Q^2)$ under control -- as is possible when $d=2$ -- cannot be obtained for $d=3$. As a consequence, deriving a higher-order energy inequality becomes challenging. To circumvent this obstacle, we have to make an additional assumption akin to the methodology outlined in \cite{chen2016existence}. We require that $L_1$ is significantly larger than $(L_2+L_3)$, that is:
\begin{equation}
\label{eq:L_1_large_assumption}
L_1\gg (L_2+L_3).
\end{equation}
In the following, we will always work under assumption \eqref{eq:L_1_large_assumption}.  This way, we can establish a connection between the $H^2$-norm of $Q$ and the $L^2$-norm of $H(Q)$ via the subsequent lemma, simplifying the derivation of an $H^2$-bound for $Q$ by estimating $H(Q)$ in its place.
\begin{lemma}\label{lem:Q_H_control}
	Assume assumption \eqref{eq:L_1_large_assumption} holds. Then for any fixed $T>0$, there exists $C>0$ independent of $Q$ such that
	\begin{equation}\label{eq:QH2}
	\|\Delta Q\|_{L^2}\leq C\left(z\left(\norm{Q}_{H^1}\right)+\|H(Q)\|_{L^2}\right),
	\end{equation}
 where $z(\cdot)$ is a polynomial. When we have that $\|Q\|_{H^1}$ is uniformly bounded, \eqref{eq:QH2} will reduce to
 	\begin{equation}\label{eq:QH2_easy}
	\|\Delta Q\|_{L^2}\leq C\left(1+\|H(Q)\|_{L^2}\right).
	\end{equation}
\end{lemma}
\begin{proof}
	Taking the inner product of equation~\eqref{eq:H}  with $\Delta Q$ in $L^2(\Omega)$ yields
	\begin{equation}
	\label{eq:deltaQestimate}
	L_1\|\Delta Q\|_{L^2}^2+\frac{L_2+L_3}{2}\int_\Omega \alpha(Q):\Delta Q\, dx= \int_\Omega \left(f(Q):\Delta Q+H(Q):\Delta Q\right)\, dx.
	\end{equation}
Using the definition~\eqref{eq:alpha} of $\alpha(Q)$, Lemma \ref{lem:Q_H2_Laplacian} and Assumption \eqref{eq:L_1_large_assumption}, we can control the term $\int_\Omega \alpha(Q):\Delta Q\,dx$ with the Cauchy-Schwarz inequality as
	\begin{equation*}
	\label{eq:alpha_estimate}
	\begin{aligned}
	\frac{L_2+L_3}{2}\int_\Omega \alpha(Q):\Delta Q\, dx&\leq   \frac{L_2+L_3}{2}\|\alpha(Q)\|_{L^2}\|\Delta Q\|_{L_2}\\
	&\leq C(L_2+L_3)\|Q\|_{H^2}\|\Delta Q\|_{L^2 }\\
 &\leq C(L_2+L_3)\left(\|Q\|_{L^2}+\|\Delta Q\|_{L^2}\right)\|\Delta Q\|_{L^2 } \\
 &\leq C(L_2+L_3)\left (\|Q\|_{L^2}^2+\|\Delta Q\|_{L^2 }^2\right)\\
 &\leq \frac{L_1}{8}\left (\|Q\|_{L^2}^2+\|\Delta Q\|_{L^2 }^2\right)\\
 &\leq \frac{L_1}{8}\|\Delta Q\|^2_{L^2} + C\|Q\|_{H^1}^2.
	\end{aligned}
	\end{equation*}
	The term involving $f(Q)$ on the right-hand side of ~\eqref{eq:deltaQestimate} can be estimated as
	\begin{equation*}
	\label{eq:fQ_estimate}
	\begin{aligned}
	\int_\Omega f(Q):\Delta Q\,dx&\leq \frac{L_1}{4}\|\Delta Q\|_{L^2}^2+C\|f(Q)\|_{L^2}^2\\
	&\leq \frac{L_1}{4}\|\Delta Q\|_{L^2}^2 + C\left(\|Q\|_{L^2}^2+\|Q\|_{L^4}^4+\|Q\|_{L^6}^6\right)\\
	&\leq \frac{L_1}{4}\|\Delta Q\|^2_{L^2} + C(\norm{Q}^6_{H^1}+\norm{Q}^2_{H^1}),
	\end{aligned}
	\end{equation*}
	where we have used the Sobolev embedding theorem from the second line to the third line. We can control the last term in~\eqref{eq:deltaQestimate} by Young's inequality as
	\begin{equation*}
	\label{eq:HQ_estimate}
	\int_\Omega H(Q):\Delta Q\,dx\leq \frac{L_1}{2}\|\Delta Q\|^2_{L^2}+C\|H(Q)\|_{L^2}^2
	\end{equation*}
	Combining these estimates, we obtain that
	\begin{equation*}
	\label{eq:DeltaQHQ}
	\begin{aligned}
	L_1\|\Delta Q\|_{L^2}^2 \leq \frac{7\,L_1}{8}\|\Delta Q\|_{L^2}^2+C\left(\norm{Q}^6_{H^1}+\norm{Q}^2_{H^1}\right) + C\|H(Q)\|_{L^2}^2.
	\end{aligned}
	\end{equation*}
	Then ~\eqref{eq:QH2} follows and $\eqref{eq:QH2_easy}$ is a direct result of $\eqref{eq:QH2}$ if $\|Q\|^2_{H^1}$ is bounded.
\end{proof}

\begin{remark}
	\label{rem:H_L2_Q_H2}
	Due to Lemma \ref{lem:Q_H_control} and the energy inequality, if there exists a smooth solution $Q$ satisfying \eqref{eq:Qt}, the $H^2$-norm of $Q$ is controlled by the $L^2$-norm of $H(Q)$. Conversely, the $L^2$-norm of $H(Q)$ is clearly governed by the $H^2$-norm of $Q$, a direct consequence of Sobolev's inequality. Hence, it can be asserted that the $L^2$-norm of $H(Q)$ is equivalent to the $H^2$-norm of $Q$. Subsequently, when attempting to derive an $H^2$ estimate of $Q$, we typically prefer conducting estimates on the $L^2$-norm of $H(Q)$, as it often simplifies computations.
\end{remark}

Next, we recall Agmon's inequality.
\begin{lemma}\label{lem:agmon}
Assume $s>\frac{n}{2}$ and let $s'$ be an integer such that $s'<\frac{n}{2}$. Then there exists a constant $C>0$, such that for any $f\in H^s(\Omega)$,
\begin{equation}
    \label{eq:interpolation_inequality}
    \|f\|_{L^\infty(\Omega)}\leq C \|f\|_{H^{s'}(\Omega)}^{1-t} \|f\|_{H^s(\Omega)}^t,
\end{equation}
where $t=\frac{\frac{n}{2}-s'}{s-s'}$. In particular, for the case $n=2$, $s'=0$, $s=2$, we have
\begin{equation}
    \label{eq:agmon_inequality_2d}
    \|f\|_{L^\infty(\Omega)}\leq C\|f\|_{L^2(\Omega)}^\frac{1}{2}\,\|f\|_{H^2(\Omega)}^{\frac{1}{2}},
\end{equation}
and for the case $n=3$, $s'=1$, $s=2$, we have
\begin{equation}
    \label{eq:agmon_inequality}
    \|f\|_{L^\infty(\Omega)}\leq C\|f\|_{H^1(\Omega)}^\frac{1}{2}\,\|f\|_{H^2(\Omega)}^{\frac{1}{2}}.
\end{equation}
\end{lemma}
\begin{proof}
See ~\cite[Lemma 4.10]{ConstantinFoias}.
\end{proof}

With the help of this lemma, we will prove an estimate on the growth of $H$ which will be key to proving the higher order regularity of the solution.

\begin{proposition}
\label{prop:energyinequality}
Assume $d=2,3$. Define $$W(t)=\frac{1}{2} \left\|H(t)\right\|_{L^2}^2,\quad\quad Z(t)=\frac{1}{2}\left( L_1\|\nabla Q_t(t)\|_{L^2}^2 +(L_2+L_3)\|\Div Q_t(t)\|_{L^2}^2 \right),$$ for $t\in[0,T]$. Then they satisfy the following inequality:
\begin{equation}
    \label{eq:higher_order_energy_inequality_W_Z}
     \frac{d}{dt}(W+\sigma Z)\leq \frac{C}{\sqrt{\sigma}}(1+W).
\end{equation}

\end{proposition}
\begin{proof}
By taking the derivative on both sides of ~\eqref{eq:H}, we obtain
\begin{equation*}
    \label{eq:Ht}
    H_t=L_1\Delta Q_t+\frac{L_2+L_3}{2}\alpha(Q)_t-f(Q)_t.
\end{equation*}
Using the definition of $H(Q)$ and integrating by parts, we have
\begin{equation*}
    \begin{aligned}
    \frac{dW}{dt}&=\int_\Omega  H:H_t\,dx\\
    &=\int_\Omega\left(Q_t+\sigma Q_{tt}\right):\left[L_1\Delta Q_t+\frac{L_2+L_3}{2}\alpha(Q)_t \right]\,dx-\int_\Omega  H:f(Q)_t\,dx\\
    &=-L_1\|\nabla Q_t\|_{L^2}^2-(L_2+L_3)\|\Div Q_t\|_{L^2}^2-\frac{\sigma}{2}\frac{d}{dt}\int_\Omega \left[L_1 \lvert \nabla Q_t\rvert^2+(L_2+L_3)\lvert \Div Q_t\rvert^2 \right]\,dx-\int_\Omega H:f(Q)_t\,dx\\
    &=-2Z-\sigma \frac{dZ}{dt}-\int_\Omega  H:f(Q)_t\,dx.
    \end{aligned}
\end{equation*}
It remains to control the last term $\int_\Omega H:f(Q)_t,dx$. From the energy dissipation law \eqref{eq:dissipationlaw}, we have that  $\sqrt{\sigma}\|Q_t\|_{L^2}$ is uniformly bounded. Thus, following from the definition of $f(Q)$ in \eqref{eq:f}, the Sobolev inequality, Lemma~\ref{lem:Q_H_control}, and Agmon's inequality, Lemma~\ref{lem:agmon} and the energy inequality so $Q\in H^1$ uniformly, we obtain
\begin{equation}\label{eq:f(Q)_t_estimate}
\begin{aligned}
    \|f(Q)_t\|_{L^2}
    &\leq C \left(\|Q_t\|_{L^2} + \|Q\|_{L^\infty}\|Q_t\|_{L^2}+\|Q\|^2_{L^\infty}\|Q_t\|_{L^2}\right)\\
    &\leq \frac{C}{\sqrt{\sigma}}(1+\|Q\|_{L^\infty}^2)\leq \frac{C}{\sqrt{\sigma}}(1+\|Q\|_{H^2})\leq \frac{C}{\sqrt{\sigma}} (1+\|H\|_{L^2}).
\end{aligned}
\end{equation}
Applying Cauchy-Schwarz inequality, we have
\begin{align*}
    \int_\Omega \lvert  H:f(Q)_t\rvert \,dx&\leq \frac{C}{\sqrt{\sigma}}\|H\|_{L^2}\left( 1+\|H\|_{L^2}\right)\leq \frac{C}{\sqrt{\sigma}}(1+\|H\|_{L^2}^2)\leq \frac{C}{\sqrt{\sigma}}(1+W).
\end{align*}
This leads to the desired inequality
\begin{equation*}
     \frac{d}{dt}(W+\sigma Z)=-2Z-\int_\Omega  H:f(Q)_t\,dx\leq \frac{C}{\sqrt{\sigma}}(1+W).
\end{equation*}
 
\end{proof}

Using this proposition, we can show that if the initial data possesses enough regularity, then so will the solution at a later time.

\begin{corollary}
    \label{cor:Q_H2_estimate}
   Given $Q_0\in H^2(\Omega)$ and $Q_{t}(0)\in H^1(\Omega)$, the $H^2$-norm of $Q$ is  uniformly bounded in time, that is, for some constant $C>0$,
    \begin{equation*}
        \|Q(t)\|_{H^2}\leq C(\sigma), \quad \|Q_t(t)\|_{H^1}\leq C(\sigma)\quad \forall t\in [0,T].
    \end{equation*}
\end{corollary}
\begin{proof}

If $Q_0\in H^2(\Omega)$ and $Q_t(0)=Q_{t,0}\in H^1(\Omega)$, then
\begin{equation*}
    W(0)=\frac{1}{2} \|H(0)\|^2_{L^2}\leq C\left(\|Q_0\|_{H^2}^2+\|Q_0\|_{L^2}^2+\|Q_0\|_{L^6}^6\right)\leq C,
\end{equation*}
and
\begin{equation*}
    Z(0)=\frac{1}{2}\left( L_1\|\nabla Q_t(0)\|_{L^2}^2 +(L_2+L_3)\|\Div Q_t(0)\|_{L^2}^2 \right)\leq C.
\end{equation*}
From Proposition \ref{prop:energyinequality}, we have
\begin{equation*}
    \frac{d}{dt}(W+\sigma Z)\leq \frac{C}{\sqrt{\sigma}}(1+W).
\end{equation*}
Then Gr\"{o}nwall's inequality implies that
\begin{equation*}
    W(t)+\sigma Z(t)\leq C(T,\sigma) \left( W(0)+\sigma Z(0)  \right)\leq C(\sigma),\quad t\in [0,T].
\end{equation*}
Uniform boundedness of $Q$ in the $H^2$-norm then follows from Lemma \ref{lem:Q_H_control} and Remark~\ref{rem:H_L2_Q_H2}.
\end{proof}

Moreover, in situations where the inertia effect cannot be neglected, which corresponds to $\sigma > 0$, we can also obtain an energy estimate for $Q_{tt}$, as presented in the following corollary.

\begin{corollary}
    \label{cor:Q_tt_estimate}
    There is a constant $C(\sigma)>0$ such that the second derivative with respect to time of $Q$ satisfies
\begin{equation*}
    \|Q_{tt}(t)\|_{L^2}\leq C(\sigma)\quad\forall t\in [0,T], \quad \text{and}\quad \|Q_{tt}\|_{L^2([0,T]\times\Omega)}\leq C(\sigma).
\end{equation*}
\end{corollary}

\begin{proof}
We differentiate \eqref{eq:Qt} and obtain
\begin{align*}
    Q_{tt}&=-\sigma Q_{ttt}+ H_t\\
    &=-\sigma Q_{ttt}+ \left(L_1\Delta Q_t+\frac{L_2+L_3}{2}\alpha(Q)_t-f(Q)_t\right)
\end{align*}

Taking inner product with $-Q_{tt}$ on both sides and integrating and integrating by parts, we get
\begin{equation*}
    -\int_\Omega \lvert Q_{tt}\rvert^2\,dx=\frac{\sigma}{2}\frac{d}{dt}\int_\Omega \lvert Q_{tt}\rvert^2_F\,dx+\frac{1}{2}\frac{d}{dt}\int_\Omega \left[L_1\lvert \nabla Q_t\rvert^2 + (L_2+L_3)\lvert \Div Q_t\rvert^2\right]\,dx+ \int_\Omega f(Q)_t:Q_{tt}\,dx
\end{equation*}
Combining \eqref{eq:f(Q)_t_estimate} with Corollary~\ref{cor:Q_H2_estimate}, we obtain
\begin{equation*}
    \|f(Q)_t\|_{L^2}\leq \frac{C}{\sqrt{\sigma}}(1+\|H\|_{L^2})\leq C(\sigma).
\end{equation*}
Then we can use the Cauchy-Schwarz inequality to get
\begin{align*}
    &\frac{1}{2}\frac{d}{dt}\int_\Omega \left[\sigma\lvert Q_{tt}\rvert^2+ L_1\lvert \nabla Q_t\rvert^2+ (L_2+L_3)\lvert \Div Q_t\rvert^2\right ]\,dx\\&=-\int_\Omega \lvert Q_{tt}\rvert^2\,dx- \int_\Omega f(Q)_t:Q_{tt}\,dx\\
    &\leq -\int_\Omega \lvert Q_{tt}\rvert^2\,dx+\frac{1}{2}\int_\Omega \lvert Q_{tt}\rvert^2\,dx+ C\|f(Q)_t\|^2_{L^2}\\
    &\leq -\frac{1}{2}\int_\Omega \lvert Q_{tt}\rvert^2\,dx+C(\sigma)
\end{align*}
Integrating both sides with respect to time yields
\begin{align*}
    &\quad\,\,\sigma \|Q_{tt}(t)\|_{L^2}^2+ L_1\|\nabla Q_t(t)\|_{L^2}^2+ (L_2+L_3)\|\Div Q_t(t)\|_{L^2}^2+\|Q_{tt}\|_{L^2([0,T]\times\Omega)}^2\\
    &\leq C(\sigma)T+\sigma \|Q_{tt}(0)\|^2_{L^2}+ L_1\|\nabla Q_t(0)\|_{L^2}^2+ (L_2+L_3)\|\Div Q_t(0)\|_{L^2}^2.
\end{align*}
We note that $\sigma Q_{tt}(0)=-Q_t(0)+ H(0)$ implies that $\sigma Q_{tt}(0)\in L^2(\dom)$ and so we can conclude that
\begin{equation*}
    \sigma \|Q_{tt}(t)\|_{L^2}^2+ L_1\|\nabla Q_t(t)\|_{L^2}^2+ (L_2+L_3)\|\Div Q_t(t)\|_{L^2}^2+\|Q_{tt}\|_{L^2([0,T]\times\Omega)}^2\leq C\left({\sigma}\right).
\end{equation*}
\end{proof}

\subsection{Galerkin Approximation}\label{subsection:galerkin}
Finally, to prove the existence of strong solutions, we will use a standard Galerkin approximation scheme for~\eqref{eq:Qt}, as was done, for example, in~\cite{Lin1995NonparabolicDS}. Our treatment focuses on the Dirichlet boundary condition case; Neumann boundary conditions can be handled similarly.

To establish the Galerkin approximation scheme, we first introduce a sequence of matrix functions $\{\omega_s\}_{s=1}^{\infty}$, where each $\omega_s$ is a symmetric, trace-free eigenfunction of the operator $L=-\Delta$, subject to the Dirichlet boundary conditions. Specifically, we have

\begin{equation}\label{eq:Laplace_eigen}
-\Delta \omega_s = \lambda_s\omega_s,
\end{equation}
where $\{\lambda_s\}_{s=1}^{\infty}$ is a sequence of eigenvalues satisfying $0\leq\lambda_1\leq\lambda_2\leq\cdots\leq \lambda_n\leq\cdots$ with $\lim\limits_{n\to\infty}\lambda_n=\infty$. This choice of $\{\omega_s\}_{s=1}^\infty$ constitutes a family of smooth functions that forms an orthogonal basis for the spaces $L^2(\Omega)$ and $H^1_0(\Omega)$~\cite[Section 6.5.1]{Evans}.

Given a positive integer $N$, we denote by $U_N=\text{span}\{\omega_1,\omega_2,\cdots,\omega_N\}$ the finite-dimensional subspace generated by the first $N$ eigenfunctions. We also introduce $P_N:H_0^1(\Omega)\cap H^r(\Omega)\to U_N$, the orthogonal projection operator from $H^1_0(\Omega)$ onto the finite dimensional subspace $U_N$ (the choice or $r$ here depends on the requirement of the regularity we need, for example, in this proof, we can choose $r=3$).

Now we construct a smooth function $Q^N: [0,T]\times \Omega\to \R^{d\times d}$ that belongs to $U_N$ at every point in time and approximates the solution to our original problem. The Galerkin approximation problem is thus formulated as follows:

\begin{subequations}\label{eq:Galerkin_scheme}
\begin{empheq}[left=\empheqlbrace]{align}
 &Q^N_t=-\sigma Q^N_{tt}+ \left(L_1\Delta Q^N+\frac{L_2+L_3}{2}\alpha(Q^N)-f(Q^N)\right):= -\sigma Q^N_{tt}+ H(Q^N),\label{eq:Q_Galerkin_formula}\\
   & Q^N(0,x) = P_N\,Q_0(x),\\ &Q^N_t(0,x) = P_N\, Q_{t,0}(x),\label{eq:Q0_Garlerkin}
\end{empheq}
\end{subequations}
where $\alpha(Q^N)$ and $f(Q^N)$ are given by
\begin{subequations}
    \label{eq:alpha_N_Galerkin}
    \begin{align}
         &\alpha(Q^N)_{ij}= \sum_{k=1}^d\left( \partial_{ik}Q^N_{jk} + \partial_{jk}Q^N_{ik}\right) - \frac{2}{d} \sum_{k,s=1}^d \partial_{ks}Q^N_{ks} \delta_{ij}\\
 &f(Q^N)=a Q^N-b\left[(Q^N)^2-\frac{1}{d}\text{tr}\left((Q^N)^2\right)I\right]+c\text{tr}\left((Q^N)^2\right)Q^N.
    \end{align}
\end{subequations}

Next we prove the existence and uniqueness of a solution of~\eqref{eq:Galerkin_scheme}.

 \begin{lemma}
 \label{lem:existence_approximate_solution}
 For any $N>0$, there exists $T_N>0$ such that the problem ~\eqref{eq:Q_Galerkin_formula} -~\eqref{eq:Q0_Garlerkin} admits a unique smooth solution in $\Omega_{T_N}=[0,T_N)\times\Omega$.
 \end{lemma}
 \begin{proof}
Since we seek $Q^N\in U_N$,  we  make the ansatz
\begin{equation*}\label{eq:Q_Galerkin_approximation}
    Q^N=\sum_{s=1}^N\, g^N_s(t)\,\omega_s(x)
\end{equation*}
for suitable functions $\{g_s(t)\}_{s=1}^n$ such that~\eqref{eq:Galerkin_scheme} is satisfied. In fact, by testing \eqref{eq:Q_Galerkin_formula} with $\omega_s$ for $s=1,2,\cdots,N$, and using the orthogonality of the $\omega_s$ in $L^2$ and $H^1_0$, we find that \eqref{eq:Q_Galerkin_formula} can be rewritten as an equivalent system of ordinary differential equations:

 \begin{equation*}
     \label{eq:g_i_formula}
     \begin{aligned}
     \frac{d}{dt}\,g^N_s(t)=&-\sigma \frac{d^2}{dt^2}\, g_s^N(t)   -L_1\lambda_s\, g^N_{s}(t)-({L_2+L_3})\sum_{m=1}^NA_{m,s}\,g^N_m (t) \\
     &\phantom{-\sigma \frac{d^2}{dt^2}\, g_s^N(t) +  -} - ag^N_s(t) + b  \sum_{k,m=1}^N B_{k,m,s}\,g^N_k(t)\,g^N_{m}(t) - c \sum_{k,m,n=1}^N C_{k,m,n,s}\,g^N_k(t)\,g^N_m(t)\,g_n(t) ,
     \end{aligned}
 \end{equation*}
where
\begin{equation*}\label{eq:cofficients_for_g}
 \begin{aligned}
  A_{m,s}=\int_\Omega\Div \omega_m\cdot\Div\omega_s,\quad B_{k,m,s}=\int_\Omega \tr(\omega_k\omega_m\omega_s),\quad C_{k,m,n,s}=\int_\Omega \tr(\omega_k\omega_m)\tr(\omega_n\omega_s).
 \end{aligned}
 \end{equation*}
The initial conditions~\eqref{eq:Q0_Garlerkin} are
 \begin{equation*}
     \label{eq:g_i0}
     g^N_s(0)=\langle Q_0, \omega_s \rangle,\quad \frac{dg_s^N}{dt}(0)=\langle Q_{t,0}, \omega_s \rangle,
 \end{equation*}
 for $s=1,2,\cdots,N$. Now, the standard existence theory of ordinary differential equations ensures the existence of a unique solution $g^N(t)=\left(g^N_1(t), g^N_2(t),\cdots,g^N_N(t)\right)$ satisfying the corresponding initial conditions on the interval $[0,T_N)$, for some $T_N>0$ (for example, see \cite[Theorem 2.2]{teschl2012ordinary}). A bootstrap argument confirms the smoothness of $g^N(t)$ and, hence, of $Q^N$, thereby completing the proof of Lemma~\ref{lem:existence_approximate_solution}.
\end{proof}

\begin{remark}
    It is straightforward to prove the energy dissipation law~\eqref{eq:dissipationlaw}, Proposition \ref{prop:energyinequality}, and Corollary~\ref{cor:Q_tt_estimate} for the approximate solution $Q^N$. Consequently, for fixed $N$ and any $t\in[0,t_N]$, we have
\begin{equation}
\label{eq:Q^N_estimates_uniform}
\|Q^N(t)\|_{H^2\cap H^1_0}\leq C,   \quad \|Q^N_t(t)\|_{H^1}
\leq C, \quad \|Q^N_{tt}(t)\|_{L^2}\leq C,
\end{equation}
where $C$ is independent of $N$ and $t_N$ since $\|P_NQ_0\|_{H^r} \le \|Q_0\|_{H^r}$ for $r\geq 0$.
These bounds imply that $Q^N$ does not blow up in finite time and, hence, that the solution $Q^N$ can be extended to $[0,T]$ for arbitrarily large $T>0$. 
\end{remark}

\subsection{Existence and uniqueness results}
We have thus far constructed a sequence of approximating problems and have proven their well-posedness. This results in a sequence of approximating solutions $\{Q^N\}_{N\in\mathbb{N}}$, each solving system \eqref{eq:Galerkin_scheme} for $N\in\N$. Finally, we show that the limit of these approximate solutions is a strong solution that satisfies Definition \ref{def:strong_solution} for the original system \eqref{eq:Qt} and is indeed unique. 

\begin{proof}[Proof of Theorem~\ref{thm:Existence_Strong_solution}]

We infer from~\eqref{eq:Q^N_estimates_uniform} that 
\begin{equation}
    \label{eq:Q^N_estimates}
    Q^N\in L^\infty\left(0,T;H^2(\Omega)\cap H^1_0(\Omega)\right),\quad Q^N_t\in L^\infty(0,T;H^1(\Omega)), \quad Q_{tt}^N\in L^\infty(0,T;L^2(\Omega)).
\end{equation}
Using the Banach-Alaoglu theorem and the Aubin-Lions lemma~\cite{Simon}, we can pass the limit, up to a subsequence, which we still denote by $\{Q^N\}_{N\in\N}$, and conclude that there exists 
\begin{equation*}
\label{eq:existence_strong_solution}
    Q\in L^\infty\left(0,T;H^2(\Omega)\cap H^1_0(\Omega)\right)\,
\end{equation*}
with its weak derivatives satisfying
\begin{equation*}
    Q_t\in L^\infty\left(0,T;H^1(\Omega)\right),\quad Q_{tt}\in L^\infty\left(0,T;L^2(\Omega)\right),
\end{equation*}
such that
\begin{equation}
    \label{eq:passing_the_limit}
  \begin{aligned}
    Q^N\overset{\ast} {\rightharpoonup} Q \quad\text{in}\,\,L^\infty(0,T;H^2(\Omega)\cap H^1_0(\Omega)),\quad\quad & Q^N\to Q \quad\text{in}\,\, C([0,T];H^1(\Omega)),
    \\ 
    Q^N_t\to Q_t\quad\text{in}\,\,C([0,T];L^2(\Omega)),\quad\quad &Q^N_{tt}\overset{\ast} {\rightharpoonup} Q_{tt}\quad\text{in}\,\,L^\infty(0,T;L^2(\Omega)).
  \end{aligned}
\end{equation}
It remains to show that $Q$ solves~\eqref{eq:Qt} and satisfies initial condition~\eqref{eq:initial_condition}. Note that the left-hand side of~\eqref{eq:Q_Galerkin_formula} converges to $Q_t$ in $L^2(0,T;L^2(\Omega))$ in the norm sense. Using the strong convergence of $Q^N$ in $C(0,T;H^1(\Omega))$ and the weak star convergence of $Q^N_{tt}$ and $\Delta Q^N$ in $L^\infty(0,T;L^2(\Omega))$, we obtain in the limit $N\to\infty$,
\begin{equation*}
    \label{eq:Qt_limit}
    Q_t = -\sigma Q_{tt}+ H(Q)\quad \text{in }\, L^2\left(0,T;L^2(\Omega)\right).
\end{equation*}
 Finally, thanks to the uniform time continuity of the $Q^N$, we have $Q^N(0)\to Q(0)$ in $H^1(\Omega)$ and $Q^N(0)\to Q_0$ in $H^1(\Omega)$ and hence $Q(0)=Q_0$. By passing to the limit in $C([0,T];L^2(\Omega))$, $Q^N_t(0)\to Q_{t}(0)$ in $L^2(\Omega)$ and $Q^N_t(0)\to Q_{t,0}$ we obtain that $Q_t(0)=Q_{t,0}$. The uniform $L^\infty$ estimate for $Q$ is a direct consequence of Lemma \ref{lem:agmon}. This proves the existence of a strong solution.
 Next, we show that this solution is unique among all weak solutions as in Definition~\ref{def:Q_inertia_weakQr} (or Definition~\ref{def:Q_inertia_weakQ}). To do so, we use the relative entropy method, which was introduced by Dafermos~\cite{Dafermos1979,Dafermos1979a,Dafermos2016book} and DiPerna~\cite{Diperna1979} in the context of hyperbolic conservation laws.
 
 Let us denote $(Q,r)$ as a weak solution satisfying Definition \ref{def:Q_inertia_weakQr} and $\hat{Q}$ as a strong solution fulfilling Definition \ref{def:strong_solution}. As explained earlier, by setting $\hat{r}=r(\hat{Q})$, we obtain a strong solution for the reformulated system \eqref{eq:reformulation}. We assume that both $(Q,r)$ and $(\hat{Q},\hat{r})$ have the same initial data.

Then, based on the definition of energy in \eqref{eq:dissipationlawr}, given the weak solution $(Q,r)$ and the strong solution $(\hat{Q},\hat{r})$, we define the relative energy as
\begin{equation}\label{eq:relative_energy_definition}
\begin{aligned}
    E_{rel}(t)= \frac{ L_1}{2}\|\nabla Q(t,\cdot)-\nabla \hat{Q}(t,\cdot)\|^2_{L^2}&+ \frac{L_2+L_3}{2}\|\Div Q(t,\cdot)-\Div \hat{Q}(t,\cdot)\|_{L^2}^2\\&+\frac{1}{2}\|r(t,\cdot)-\hat{r}(t,\cdot)\|^2_{L^2}+\frac{\sigma}{2}\|Q_t(t,\cdot)-\hat{Q}_t(t,\cdot)\|^2_{L^2}.
    \end{aligned}
\end{equation}
Expanding this expression and using the energy inequality we obtain
\begin{align*}
    E_{rel}(t)&=\frac{ L_1}{2}\|\nabla Q(t,\cdot)\|_{L^2}^2+ \frac{L_2+L_3}{2}\|\Div Q(t,\cdot)\|^2_{L^2}+\frac{1}{2}\|r(t,\cdot)\|_{L^2}^2+\frac{\sigma}{2}\|Q_t(t,\cdot)\|^2_{L^2}\\
    &\qquad+\frac{ L_1}{2}\|\nabla \hat{Q}(t,\cdot)\|_{L^2}^2+ \frac{L_2+L_3}{2}\|\Div \hat{Q}(t,\cdot)\|^2_{L^2}+\frac{1}{2}\|\hat{r}(t,\cdot)\|_{L^2}^2+\frac{\sigma}{2}\|\hat{Q}_t(t,\cdot)\|^2_{L^2}\\
    &\qquad- L_1\int_\Omega \nabla Q(t,x)\cdot\nabla \hat{Q}(t,x)\,dx- (L_2+L_3)\int_\Omega \Div Q(t,x)\cdot\Div \hat{Q}(t,x)\,dx\\
    &\qquad- \int_\Omega r(t,x)\,\hat{r}(t,x)\,dx-\sigma\int_\Omega Q_t(t,x):\hat{Q}_t(t,x)\,dx\\
    &\leq \frac{ L_1}{2}\|\nabla Q(0,\cdot)\|_{L^2}^2+ \frac{L_2+L_3}{2}\|\Div Q(0,\cdot)\|^2_{L^2}+\frac{1}{2}\|r(0,\cdot)\|_{L^2}^2+\frac{\sigma}{2}\|Q_t(0,\cdot)\|^2_{L^2}\\
    &\qquad+ \frac{L_1}{2}\|\nabla \hat{Q}(0,\cdot)\|_{L^2}^2+ \frac{L_2+L_3}{2}\|\Div \hat{Q}(0,\cdot)\|^2_{L^2}+\frac{1}{2}\|\hat{r}(0,\cdot)\|_{L^2}^2+\frac{\sigma}{2}\|\hat{Q}_t(0,\cdot)\|^2_{L^2}\\
    &\qquad-\int_0^t\|Q_t(\tau,\cdot)\|_{L^2}^2\,d\tau-\int_0^t\|\hat{Q}_t(\tau,\cdot)\|_{L^2}^2\,d\tau\\
    &\qquad- L_1\int_\Omega \nabla Q(t,x)\cdot\nabla \hat{Q}(t,x)\,dx- (L_2+L_3)\int_\Omega \Div Q(t,x)\cdot\Div \hat{Q}(t,x)dx\\
    &\qquad- \int_\Omega r(t,x)\,\hat{r}(t,x)\,dx-\sigma\int_\Omega Q_t(t,x):\hat{Q}_t(t,x)\,dx\\
     &={ L_1}\|\nabla \hat{Q}(0,\cdot)\|_{L^2}^2+ (L_2+L_3)\|\Div \hat{Q}(0,\cdot)\|^2_{L^2}+ \|\hat{r}(0,\cdot)\|_{L^2}^2+\sigma\|\hat{Q}_t(0,\cdot)\|^2_{L^2}\\
      &\qquad-\int_0^t\|Q_t(\tau,\cdot)\|_{L^2}^2 d\tau-\int_0^t\|\hat{Q}_t(\tau,\cdot)\|_{L^2}^2\,d\tau\\
     &\qquad- L_1\int_\Omega \nabla Q(t,x)\cdot\nabla \hat{Q}(t,x)\,dx- (L_2+L_3)\int_\Omega \Div Q(t,x)\cdot\Div \hat{Q}(t,x)dx\\
    &\qquad- \int_\Omega r(t,x)\,\hat{r}(t,x)\,dx-\sigma\int_\Omega Q_t(t,x):\hat{Q}_t(t,x)\,dx
\end{align*}
As it is stated in Remark \ref{rem:L2_test_functions}, we can choose the test function in the weak formulation \eqref{eq:weakform} to be $\hat{Q}_t$ given the regularity satisfied by the strong solution. We obtain 
\begin{equation}\label{eq:drittseck}
		\begin{split}
		&\quad\,\,-\int_0^t\int_\Omega Q_t(\tau,x):\hat{Q}_t(\tau,x)\,dx\,d\tau\\
   &= -\sigma\int_{0}^{t}\int_{\Omega} Q_{t}(\tau,x):\hat{Q}_{tt}(\tau,x)\, dx\, d\tau + \sigma\int_{\Omega} Q_{t}(t, x):\hat{Q}_t(t, x)\, dx - \sigma\int_{\Omega} Q_{t}(0, x):\hat{Q}_t(0, x)\, dx\\ 
		& \qquad+ L_1\int_0^t\int_\Omega \nabla Q(\tau,x)\cdot\nabla\hat{Q}_t(\tau,x)\,dx\,d\tau+{ (L_2+L_3)} \int_0^t\int_\Omega \Div Q(\tau,x)\cdot\Div \hat{Q}_t(\tau,x)\,dx\,d\tau \\
		&\qquad+ \int_{0}^{t}\int_{\Omega} r(\tau,x)P\left(Q(\tau,x)\right):\hat{Q}_t(\tau,x)\, dx\, d\tau.
		\end{split}
		\end{equation}
We proceed now to demonstrate an alternative method to compute the integral $\int_0^t\int_\Omega Q_t(\tau,x):\hat{Q}_t(\tau,x)dxd\tau$. Prior to delving into this, let us first establish two identities that will facilitate our subsequent calculations. Suppose $\varphi$ is a smooth function with compact support. Using integration by parts and the weak formulation for $Q_t$, yields the following relationships:
\begin{equation*}
    \begin{aligned}
       &\quad\,\, \int_0^t\int_\Omega Q_t(\tau,x):\Delta\varphi(\tau,x)\,dx\,d\tau-\int_0^t\int_\Omega \nabla Q(\tau,x):\nabla\varphi_t(\tau,x)\,dx\,d\tau\\
       &=\int_0^t\int_\Omega Q_t(\tau,x):\Delta\varphi(\tau,x)\,dx\,d\tau+\int_0^t\int_\Omega Q(\tau,x):\Delta\varphi_t(\tau,x)\,dx\,d\tau\\
       &=\int_\Omega Q(t,x):\Delta \varphi(t,x)\,dx-\int_\Omega Q(0,x):\Delta \varphi(0,x)\,dx\\
       &=-\int_\Omega \nabla Q(t,x)\cdot\nabla\varphi(t,x)\,dx+\int_\Omega \nabla Q(0,x)\cdot\nabla\varphi(0,x)\,dx,
    \end{aligned}
\end{equation*}
and 
\begin{equation*}
    \begin{aligned}
        &\quad\,\, \int_0^t\int_\Omega Q_t(\tau,x):\alpha \left( \varphi(\tau,x) \right)\,dx\,d\tau-2\int_0^t\int_\Omega \Div Q(\tau,x)\cdot\Div \varphi_t(\tau,x)\,dx\,d\tau\\
        &= \int_0^t\int_\Omega Q_t(\tau,x):\alpha \left( \varphi(\tau,x) \right)\,dx\,d\tau+\int_0^t\int_\Omega Q(\tau,x):\alpha \left(\varphi_t(\tau,x)\right)\,dx\,d\tau\\
        &=\int_\Omega Q(t,x):\alpha \left( \varphi (t,x) \right)\,dx-\int_\Omega Q(0,x):\alpha \left( \varphi (0,x) \right)\,dx\\
        &=-2\int_\Omega \Div Q(t,x)\cdot \Div \varphi(t,x)\,dx+2\int_\Omega \Div Q(0,x)\cdot \Div \varphi(0,x)\,dx.
    \end{aligned}
\end{equation*}
These two identities retain their validity even if $\varphi$ is replaced by the strong solution $\hat{Q}$, as per the standard density argument and the regularity condition. Upon substituting $\hat{Q}_t$ with $\left(-\sigma \hat{Q}_{tt}+ L_1\Delta \hat{Q}+ \frac{L_2+L_3}{2}\alpha(\hat{Q})- \hat{r}P(\hat{Q})\right)$ in the equation $\int_0^t\int_\Omega Q_t(\tau,x):\hat{Q}_t(\tau,x)dxd\tau$ and invoking the identities derived above, we obtain
\begin{align}\label{eq:drittseck2}
    \begin{split}
        &-\int_0^t\int_\Omega Q_t(\tau,x):\hat{Q}_t(\tau,x)\,dx\,d\tau\\
        &=-\int_0^t\int_\Omega Q_t(\tau,x):\left(-\sigma \hat{Q}_{tt}(\tau,x)+ L_1\Delta \hat{Q}(\tau,x)+ \frac{L_2+L_3}{2}\alpha\left(\hat{Q}(\tau,x)\right)- \hat{r}(\tau,x)P\left(\hat{Q}(\tau,x)\right)\right)\,dx\,d\tau\\
        &=\sigma\int_0^t\int_\Omega Q_{t}(\tau,x):\hat{Q}_{tt}(\tau,x)\, dx\, d\tau- L_1\int_0^t\int_\Omega \nabla Q(\tau,x):\nabla\hat{Q}_t(\tau,x)\,dx\,d\tau\\
        &\qquad-{(L_2+L_3)} \int_0^t\int_\Omega \Div Q(\tau,x)\cdot\Div \hat{Q}_t(\tau,x)\,dx\,d\tau +  L_1\int_\Omega\nabla Q(t,x): \nabla\hat{Q}(t,x)\,dx\\
        &\qquad - L_1\int_\Omega \nabla Q(0,x):\nabla \hat{Q}(0,x)\,dx+  (L_2+L_3)\int_\Omega \Div Q(t,x)\cdot\Div \hat{Q}(t,x)\,dx\\
        &\qquad  -  (L_2+L_3)\int_\Omega \Div Q(0,x)\cdot\Div \hat{Q}(0,x)\,dx + \int_0^t\int_\Omega \hat{r}(\tau,x)P\left( \hat{Q}(\tau,x)\right):{Q}_t(\tau,x)\,dx\,d\tau,
    \end{split}
\end{align}
where we have applied integration by parts.
Adding~\eqref{eq:drittseck} and~\eqref{eq:drittseck2} and rearranging terms, we have
\begin{multline}\label{eq:faen}
-\sigma\int_{\dom}Q_t(t,x):\hat{Q}_t(t,x) dx -L_1 \int_{\dom}\Grad Q(t,x):\Grad \hat{Q}(t,x) dx - (L_2+L_3)\int_{\dom} \Div Q(t,x)\cdot \Div \hat{Q}(t,x) dx \\
= 2 \int_0^t \int_{\dom} Q_t (\tau,x):\hat{Q}(\tau,x)dx d\tau - \sigma\int_{\dom} Q_t(0,x):\hat{Q}_t(0,x)dx -L_1 \int_{\dom}\Grad Q(0,x):\Grad \hat{Q}(0,x) dx\\
-(L_2+L_3)\int_{\dom}\Div Q(0,x)\cdot \Div \hat{Q}(0,x) dx\\ +\int_0^t\int_{\dom}\left(r(\tau,x)P(Q(\tau,x)):\hat{Q}_t(\tau,x)+\hat{r}(\tau,x)P(\hat{Q}(\tau,x)):Q_t(\tau,x)\right) dx d\tau.
\end{multline}
Taking $\phi=\hat{r}$ as a test function in~\eqref{eq:weakr_reformulation_inertia}, and using that $\hat{r}$ is a strong solution, we also have
\begin{multline}\label{eq:faen2}
-\int_{\dom} r(t,x)\hat{r}(t,x) dx = -\int_0^t\int_{\dom} r(\tau,x)P(\hat{Q}(\tau,x)):\hat{Q}_t(\tau,x) dx d\tau\\
 -\int_{\dom} r_0(x)^2 dx - \int_0^t\int_{\dom} P(Q(\tau,x)):Q_t(\tau,x)\hat{r}(\tau,x)dx d\tau.
\end{multline}
Using~\eqref{eq:faen} and~\eqref{eq:faen2}, we can estimate the relative entropy as follows:
\begin{equation*}
    \begin{aligned}
        E_{rel}(t)&\leq-\int_0^t\|Q_t(\tau,\cdot)\|_{L^2}^2\,d\tau-\int_0^t\|\hat{Q}_t(\tau,\cdot)\|_{L^2}^2\,d\tau+2\int_0^t\int_\Omega Q_t(\tau,x):\hat{Q}_t(\tau,x)\,dx\,d\tau\\
     &\quad+ \int_0^t\int_\Omega r(\tau,x)P\left(Q(\tau,x)\right):\hat{Q}_t(\tau,x)\,dx\,d\tau- \int_0^t\int_\Omega r(\tau,x)P\left(\hat{Q}(\tau,x)\right):\hat{Q}_t(\tau,x)\,dx\,d\tau\\
    &\quad+ \int_0^t\int_\Omega \hat{r}(\tau,x)P\left(\hat{Q}(\tau,x)\right):Q_t(\tau,x)\,dx\,d\tau- \int_0^t\int_\Omega \hat{r}(\tau,x)P\left(Q(\tau,x)\right):Q_t(\tau,x)\,dx\,d\tau\\
    &=-\int_0^t\|Q_t(\tau,\cdot)-\hat{Q}_t(\tau,\cdot)\|_{L^2}^2\,d\tau\\
    &\qquad+ {\int_0^t\int_\Omega\left[\hat{r}(\tau,x)Q_t(\tau,x)-r(\tau,x)\hat{Q}_t(\tau,x)\right]:\left[P\left(\hat{Q}(\tau,x)\right)-P\left(Q(\tau,x)\right)\right]}\,dx\,d\tau.
    \end{aligned}
\end{equation*}
The remaining task is to find an appropriate estimate for the last term above. For simplicity of notation, we use $J$ to denote ${\int_0^t\int_\Omega\left[\hat{r}(\tau,x)Q_t(\tau,x)-r(\tau,x)\hat{Q}_t(\tau,x)\right]:\left[P\left(\hat{Q}(\tau,x)\right)-P\left(Q(\tau,x)\right)\right]}\,dx\,d\tau$. To obtain control of $J$, we first provide an estimate concerning $\int_0^t\|P\big(\hat{Q}(\tau,\cdot)\big)-P\big(Q(\tau,\cdot)\big)\|_{L^2}^2\,d\tau$. Recalling the Lipschitz continuity of $P(Q)$ stated in Lemma \ref{lem:PLipschitz}, we have
\begin{equation*}\label{eq:PQL2}
\begin{aligned}
 \vert P\big(\hat{Q}(\tau,x)\big)-P\big(Q(\tau,x)\big) \vert_F&\leq C \vert \hat{Q}(\tau,x)-Q(\tau,x)\vert_F,
\end{aligned}
\end{equation*}
for some constant $C>0$ for all $(\tau,x)\in[0,T]\times\dom$. Furthermore, we can rewrite $J$ as
\begin{equation*}
    \begin{aligned}
    J&=\underbrace{\int_0^t\int_\Omega\hat{r}(\tau,x)\left[P\left(\hat{Q}(\tau,x)\right)-P\left(Q(\tau,x)\right)\right]:\left[Q_t(\tau,x)-\hat{Q}_t(\tau,x)\right]\,dx\,d\tau}_{K_1}\\
    &\quad+\underbrace{\int_0^t\int_\Omega\left[\hat{r}(\tau,x)-{r}(\tau,x)\right]\left[P\left(\hat{Q}(\tau,x)\right)-P\left(Q(\tau,x)\right)\right]:\hat{Q}_t(\tau,x)\,dx\,d\tau}_{K_2}:=K_1+K_2.
    \end{aligned}
\end{equation*}
For $K_1$, by using the fact that $\hat{r}\in H^1(\Omega)$, H\"older's inequality, the Lipschitz continuity of $P$, Sobolev inequality, Poincare inequality, and Young's inequality, we get:
\begin{equation*}
    \begin{aligned}
        &\quad\,\,\lvert K_1\rvert\\&\leq \left(\int_0^t\int_\Omega\vert Q_t(\tau,x)-\hat{Q}_t(\tau,x)\vert_F^2\,dx\,d\tau\right)^{\frac{1}{2}} \left(\int_0^t\int_\Omega\vert P\big(\hat{Q}(\tau,x)\big)-P\big(Q(\tau,x)\big)\vert_F^4\,dx\,d\tau\right)^{\frac{1}{4}}\left(\int_0^t\int_\Omega \lvert\hat{r}(\tau,x)\rvert^4\,dx\,d\tau\right)^{\frac{1}{4}}\\
        &\leq C\left(\int_0^t\int_\Omega\vert Q_t(\tau,x)-\hat{Q}_t(\tau,x)\vert_F^2\,dx\,d\tau\right)^{\frac{1}{2}} \left(\int_0^t\int_\Omega\vert \hat{Q}(\tau,x)-Q(\tau,x)\vert_F^4\,dx\,d\tau\right)^{\frac{1}{4}}\left(\int_0^t\int_\Omega \lvert\hat{r}(\tau,x)\rvert^4\,dx\,d\tau\right)^{\frac{1}{4}}\\
        &\leq C\left(\int_0^t\int_\Omega\vert Q_t(\tau,x)-\hat{Q}_t(\tau,x)\vert_F^2\,dx\,d\tau\right)^{\frac{1}{2}} \left(\int_0^t\int_\Omega\vert \nabla\hat{Q}(\tau,x)-\nabla Q(\tau,x)\vert_F^2\,dx\,d\tau\right)^{\frac{1}{2}}\left(\int_0^t\int_\Omega \lvert\nabla\hat{r}(\tau,x)\rvert^2\,dx\,d\tau\right)^{\frac{1}{2}}\\
        &\leq \int_0^t\int_\Omega \vert Q_t(\tau,x)-\hat{Q}_t(\tau,x)\vert_F^2\,dx\,d\tau+C\int_0^t\int_\Omega \vert \nabla\hat{Q}(\tau,x)-\nabla Q(\tau,x)\vert_F^2\,dx\,d\tau,
    \end{aligned}
\end{equation*}
where $C$ is a constant depending on $\|\hat{r}\|_{L^2(0,T;H^1)}$ and the Lipschitz constant $D$. For $K_2$, we have
\begin{equation*}
    \begin{aligned}
    &\quad\,\,|K_2|\\&\leq \left(\int_0^t\int_\Omega\vert r(\tau,x)-\hat{r}(\tau,x)\vert_F^2\,dx\,d\tau\right)^{\frac{1}{2}} \left(\int_0^t\int_\Omega\vert P\big(\hat{Q}(\tau,x)\big)-P\big(Q(\tau,x)\big)\vert_F^4\,dx\,d\tau\right)^{\frac{1}{4}}\left(\int_0^T\int_\Omega \lvert\hat{Q}_t(\tau,x)\rvert^4\,dx\,d\tau\right)^{\frac{1}{4}}  \\
    &\leq C\left(\int_0^t\int_\Omega\vert r(\tau,x)-\hat{r}(\tau,x)\vert_F^2\,dx\,d\tau\right)^{\frac{1}{2}} \left(\int_0^t\int_\Omega\vert \nabla\hat{Q}(\tau,x)-\nabla Q(\tau,x)\vert_F^2\,dx\,d\tau\right)^{\frac{1}{2}}\left(\int_0^t\int_\Omega \lvert\nabla\hat{Q}_t(\tau,x)\rvert^2\,dx\,d\tau\right)^{\frac{1}{2}}\\ &\leq\int_0^t\int_\Omega\vert r(\tau,x)-\hat{r}(\tau,x)\vert^2\,dx\,d\tau+C\int_0^t\int_\Omega\vert \nabla\hat{Q}(\tau,x)-\nabla Q(\tau,x)\vert_F^2\,dx\,d\tau,
    \end{aligned}
\end{equation*}
where $C$ is a constant depending on $\|\hat{Q}_t\|_{L^2(0,T;H^1)}$ which is bounded by Corollary~\ref{cor:Q_H2_estimate}.
Combining the estimates for $K_1$ and $K_2$, we deduce that
\begin{align*}
    |J|\leq \int_0^t\int_\Omega \vert Q_t(\tau,x)&-\hat{Q}_t(\tau,x)\vert_F^2\,dx\,d\tau\\&+\int_0^t\int_\Omega\vert r(\tau,x)-\hat{r}(\tau,x)\vert^2\,dx\,d\tau+C\int_0^t\int_\Omega \vert \nabla\hat{Q}(\tau,x)-\nabla Q(\tau,x)\vert_F^2\,dx\,d\tau.
\end{align*}
Therefore, by substituting this inequality back to the estimate for $E_{rel}(t)$, we conclude that
\begin{equation*}
    E_{rel}(t)\leq C\int_0^t E_{rel}(\tau)\,d\tau.
\end{equation*}
Using the Gr\"onwall inequality and the fact that $E_{rel}(0)=0$, we obtain that $E_{rel}(t)\equiv 0$ for all $t\in[0,T]$ and  thus any weak solution of~\eqref{eq:qtensorflow} is in fact the unique strong solution.
\end{proof}

\section{Finite Element Scheme}
\label{sec:fem}
This section introduces a numerical approach for solving the reformulated system \eqref{eq:reformulation}. We will discuss a fully-discrete scheme using mass lumping, which allows the derivation of an energy dissipation law.

Let $\mathfrak{T}_h = \{R_{k}\}_{k}$ be a conforming shape-regular and quasi-uniform triangulation  of $\Omega$ with mesh size $h$, and let $\{\mathcal{N}_{z}\}_{z\in\mathcal{N}^{h}}$ be the set of interior nodes of the mesh with corresponding standard piecewise linear hat functions $\{\varphi_{z}\}_{z\in\mathcal{N}^{h}}$, where $\varphi_{z}(\mathcal{N}_{y}) = \delta_{yz}$ with $y\in\mathcal{N}^{h}$ (in an abuse of notation, we will write this as $\varphi_{z}(y)$).

Then we define
\begin{equation}
    \begin{split}
    \mathcal{S}^{h} &= \{W \in H_0^1(\Omega;\mathbb{R}^{d\times d}): \left.(W)_{ij}\right\rvert_{R_k} \in P_1(R_k) \text{ for all $R_k\in\mathfrak{T}_h, 0\le i,j\le d$}\},\\
    \mathcal{T}^{h} &= \{U\in H_0^1(\Omega): \left.U\right\rvert_{R_k} \in P_1(R_k) \text{ for all $R_k\in\mathfrak{T}_h$}\},
    \end{split}
\end{equation}
where $P_1(R_k)$ are the polynomials of degree up to 1 on $R_k$.
Let $E_{ij}$ be the $d\times d$ matrix with a one in entry $(i,j)$ and zeros elsewhere. Note that $\bigcup_{i,j=1}^{d} \{E_{ij}\varphi_{z}\}_{z\in\mathcal{N}^{h}}$ is a basis for $\mathcal{S}^{h}$ and $\{\varphi_{z}\}_{z\in\mathcal{N}^{h}}$ is a basis for $\mathcal{T}^{h}$.

We now define the mass lumping inner products and norms. Let $W, W^{1}, W^{2} \in \mathcal{S}^{h}$. Then we define
\[
    \langle W^{1}, W^{2}\rangle_{h} = \int_{\Omega} \mathcal{I}_{h}(W^{1}:W^{2})\, dx = \sum_{z\in\mathcal{N}_{h}} \gamma_{z} W^{1}(z):W^{2}(z),\quad \text{where }\gamma_{z} = \int_{\Omega} \varphi_{z}\, dx,
\]
and we have the induced norm
\[
    \lVert W \rVert_{h}^{2} = \sum_{z\in\mathcal{N}^{h}} \gamma_{z} \lvert W(z) \rvert_{F}^{2}.
\]
In an abuse of notation, we may write $\mathcal{I}_h W$ for $W\in\mathcal{S}^h$ to mean that $\mathcal{I}_h$ is applied to each entry of $W$. For $U, U^{1}, U^{2}\in\mathcal{T}^{h}$, we define
\[
    \langle U^{1}, U^{2}\rangle_{h} = \int_{\Omega} \mathcal{I}_{h}(U^{1}U^{2})\, dx = \sum_{z\in\mathcal{N}^{h}} \gamma_{z}U^{1}(z)U^{2}(z),
\]
with corresponding norm
\[
    \lVert U\rVert_{h}^{2} = \sum_{z\in\mathcal{N}^{h}} \gamma_{z}\lvert U(z)\rvert^{2}.
\]
Further define
\begin{equation}
    \langle \alpha(W^1), W^2\rangle_{h} 
        := -\int_{\Omega} \sum_{i,j,k=1}^{d}\left(\partial_{k} W_{jk}^1\cdot\partial_{i} W_{ij}^2 + \partial_{k} W_{ik}^1\cdot\partial_{j}W_{ij}^2\right)\, dx + \frac{2}{d}\int_\Omega \sum_{i,j,k,\ell=1}^d \partial_{\ell} W_{k\ell}^1 \cdot \partial_{k} W_{ij}^2\delta_{ij}\, dx.
\end{equation}
In what follows, we let $\Delta t > 0$, $t^n=n\Delta t$ and $T = N\Delta t$ for a final time of computation $T>0$ and number of time steps $N$. Furthermore, we have the forward and backward difference operators
\[
    D_t^+ g_h^n := \frac{g_h^{n+1} - g_h^n}{\Delta t},\quad D_t^- g_h^n := \frac{g_h^{n} - g_h^{n-1}}{\Delta t},
\]
where the superscript $n$ denotes the numerical approximations at time $n\Delta t$, and $n+1/2$ denotes the average of time steps:
\[
    Q_h^{n+1/2} = \frac{Q_h^{n+1} + Q_h^n}{2},\quad r_h^{n+1/2} = \frac{r_h^{n+1} + r_h^n}{2}.
\]
Lastly, $\alpha_h^{n+1/2}$ is defined as
\[
    \alpha_h^{n+1/2} = \frac{\alpha(Q_h^{n+1}) + \alpha(Q_h^n)}{2}.
\]
Now the discrete problem becomes: Find $Q_{h}^{n+1} \in\mathcal{S}^h$ and $r^{n+1}_h\in \mathcal{T}_h$ such that
\begin{equation}
\label{eq:fully-discrete-mass-lumping-fem-scheme}
    \begin{cases}
        \langle D_{t}^{+} Q_{h}^{n}, \Phi_h\rangle_{h}  &= -\sigma\langle D_{t}^{-}D_{t}^{+} Q_{h}^{n}, \Phi_h\rangle_{h} - L_{1}\langle \Grad Q_{h}^{n+\frac{1}{2}}, \Grad \Phi_h\rangle_{L^2(\Omega)}\\ 
        &\quad\quad+ \frac{(L_2+L_3)}{2}\langle \alpha_h^{n+1/2}, \Phi_h\rangle_h - \langle r_{h}^{n+\frac{1}{2}}P(Q_h^n), \Phi_h\rangle_{h}\\
        \langle r_{h}^{n+1} - r_{h}^{n}, \psi_h\rangle_{h} &= \langle P(Q_h^n) : (Q_{h}^{n+1} - Q_{h}^{n}), \psi_h\rangle_{h}
    \end{cases}
\end{equation}
for all $n=1,2,\dots,N$ and all $\Phi_h\in\mathcal{S}^h$, $\psi_h\in\mathcal{T}^h$, where
\begin{equation}
    Q_h^0 = \mathcal{I}_h Q_0^\sigma,\quad Q_h^1 = Q_h^0 + \Delta t \mathcal{I}_h Q_{t,0}^\sigma,\quad r_h^0 = \mathcal{I}_h r(Q_h^0),
\end{equation}
and $r_h^1$ satisfies
\[
    \langle r_h^1 - r_h^0, \psi_h\rangle_h = \langle P(Q_h^0) : (Q_h^1 - Q_h^0), \psi_h\rangle_h
\]
for all $\psi_h\in\mathcal{T}^h$. We will derive an energy estimate for this scheme, but first we will need to show that the scheme preserves the trace-free and symmetry properties of the Q-tensors. Note that $Q_h^0$ and $Q_h^1$ are trace-free and symmetric because they only involve interpolations of given trace-free and symmetric functions. 

\begin{lemma}[Trace-free and Symmetry Preservation]
    Assume that $Q_h^{n-1}$ and $Q_h^n$ are trace-free and symmetric. Then $Q_h^{n+1}$ computed via scheme \eqref{eq:fully-discrete-mass-lumping-fem-scheme} is trace-free and symmetric.
\end{lemma}
\begin{proof}
    Take $\Phi_h = Q_h^{n+1} - (Q_h^{n+1})^\top$ in equation \eqref{eq:fully-discrete-mass-lumping-fem-scheme}. Then since $\Phi_h$ is skew-symmetric and $Q_h^n$ is symmetric, we have that the left hand side of \eqref{eq:fully-discrete-mass-lumping-fem-scheme} becomes
    \[
        \langle D_{t}^{+} Q_{h}^{n}, \Phi_h\rangle_{h} = \frac{1}{\Delta t}\langle Q_h^{n+1}, \Phi_h\rangle_h.
    \]
    Similarly, the right hand side of \eqref{eq:fully-discrete-mass-lumping-fem-scheme} becomes
    \[
        -\frac{\sigma}{\Delta t^2}\langle Q_h^{n+1}, \Phi_h\rangle_h - \frac{ L_1}{2}\langle \Grad Q_h^{n+1}, \Grad\Phi_h\rangle_{L^2}+ \frac{ (L_2+L_3)}{4}\langle \alpha_h^{n+1}, \Phi_h\rangle_h,
    \]
    where we have used that $\alpha(Q_h^n)$ and $P(Q_h^n)$ are symmetric since $Q_h^n$ is symmetric. Thus
    \begin{equation}
        \frac{1}{\Delta t}\langle Q_h^{n+1}, \Phi_h\rangle_h = -\frac{\sigma}{\Delta t^2}\langle Q_h^{n+1}, \Phi_h\rangle_h - \frac{ L_1}{2}\langle \Grad Q_h^{n+1}, \Grad\Phi_h\rangle_{L^2}+ \frac{ (L_2+L_3)}{4}\langle \alpha(Q_h^{n+1}), \Phi_h\rangle_h.
    \end{equation}
    It follows that
    \begin{equation*}
        \frac{1}{\Delta t}\langle \Phi_h, \Phi_h\rangle_h = -\frac{\sigma}{\Delta t^2}\langle \Phi_h, \Phi_h\rangle_h - \frac{ L_1}{2}\langle \Grad \Phi_h, \Grad\Phi_h\rangle_{L^2}+ \frac{ (L_2+L_3)}{4}\langle \alpha(\Phi_h), \Phi_h\rangle_h,
    \end{equation*}
    in other words,
    \begin{equation}
        \left(\frac{1}{\Delta t}+\frac{\sigma}{\Delta t^2}\right)\|\Phi_h\|_h^2 + \frac{ L_1}{2}\|\Grad\Phi_h\|_{L^2}^2 = \frac{ (L_2+L_3)}{4}\langle\alpha(\Phi_h),\Phi_h\rangle_h.
    \end{equation}
    Now, using the fact that $\Phi$ is skew-symmetric, we have
    \begin{align*}
        \langle \alpha(\Phi_h), \Phi_h\rangle_h
            &= -\int_\Omega \sum_{i,j,k=1}^d \left(\partial_k\Phi_{jk}\partial_i\Phi_{ij} + \partial_k\Phi_{ik}\partial_j\Phi_{ij}\right)\, dx + \frac2d\int_\Omega\sum_{i,j,k,\ell=1}^d \partial_\ell\Phi_{k\ell}\partial_k\Phi_{ij}\delta_{ij}\, dx\\
            &= \int_\Omega \sum_{i,j,k=1}^d \partial_k \Phi_{kj}\partial_i\Phi_{ij}\, dx - \int_\Omega\sum_{i,j,k=1}^d \partial_k\Phi_{ik}\partial_j\Phi_{ij}\, dx + \frac2d\int_\Omega\sum_{i,k,\ell=1}^d \partial_\ell\Phi_{k\ell}\partial_k\Phi_{ii}\, dx\\
            &= \int_\Omega \lvert\operatorname{div}{\Phi_h}\rvert^2\, dx - \int_\Omega\lvert\operatorname{div}{\Phi_h}\rvert^2\, dx + 0\\
            &= 0.
    \end{align*}
    Combining this with the previous equation gives
    \[
        \left(\frac{1}{\Delta t}+\frac{\sigma}{\Delta t^2}\right)\|\Phi_h\|_h^2 + \frac{ L_1}{2}\|\Grad\Phi_h\|_{L^2}^2 = 0,
    \]
    which implies that $\Phi_h \equiv 0$, or $Q_h^{n+1} = (Q_h^{n+1})^\top$ so that $Q_h^{n+1}$ is symmetric.

    Now to show the trace-free property, take $\Phi_h = \operatorname{tr}(Q_h^{n+1})I$ in equation \eqref{eq:fully-discrete-mass-lumping-fem-scheme}. Because $\alpha$ is trace-free with no assumptions, and $P(Q_h^n)$ is trace-free since $Q_h^n$ is trace-free, we obtain
    \[
        \left(\frac{1}{\Delta t} + \frac{\sigma}{\Delta t^2}\right)\|\operatorname{tr}(Q_h^{n+1})\|_h^2 + \frac{L_1}{2}\|\Grad \operatorname{tr}(Q_h^{n+1})\|_{L^2}^2 = 0,
    \]
    which implies that $\operatorname{tr}(Q_h^{n+1})\equiv 0$, as desired.
\end{proof}
Having proved the trace-free and symmetry preservation properties of our scheme, we are now able to prove the following stability result:
\begin{theorem}[Energy Stability]
\label{thm:Fully-discrete-energy}
Define the energy
\begin{equation}
    E^{n+1} = \frac{\sigma}{2}\lVert D_{t}^{+} Q_{h}^{n}\rVert_{h}^{2} + \frac{ L_{1}}{2}\lVert \nabla Q_{h}^{n+1}\rVert_{L^{2}}^{2} + \frac{ (L_2+L_3)}{2}\|\operatorname{div}(Q_h^{n+1})\|_{L^2}^2 + \frac{1}{2} \lVert r_{h}^{n+1}\rVert_{h}^{2}.
\end{equation}
Then
\[
    E^{n+1} - E^{n} = -\Delta t \lVert D_{t}^{+} Q_{h}^{n}\rVert_{h}^{2} - \frac{\sigma}{2}\lVert D_{t}^{+} Q_{h}^{n} - D_{t}^{+} Q_{h}^{n-1}\rVert_{h}^{2}.
\]
\end{theorem}
\begin{proof}
Take $\Phi_h = -\Delta t D_{t}^{+} Q_{h}^{n}$ in equation \eqref{eq:fully-discrete-mass-lumping-fem-scheme}. Doing so, we obtain for the LHS
\[
    \langle D_{t}^{+} Q_{h}^{n}, -\Delta t D_{t}^{+} Q_{h}^{n}\rangle_{h} = -\Delta t\lVert D_{t}^{+} Q_{h}^{n}\rVert_{h}^{2}.
\]
For the first term on the RHS, we obtain
\begin{align*}
    \begin{split}
        -\sigma\langle D_{t}^{-} D_{t}^{+} Q_{h}^{n}, -\Delta t D_{t}^{+} Q_{h}^{n}\rangle_{h}
            &= \sigma \langle D_{t}^{+} Q_{h}^{n} - D_{t}^{+} Q_{h}^{n-1}, D_{t}^{+} Q_{h}^{n}\rangle_{h}\\
            &= \frac{\sigma}{2}\left(\langle D_{t}^{+} Q_{h}^{n} - D_{t}^{+} Q_{h}^{n-1}, D_{t}^{+} Q_{h}^{n} - D_{t}^{+} Q_{h}^{n-1}\rangle_{h}\right.\\ 
            &\quad\quad\left.+ \langle D_{t}^{+} Q_{h}^{n} - D_{t}^{+} Q_{h}^{n-1}, D_{t}^{+} Q_{h}^{n} + D_{t}^{+} Q_{h}^{n-1}\rangle_{h}\right)\\
            &= \frac{\sigma}{2}\left(\lVert D_{t}^{+} Q_{h}^{n} - D_{t}^{+} Q_{h}^{n-1}\rVert_{h}^{2} + \lVert D_{t}^{+} Q_{h}^{n}\rVert_{h}^{2} - \lVert D_{t}^{+} Q_{h}^{n-1}\rVert_{h}^{2}\right).
    \end{split}
\end{align*}
For the second term on the RHS, we have
\begin{align*}
    \begin{split}
        - L_1\langle \Grad Q_{h}^{n+\frac{1}{2}}, \Grad(-\Delta t D_{t}^{+} Q_{h}^{n})\rangle_{h}
            &= \frac{ L_{1}}{2} \langle \nabla Q_{h}^{n+1} + \nabla Q_{h}^{n}, \nabla Q_{h}^{n+1} - \nabla Q_{h}^{n}\rangle_{L^{2}}\\
            &= \frac{L_{1}}{2}\lVert \nabla Q_{h}^{n+1}\rVert_{L^{2}}^{2} - \frac{L_{1}}{2}\lVert \nabla Q_{h}^{n}\rVert_{L^{2}}^{2}.
    \end{split}
\end{align*}

For the $\alpha$ term, using the notation $Q_{ij}^n := (Q_h^n)_{ij}$, we have
\begin{align*}
    &\frac{(L_2+L_3)}{2}\langle \alpha_h^{n+1/2},-\Delta t D_t^+ Q_h^n\rangle_h\\
        &= -\frac{(L_2+L_3)}{2}\int_\Omega \sum_{i,j,k=1}^d (\partial_k Q_{jk}^{n+1/2} \partial_i (-\Delta t D_t^+ Q_{ij}^n) + \partial_k Q_{ik}^{n+1/2} \partial_j (-\Delta t D_t^+ Q_{ij}^n))\, dx\\ 
        &\quad+ \frac{(L_2+L_3)}{2}\cdot\frac{2}{d}\int_\Omega \sum_{i,j,k,\ell=1}^d \partial_\ell Q_{k\ell}^{n+1/2}\partial_k (-\Delta t D_t^+ Q_{ij}^n)\delta_{ij}\, dx\\
        &= \frac{(L_2+L_3)}{2}\int_\Omega \sum_{i,j,k=1}^d (\partial_k Q_{jk}^{n+1/2} \partial_i (Q_{ij}^{n+1}-Q_{ij}^n) + \partial_k Q_{ik}^{n+1/2} \partial_j (Q_{ij}^{n+1}-Q_{ij}^n))\, dx\\ 
        &\quad- \frac{(L_2+L_3)}{d}\int_\Omega \sum_{i,j,k,\ell=1}^d \partial_\ell Q_{k\ell}^{n+1/2}\partial_k (Q_{ii}^{n+1}-Q_{ii}^n)\, dx\\
        &= \frac{(L_2+L_3)}{2}\int_\Omega \sum_{i,j,k=1}^d (\partial_k Q_{jk}^{n+1/2} \partial_i (Q_{ij}^{n+1}-Q_{ij}^n) + \partial_k Q_{ik}^{n+1/2} \partial_j (Q_{ij}^{n+1}-Q_{ij}^n))\, dx\\
        &= \frac{(L_2+L_3)}{2}\int_\Omega \left(\lvert \operatorname{div}(Q_h^{n+1})\rvert^2 - \lvert \operatorname{div}(Q_h^{n})\rvert^2\right)\, dx,
\end{align*}
where in the last two lines, we used that $Q_h^{n+1}$ and $Q_h^n$ are trace-free and symmetric by the previous lemma.

Now for the last term on the RHS, first note that, taking $\psi_h = r_{h}^{n+\frac{1}{2}}$ in equation \eqref{eq:fully-discrete-mass-lumping-fem-scheme}, we have
\[
    \frac{1}{2}\lVert r_{h}^{n+1}\rVert_{h}^{2} - \frac{1}{2}\lVert r_{h}^{n}\rVert_{h}^{2} = \sum_{z\in\mathcal{N}^{h}} \gamma_{z} r_{h}^{n+\frac{1}{2}}(z)P(Q_h^n(z)):(Q_{h}^{n+1}(z) - Q_{h}^{n}(z)).
\]
Thus, we have
\begin{align*}
    \begin{split}
        -\langle r_{h}^{n+\frac{1}{2}}P(Q_h^n), -\Delta t D_{t}^{+} Q_{h}^{n}\rangle_{h}
            &= \sum_{z\in\mathcal{N}^{h}} \gamma_{z} r_{h}^{n+\frac{1}{2}}(z)P(Q_h^n(z)):(Q_{h}^{n+1}(z) - Q_{h}^{n}(z))\\
            &= \frac{1}{2}\lVert r_{h}^{n+1}\rVert_{h}^{2} - \frac{1}{2}\lVert r_{h}^{n}\rVert_{h}^{2}.
    \end{split}
\end{align*}
In summary,
\[
    \begin{split}
    -\Delta t \lVert D_{t}^{+} Q_{h}^{n}\rVert_{h}^{2} 
        &= \frac{\sigma}{2}\left(\lVert D_{t}^{+} Q_{h}^{n} - D_{t}^{+} Q_{h}^{n-1}\rVert_{h}^{2} + \lVert D_{t}^{+} Q_{h}^{n}\rVert_{h}^{2} - \lVert D_{t}^{+} Q_{h}^{n-1}\rVert_{h}^{2}\right)\\ 
        &\quad\quad + \frac{L_{1}}{2}\lVert \nabla Q_{h}^{n+1}\rVert_{L^{2}}^{2} - \frac{L_{1}}{2} \lVert \nabla Q_{h}^{n}\rVert_{L^{2}}^{2}\\ 
        &\quad\quad + \frac{(L_2+L_3)}{2}\|\operatorname{div}(Q_h^{n+1})\|_{L^2}^2 - \frac{(L_2+L_3)}{2}\|\operatorname{div}(Q_h^{n})\|_{L^2}^2\\
        &\quad\quad+ \frac{1}{2}\lVert r_{h}^{n+1}\rVert_{h}^{2} - \frac{1}{2}\lVert r_{h}^{n}\rVert_{h}^{2}
    \end{split}
\]
as desired.
\end{proof}

\subsection{Convergence to the strong solution} Building upon the discussion above, our goal is to demonstrate the convergence of the fully-discrete scheme. First, we must construct an approximating numerical solution from the discrete values the fully-discrete system obtains. To achieve this, we employ standard piecewise constant interpolation to construct the following functions:
\begin{align*}
    \label{eq:Qrsol}
    &Q_{h,\Delta t}(t,x) = \sum_{n=0}^{N-1}  Q_h^n\chi_{S_n}(t),\\
    &r_{h,\Delta t}(t,x) = \sum_{n=0}^{N-1} r_h^n(x)\chi_{S_n}(t)  
\end{align*}
Here $S_n=[n\Delta t, (n+1)\Delta t)$, and $\chi_{S_n}$ is the characteristic function defined on the interval $S_n$.

\begin{theorem}[Converge Rate]\label{thm:convergence-rate}
    Let $Q_0 \in H^4(\Omega)$ and $Q_{t,0} \in H^3(\Omega)$. Then the numerical solution $(Q_h, r_h)$ converges to the strong solution $(\hat Q, \hat r)$ with
    \[
        \|Q_h(t) - \hat Q(t)\|_{H^1} = O(\Delta t^{1/2} + h),\quad \|r_h(t) - \hat r(t)\|_{L^2} = O(\Delta t^{1/2} + h),\quad \|D_t^- Q_h(t) - \hat Q_t(t)\|_{L^2} = O(\Delta t^{1/2} + h)
    \]
    for any $t\in[0,T]$.
\end{theorem}

\begin{proof}
Let $U_h = (Q_h, r_h)$ be the discrete solution given by the scheme \eqref{eq:fully-discrete-mass-lumping-fem-scheme}, constantly interpolated in time, and let $\hat{U} = (\hat{Q}, \hat{r})$ be the strong solution of \eqref{eq:qtensorflow}. Define the relative entropy
\[
    \mathcal{E}(U_h \vert \hat{U}) = E(U_h) - E(\hat{U}) - dE(\hat{U})(U_h - \hat{U}),
\]    
where
\[
    E(U_h(t)) = \int_\Omega \left(\frac{L_1}{2}|\nabla Q_h(t)|^2 + \frac{L_2+L_3}{2}|\Div Q_h(t)|^2 + \frac{\sigma}{2}\mathcal{I}_h |D_t^c Q_h(t)|^2 + \frac12 \mathcal{I}_h|r_h(t)|^2\right)\, dx,
\]
with $\mathcal{I}_h$ is a piecewise linear interpolation in space on the finite element mesh and
\[
    D_t^c Q_h(t) = \frac{Q_h(t+\Delta t/2) - Q_h(t - \Delta t/2)}{\Delta t}.
\]
We can extend the strong solution using Taylor series and the discrete solution using interpolating polynomials to ensure that $E(U_h(t))$ and $E(\hat{U}(t))$ are defined also for $t\in[0, \Delta t/2)$ so that the extensions equal $Q_0$ up to $O(\Delta t^2)$. Then computing $dE(\hat{U})(U_h - \hat{U})$, we have
\begin{align*}
    &dE(\hat{U})(U_h - \hat{U})\\
        &= \int_\Omega \Bigg(L_1 \nabla \hat{Q} : (\nabla Q_h - \nabla \hat{Q}) + (L_2 + L_3) \Div \hat{Q} \cdot (\Div Q_h -\Div \hat{Q}) + \sigma \mathcal{I}_h\left(D_t^c \hat{Q} : (D_t^c Q_h - D_t^c \hat{Q})\right)\\
        &\qquad \quad+ \mathcal{I}_h\left(\hat{r}(r_h - \hat{r})\right)\Bigg)\, dx
\end{align*}
so that
\begin{equation}
\label{eq:discrete-entropy-equality}
\begin{split}
    &\mathcal{E}(U_h\vert\hat{U})\\
        &= \int_\Omega \left(\frac{L_1}{2}|\nabla Q_h|^2 + \frac{L_2+L_3}{2}|\Div Q_h|^2 + \frac{\sigma}{2}\mathcal{I}_h |D_t^c Q_h|^2 + \frac12 \mathcal{I}_h|r_h|^2\right)\, dx\\
        &\hspace{10ex}- \int_\Omega \left(\frac{L_1}{2}|\nabla \hat{Q}|^2 + \frac{L_2+L_3}{2}|\Div \hat{Q}|^2 + \frac{\sigma}{2}\mathcal{I}_h |D_t^c \hat{Q}|^2 + \frac12 \mathcal{I}_h|\hat{r}|^2\right)\, dx\\
        &\hspace{10ex}-\int_\Omega \Bigg(L_1 \nabla \hat{Q} : (\nabla Q_h - \nabla \hat{Q}) + (L_2 + L_3) \Div \hat{Q} \cdot (\Div Q_h - \Div \hat{Q})\\ 
        &\hspace{15ex}+ \sigma \mathcal{I}_h\left(D_t^c \hat{Q} : (D_t^c Q_h - D_t^c \hat{Q})\right) + \mathcal{I}_h\left(\hat{r}(r_h - \hat{r})\right)\Bigg)\, dx\\
        &= \frac{L_1}{2}\int_\Omega |\nabla Q_h - \nabla \hat{Q}|^2\, dx + \frac{L_2+L_3}{2}\int_\Omega |\Div Q_h - \Div \hat{Q}|^2\, dx + \frac{\sigma}{2}\int_\Omega \mathcal{I}_h|D_t^c Q_h - D_t^c \hat{Q}|^2\, dx \\
        &\qquad + \frac12\int_\Omega \mathcal{I}_h|r_h - \hat{r}|^2\, dx\\
        &= E(U_h(t)) + E(\hat{U}(t)) - dE(\hat{U}(t))U_h(t).
\end{split}
\end{equation}
Now, we estimate $E(U_h(t^{n+1}))$. Using the discrete energy estimate Theorem \ref{thm:Fully-discrete-energy}, we have
\begin{align*}
    E(U_h(t^{n+1}))
        &= \int_\Omega \left[\frac{L_1}{2}|\nabla Q_h(t^{n+1})|^2 + \frac{L_2+L_3}{2}|\Div Q_h(t^{n+1})|^2 + \frac{\sigma}{2}\mathcal{I}_h|D_t^c Q_h(t^{n+1})|^2 + \frac12\mathcal{I}_h|r_h(t^{n+1})|^2\right]\, dx\\
        &= E^{n+1}\\
        &= E^0 -\Delta t \sum_{k=0}^{n} \int_\Omega \lvert D_{t}^{+} Q_{h}^{k}\rvert^2\, dx - \frac{\sigma}{2}\sum_{k=0}^{n} \int_\Omega \mathcal{I}_h\lvert D_{t}^{+} Q_{h}^{n} - D_{t}^{+} Q_{h}^{n-1}\rvert^2\, dx\\
        &\le E(U_h(0)) -\Delta t \sum_{k=0}^{n} \int_\Omega \lvert D_{t}^{+} Q_{h}^{k}\rvert^2.
\end{align*}
Similarly, we have
\begin{align*}
    E(\hat{U}(t^{n+1}))
        &= \int_\Omega \left[\frac{L_1}{2}|\nabla\hat{Q}(t^{n+1})|^2 + \frac{L_2+L_3}{2}|\Div\hat{Q}(t^{n+1})|^2 + \frac{\sigma}{2}\mathcal{I}_h|D_t^c\hat{Q}(t^{n+1})|^2 + \frac12\mathcal{I}_h|\hat r(t^{n+1})|^2\right]\, dx\\
        &= O(\Delta t^2 + h^2) + \frac12\int_\Omega \left(L_1|\nabla\hat{Q}(t^{n+1})|^2 + (L_2 + L_3)|\Div\hat{Q}(t^{n+1})|^2 + \sigma|\hat{Q}_t(t^{n+1})|^2 + \hat{r}(t^{n+1})^2\right)\, dx\\
        &= O(\Delta t^2 + h^2) + \frac12\int_\Omega \left(L_1|\nabla\hat{Q}(0)|^2 + (L_2 + L_3)|\Div\hat{Q}(0)|^2 + \sigma|\hat{Q}_t(0)|^2 + \hat{r}(0)^2\right)\, dx \\
        &\qquad- \int_0^{t^{n+1}}\int_\Omega \lvert \partial_t\hat Q\rvert^2\, dx\, d\tau\\
        &= O(\Delta t^2 + h^2) + E(\hat{U}(0)) - \int_0^{t^{n+1}}\int_\Omega \lvert \partial_t\hat Q\rvert^2\, dx\, d\tau.
\end{align*}
Now it remains to bound $dE(\hat{U}(t^{n+1}))U_h(t^{n+1})$. Denoting $\hat{U}^k := \hat{U}(t^k)$ and $U_h^k = U_h(t^k)$, we have
\begin{align*}
    &dE(\hat{U}(t^{n+1}))U_h(t^{n+1})\\
        &= \int_\Omega \left(L_1\nabla\hat{Q}^{n+1} : \nabla Q_h^{n+1} + (L_2 + L_3)\Div\hat{Q}^{n+1} \cdot \Div Q_h^{n+1} + \sigma\mathcal{I}_h\left(D_t^c\hat{Q}^{n+1} : D_t^c Q_h^{n+1}\right) + \mathcal{I}_h\left(\hat{r}^{n+1}r_h^{n+1}\right)\right)\, dx\\
        &= dE(\hat{U}(0))U_h(0) + \Delta t\sum_{k=0}^n D_t^+ (dE(\hat{U}^k)U_h^k)\\
        &= dE(\hat{U}(0))U_h(0) + \Delta t\sum_{k=0}^n \int_\Omega \Bigg(L_1 D_t^+ \nabla\hat{Q}^k : \nabla Q_h^{k+1/2} + L_1 \nabla\hat{Q}^{k+1/2} : D_t^+ \nabla Q_h^k\\
        &\hspace{20ex} +(L_2+L_3) D_t^+ \Div\hat{Q}^k \cdot \Div Q_h^{k+1/2} + (L_2+L_3) \Div\hat{Q}^{k+1/2} \cdot D_t^+ \Div Q_h^k\\
        &\hspace{20ex} + \sigma\mathcal{I}_h\left(D_t^+ D_t^c \hat{Q}^k : D_t^c Q_h^{k+1}\right) + \sigma\mathcal{I}_h\left(D_t^c \hat{Q}^{k} : D_t^+D_t^c Q_h^k\right)\\
        &\hspace{20ex} + \mathcal{I}_h\left(r_h^{k+1/2}D_t^+\hat{r}^k\right) + \mathcal{I}_h\left(\hat{r}^{k+1/2} D_t^+ r_h^k\right)\Bigg)\, dx.
\end{align*}
Now let $\Pi_h:H_0^1(\Omega, \R^{d\times d})\to\mathcal{S}^h$ be a Galerkin projection onto the finite element space defined using the inner product
\[
    \langle Q_1, Q_2\rangle_L: = L_1 \int_\Omega \nabla Q_1 : \nabla Q_2\, dx + (L_2+L_3)\int_\Omega \Div Q_1 \cdot \Div Q_2\, dx,
\]
i.e., it satisfies $\langle \Pi_h Q_1,Q_2\rangle_L=\langle Q_1,Q_2\rangle_L$ for all $Q_2\in \mathcal{S}_h$, and let $\tilde\Pi_h: H_0^1(\Omega)\to\mathcal{T}^h$ be the standard $L^2$ projection.
Then from the definition of $\Pi_h$, we have
\begin{align*}
    &\int_\Omega L_1 D_t^+ \nabla\hat{Q}^k : \nabla Q_h^{k+1/2}\, dx + \int_\Omega (L_2 + L_3)D_t^+\Div \hat{Q}^k \cdot \Div Q_h^{k+1/2}\, dx\\
        &= \int_\Omega L_1 D_t^+ \nabla\Pi_h\hat{Q}^k : \nabla Q_h^{k+1/2}\, dx + \int_\Omega (L_2 + L_3) D_t^+\Div \Pi_h \hat{Q}^k \cdot \Div Q_h^{k+1/2}\, dx.
\end{align*}
Furthermore,
\begin{align}\label{eq:doubleo}
	\begin{split}
    \int_\Omega \mathcal{I}_h\left(D_t^+ D_t^c \hat{Q}^{k} : D_t^c Q_h^{k+1}\right)\, dx
        &= O(h^2) + \int_\Omega D_t^+ D_t^c \hat{Q}^{k} : D_t^c Q_h^{k+1}\, dx\\
        &= O(h^2) + \int_\Omega D_t^+ D_t^c \Pi_h\hat{Q}^{k} : D_t^c Q_h^{k+1}\, dx\\
        &= O(h^2) + \int_\Omega \mathcal{I}_h\left(D_t^+ D_t^c \Pi_h\hat{Q}^{k} : D_t^c Q_h^{k+1}\right)\, dx,
    \end{split}
\end{align}
and
\begin{equation}
	\label{eq:singleo}
    \int_\Omega \mathcal{I}_h\left(r_h^{k+1/2}D_t^+ \hat{r}^k\right)\, dx = O(h^2) + \int_\Omega \mathcal{I}_h\left(r_h^{k+1/2}D_t^+ \tilde\Pi_h \hat{r}^k\right)\, dx.
\end{equation}

Using the property of the projection $\Pi_h$ and substituting the numerical scheme \eqref{eq:fully-discrete-mass-lumping-fem-scheme}, we have
\begin{align}\label{eq:triangle}
    \begin{split}
    &\int_\Omega \Bigg(L_1 D_t^+\nabla\hat{Q}^k: \nabla Q_h^{k+1/2} + (L_2 + L_3)D_t^+\Div\hat{Q}^k \cdot \Div Q_h^{k+1/2}\Bigg)\, dx\\
    &= \int_\Omega \Bigg(L_1 D_t^+\nabla\Pi_h \hat{Q}^k :\nabla Q_h^{k+1/2} + (L_2 + L_3)D_t^+\Div\Pi_h \hat{Q}^k \cdot \Div Q_h^{k+1/2}\Bigg)\, dx\\
    &= \int_\Omega \Bigg(-\mathcal{I}_h\left(D_t^+ \Pi_h \hat{Q}^k:D_t^+ Q_h^k\right) - \sigma \mathcal{I}_h\left(D_t^+ \Pi_h \hat{Q}^k :D_t^-D_t^+ Q_h^k\right)\\ 
    &\hspace{25ex} - \mathcal{I}_h\left((D_t^+ \Pi_h \hat{Q}^k)(r_h^{k+1/2} P(Q_h^k))\right)\Bigg)dx.  
    \end{split}
\end{align}
Then using the strong form of the solution, we have
\begin{align}\label{eq:doubletriangle}
	\begin{split}
    &\int_\Omega \left(L_1\nabla\hat Q^{k+1/2} :D_t^+ \nabla Q_h^k +(L_2+L_3)\Div \hat{Q}^{k+1/2}\cdot  D_t^+\Div Q^k_h\right) dx\\
    & =-\int_\Omega \left(L_1\Delta \hat Q^{k+1/2} :D_t^+  Q_h^k +\frac{L_2+L_3}{2}\alpha( \hat{Q}^{k+1/2}) : D_t^+ Q^k_h\right) dx\\
    &= O(\Delta t^2) + \int_\Omega \left[-\hat Q_t(t^{k+1/2}) - \sigma\hat Q_{tt}(t^{k+1/2}) - \hat r(t^{k+1/2})P(\hat Q(t^{k+1/2}))\right] : D_t^+ Q_h^k\, dx,
\end{split}
\end{align}
where $t^{k+1/2} = (t^k + t^{k+1})/2$. Thus, using the fact that $D_t^c Q_h^k = D_t^- Q_h^k$ due to the piecewise constant interpolation, and~\eqref{eq:singleo},~\eqref{eq:doubleo},~\eqref{eq:triangle} and~\eqref{eq:doubletriangle}, we have
\begin{align*}
    &dE(\hat U(t))U_h(t) - dE(\hat U(0))U_h(0)\\
        &= \Delta t\sum_{k=0}^n\Bigg[ O(\Delta t^2 + h^2) + \int_\Omega \Bigg(-\mathcal{I}_h\left(D_t^+\Pi_h\hat Q^k : D_t^+ Q_h^k\right) - \partial_t\hat Q(t^{k+1/2}) : D_t^+ Q_h^k\\
        &\hspace{20ex} + \sigma\mathcal{I}_h \left(D_t^+D_t^c\Pi_h\hat{Q}^{k} : D_t^c Q_h^{k+1}\right) - \sigma \partial_{tt}\hat Q(t^{k+1/2}) : D_t^+ Q_h^k\\
        &\hspace{20ex}+ \sigma\mathcal{I}_h\left(D_t^c \hat{Q}^{k} : D_t^+ D_t^- Q_h^{k}\right) - \sigma \mathcal{I}_h\left(D_t^+ \Pi_h \hat{Q}^k : D_t^-D_t^+ Q_h^k\right)\\
        &\hspace{20ex}+ \mathcal{I}_h\left(r_h^{k+1/2} D_t^+\tilde\Pi_h\hat{r}^k\right) - \mathcal{I}_h\left(D_t^+ \Pi_h \hat{Q}^k : (r_h^{k+1/2} P(Q_h^k))\right)\\
        &\hspace{20ex}+ \mathcal{I}_h\left(\hat{r}^{k+1/2} D_t^+ r_h^k\right) - \hat r(t^{k+1/2}) P(\hat Q(t^{k+1/2})) : D_t^+ Q_h^k\Bigg)\, dx\Bigg]
\end{align*}
so that
\begin{align*}
    \mathcal{E}(U_h(t)\vert\hat U(t))
        &\le \mathcal{E}(U_h(0)\vert\hat U(0)) + O(\Delta t^2 + h^2) - \Delta t \sum_{k=0}^{n}\int_\Omega \lvert D_t^+ Q_h^k\rvert^2\, dx - \int_0^{t^{n+1}} \int_\Omega \lvert \partial_t \hat{Q}\rvert^2\, dx\, d\tau\\ 
        &\hspace{10ex}- \Delta t\sum_{k=0}^n \int_\Omega \Bigg(-\mathcal{I}_h\left(D_t^+\Pi_h\hat Q^k : D_t^+ Q_h^k\right) - \partial_t\hat Q(t^{k+1/2}) : D_t^+ Q_h^k\\
        &\hspace{20ex} + \sigma\mathcal{I}_h \left(D_t^+D_t^c \Pi_h\hat{Q}^{k} : D_t^c Q_h^{k+1}\right) - \sigma \partial_{tt}\hat Q(t^{k+1/2}) : D_t^+ Q_h^k\\
        &\hspace{20ex}+ \sigma\mathcal{I}_h\left(D_t^c \hat{Q}^{k} :D_t^+D_t^- Q_h^{k}\right) - \sigma\mathcal{I}_h\left(D_t^+ \Pi_h \hat{Q}^k :D_t^-D_t^+ Q_h^k\right)\\
        &\hspace{20ex}+ \mathcal{I}_h\left(r_h^{k+1/2} D_t^+\tilde\Pi_h\hat{r}^k\right) - \mathcal{I}_h\left(D_t^+ \Pi_h \hat{Q}^k : (r_h^{k+1/2} P(Q_h^k))\right)\\
        &\hspace{20ex}+ \mathcal{I}_h\left(\hat{r}^{k+1/2} D_t^+ r_h^k\right) - \hat r(t^{k+1/2})P(\hat Q(t^{k+1/2})) : D_t^+ Q_h^k\Bigg)\, dx.
\end{align*}
We aim to show that the integral and sum terms are bounded by $O(\Delta t + h^2) + \Delta t\sum_{k=0}^n \mathcal{E}(U_h(t^k)\vert \hat U(t^k))$ to apply a discrete Gr\"{o}nwall inequality. First, using the midpoint rule and Lemma \ref{lem:thirdtimebound}, we have that
\begin{align*}
    \int_0^{t^{n+1}}\int_\Omega \lvert \partial_t \hat Q\rvert^2\, dx\, d\tau = O(\Delta t^2) + \Delta t\sum_{k=0}^n\|\partial_t \hat{Q}(t^{k+1/2} )\|_{L^2}^2.
\end{align*}
Thus,
\begin{align*}
    &- \Delta t \sum_{k=0}^{n}\int_\Omega \lvert D_t^+ Q_h^k\rvert^2\, dx - \int_0^{t^{n+1}} \int_\Omega \lvert \hat Q_t\rvert^2\, dx\, d\tau + \Delta t\sum_{k=0}^n \int_\Omega \left(\mathcal{I}_h\left(D_t^+\Pi_h\hat Q^k : D_t^+ Q_h^k\right) + \partial_t\hat Q(t^{k+1/2}) : D_t^+ Q_h^k\right)\, dx\\
    &\le O(\Delta t^2 + h^2) - \Delta t \sum_{k=0}^{n}\int_\Omega \lvert D_t^+ Q_h^k\rvert^2\, dx - \Delta t \sum_{k=0}^n \|\partial_t \hat Q(t^{k+1/2})\|^2\\ 
    &\hspace{10ex}+ \Delta t \sum_{k=0}^n \left[\frac12 \|D_t^+ \hat Q^k\|_{L^2}^2 + \frac12 \|\partial_t \hat Q(t^{k+1/2})\|_{L^2}^2 + \|D_t^+Q_h^k\|_{L^2}^2\right]\\
    &= O(\Delta t^2 + h^2),
\end{align*}
where we used that
\[
    \|D_t^+ \hat Q^k\|_{L^2}^2 = \|D_t^c \hat Q(t^k+\Delta t/2)\|_{L^2}^2 = \|\partial_t \hat Q(t^{k+1/2})\|_{L^2}^2 + O(\Delta t^2).
\]
Also,
\begin{align*}
    &\int_\Omega\left(\sigma\mathcal{I}_h \left(D_t^+D_t^c\Pi_h\hat{Q}^{k} :D_t^c Q_h^{k+1}\right) - \sigma \partial_{tt}\hat Q(t^{k+1/2}) :D_t^+ Q_h^k\right)\, dx\\
        &= O(h^2) + \sigma\int_\Omega \left(D_t^+ D_t^c \hat Q^k :D_t^+ Q_h^{k} - \partial_{tt}\hat Q(t^{k+1/2}) :D_t^+ Q_h^k\right)\, dx\\
        &\le O(h^2) + \sigma \|D_t^+ Q_h^k\|_{L^2}\|D_t^+ D_t^c \hat Q^k - \partial_{tt} \hat Q(t^{k+1/2})\|_{L^2}\\
        &= O(\Delta t^2 + h^2).
\end{align*}
Here we used Lemma \ref{lem:fourthtimebound} so that it follows
\begin{align*}
    D_t^+ D_t^c \hat Q^k &= D_t^+ \left[\frac{\hat Q(t^{k+1/2}) - \hat Q(t^{k+1/2}-\Delta t)}{\Delta t}\right] = \frac{\hat Q(t^{k+1/2}+\Delta t) - 2\hat Q(t^{k+1/2}) + \hat Q(t^{k+1/2}-\Delta t)}{\Delta t^2}\\
    &= \partial_{tt} \hat Q(t^{k+1/2}) + O(\Delta t^2).
\end{align*}
We also have by summation by parts that
\begin{align*}
    &\Delta t \sum_{k=0}^n \int_\Omega\left(\sigma\mathcal{I}_h\left(D_t^c \hat{Q}^{k} : D_t^+D_t^- Q_h^{k}\right) - \sigma\mathcal{I}_h\left(D_t^+ \Pi_h \hat{Q}^k : D_t^-D_t^+ Q_h^k\right)\right)\, dx\\
        &= O(h^2) + \sigma\Delta t\sum_{k=0}^n\int_\Omega D_t^+ D_t^- Q_h^k : (D_t^c \hat Q^{k} - D_t^+ \hat Q^k)\, dx\\
        &= O(h^2) + \sigma\int_\Omega D_t^+ Q_h^n : (D_t^c \hat{Q}^n - D_t^+ \hat{Q}^n)\, dx - \sigma\int_\Omega D_t^+ Q_h^{-1} : (D_t^c \hat{Q}^0 - D_t^+ \hat{Q}^0)\, dx\\ 
        &\hspace{10ex}- \sigma\sum_{k=0}^{n-1}\int_\Omega D_t^+ Q_h^k : \left((D_t^c \hat{Q}^{k+1} - D_t^+ \hat{Q}^{k+1}) - (D_t^c \hat{Q}^k - D_t^+ \hat{Q}^k)\right)\, dx\\
        &\le O(h^2) + \sigma \|D_t^+ Q_h^n\|_{L^2}\left\| D_t^c \hat{Q}^n - D_t^+ \hat{Q}^n\right\|_{L^2} + \sigma\|D_t^+ Q_h^{-1}\|_{L^2}\left\|D_t^c \hat{Q}^0 - D_t^+ \hat{Q}^0\right\|_{L^2}\\
        &\hspace{10ex} + \sigma\Delta t\sum_{k=0}^{n-1} \|D_t^+ Q_h^k\|_{L^2}\underbrace{\left\| D^+_t D_t^c \hat{Q}^{k} -D^+_t D_t^+ \hat{Q}^{k}\right\|_{L^2}}_{=O(\Delta t)}\\
        &\le O(\Delta t + h^2),
\end{align*}
where we have used Lemma~\ref{lem:thirdtimebound}.
Finally,
\begin{align*}
    &\int_\Omega \Bigg(\mathcal{I}_h\left(r_h^{k+1/2}D_t^+\tilde\Pi_h\hat{r}^k\right) - \mathcal{I}_h\left(D_t^+ \Pi_h \hat{Q}^k : (r_h^{k+1/2} P(Q_h^k))\right)\\ 
    &\hspace{10ex}+ \mathcal{I}_h\left(\hat{r}^{k+1/2} D_t^+ r_h^k\right) - \hat r(t^{k+1/2})P(\hat Q(t^{k+1/2})) : D_t^+ Q_h^k\Bigg)\, dx\\
    &= O(h^2) + \int_\Omega \Bigg(r_h^{k+1/2}D_t^+\tilde\Pi_h\hat{r}^k - D_t^+ \Pi_h \hat{Q}^k : (r_h^{k+1/2} P(Q_h^k)) + \hat{r}^{k+1/2} D_t^+ r_h^k\\ 
    &\hspace{10ex}- \hat r(t^{k+1/2})P(\hat Q(t^{k+1/2})) : D_t^+ Q_h^k\Bigg)\, dx\\
    &=  O(h^2) + \int_\Omega \Bigg(r_h^{k+1/2} D_t^+\hat{r}^k - D_t^+ \hat{Q}^k : (r_h^{k+1/2} P(Q_h^k)) + \hat{r}^{k+1/2}P(Q_h^k) : D_t^+Q_h^k\\ 
    &\hspace{10ex}- \hat r(t^{k+1/2})P(\hat Q(t^{k+1/2})) : D_t^+ Q_h^k\Bigg)\, dx\\
    &=  O(\Delta t^2 + h^2) + \int_\Omega \Bigg(\hat{r}_t(t^{k+1/2})r_h^{k+1/2} - \hat{Q}_t(t^{k+1/2}) : (r_h^{k+1/2} P(Q_h^k)) + \hat{r}^{k+1/2}P(Q_h^k) : D_t^+Q_h^k\\
    &\hspace{10ex}- \hat r(t^{k+1/2})P(\hat Q(t^{k+1/2})) : D_t^+ Q_h^k\Bigg)\, dx\\
    &= O(\Delta t^2 + h^2) + \int_\Omega\Bigg(P(\hat Q(t^{k+1/2})):\hat Q_t(t^{k+1/2}) r_h^{k+1/2} - \hat Q_t(t^{k+1/2}) : (r_h^{k+1/2} P(Q_h^k))\\ 
    &\hspace{10ex}+ \hat r^{k+1/2}P(Q_h^k):D_t^+ Q_h^k - \hat r(t^{k+1/2})P(\hat Q(t^{k+1/2})) : D_t^+ Q_h^k\Bigg)\, dx\\
    &= O(\Delta t^2 + h^2) + \int_\Omega\left((P(\hat Q(t^{k+1/2})) - P(Q_h^k)):(r_h^{k+1/2}\hat Q_t(t^{k+1/2}) - \hat r^{k+1/2} D_t^+ Q_h^k)\right)\, dx\\
    &\le O(\Delta t^2 + h^2) + \int_\Omega \lvert P(\hat Q(t^{k+1/2})) - P(Q_h^k)\rvert \cdot \lvert r_h^{k+1/2}\hat Q_t(t^{k+1/2}) - \hat r^{k+1/2} D_t^+ Q_h^k\rvert\, dx\\
    &= O(\Delta t^2 + h^2) + \int_\Omega \lvert P(\hat Q(t^{k+1/2})) - P(Q_h^k)\rvert \cdot \lvert \hat r^{k+1/2}(D_t^+ Q_h^k - \hat Q_t(t^{k+1/2})) + \hat Q_t(t^{k+1/2})(\hat r^{k+1/2} - r_h^{k+1/2})\rvert\, dx.
\end{align*}
Using Lipschitz continuity of $P$, the Poincar\'e inequality, and the regularity of the strong solution, we bound this by
\begin{align*}
    &\le O(\Delta t^2 + h^2) + C\int_\Omega\lvert \hat Q(t^{k+1/2}) - Q_h^k\rvert\|\hat r^{k+1/2}\|_{L^\infty} \lvert D_t^+ Q_h^k - \hat Q_t(t^{k+1/2})\rvert\, dx\\ 
    &\hspace{10ex}+ C\int_\Omega \lvert \hat Q(t^{k+1/2}) - Q_h^k\rvert\|\hat Q_t(t^{k+1/2})\|_{L^\infty}\lvert \hat r^{k+1/2} - r_h^{k+1/2}\rvert\, dx\\
    &\le O(\Delta t^2 + h^2) + C\|\hat r^{k+1/2}\|_{L^\infty}\|\hat Q(t^{k+1/2}) - Q_h(t^{k+1/2})\|_{L^2}^2 + C\|\hat r^{k+1/2}\|_{L^\infty} \|D_t^+ Q_h^k - D_t^+ \hat Q^k\|_{L^2}^2\\
    &\hspace{10ex}+ C\|\hat Q_t(t^{k+1/2})\|_{L^\infty}\|\hat Q(t^{k+1/2}) - Q_h(t^{k+1/2})\|_{L^2}^2 + C\|\hat Q_t(t^{k+1/2})\|_{L^\infty}\|\hat r^{k+1/2} - r_h^{k+1/2}\|_{L^2}^2\\
    &\le O(\Delta t^2 + h^2) + C\|\nabla \hat Q(t^{k+1/2}) - \nabla Q_h(t^{k+1/2})\|_{L^2}^2 + C\|D_t^+ Q_h^k - D_t^+ \hat Q^k\|_{L^2}^2 + C\|\hat r^k - r_h^k\|_{L^2}^2 \\
    &\qquad+ C\|\hat r^{k+1} - r_h^{k+1}\|_{L^2}^2\\
    &\le O(\Delta t^2 + h^2) + C\mathcal{E}(U_h(t^k)\vert \hat U(t^k)) + C\mathcal{E}(U_h(t^{k+1/2})\vert \hat U(t^{k+1/2})) + C\mathcal{E}(U_h(t^{k+1})\vert \hat U(t^{k+1})).
\end{align*}
Here we used the fact that $Q_h$ is a piecewise constant interpolation in time so that $Q_h^k = Q_h(t^k) = Q_h(t^{k+1/2})$. Combining all of these inequalities gives
\begin{align*}
    \mathcal{E}(U_h(t^{n+1})\vert\hat U(t^{n+1}))
        &\le C(\Delta t + h^2) + \mathcal{E}(U_h(0)\vert\hat U(0)) + C\Delta t \sum_{k=0}^n \left[\mathcal{E}(U_h(t^k)\vert\hat U(t^k)) + \mathcal{E}(U_h(t^{k+1/2})\vert \hat{U}(t^{k+1/2}))\right]
\end{align*}
for $\Delta t, h$ sufficiently small and some constant $C>0$. Then by the discrete Gr\"{o}nwall inequality, we have
\begin{align*}
    \mathcal{E}(U_h(t^{n+1})\vert\hat U(t^{n+1}))
        &\le \left(C(\Delta t+h^2) + \mathcal{E}(U_h(0)\vert\hat U(0))\right)(1 + CTe^{CT}) = C(T)\left(\Delta t + h^2 + \mathcal{E}(U_h(0)\vert \hat U(0))\right).
\end{align*}
Then since $Q_h(0) = \mathcal{I}_h \hat Q(0)$ and $r_h(0) = \mathcal{I}_h \hat r(0)$, we have that
\begin{align*}
    \mathcal{E}(U_h(0)\vert \hat U(0))
        &= \frac{L_1}{2}\int_\Omega \lvert \nabla Q_h(0) - \nabla \hat Q(0)\vert^2\, dx + \frac{L_2 + L_3}{2}\int_\Omega \lvert \Div Q_h(0) - \Div \hat Q(0)\rvert^2\, dx\\
        &\hspace{10ex} + \frac{\sigma}{2}\int_\Omega \mathcal{I}_h \lvert D_t^c Q_h(0) - D_t^c \hat Q(0)\rvert^2\, dx + \frac12 \int_\Omega \mathcal{I}_h \lvert r_h(0) - \hat r(0)\rvert^2\, dx\\
        &= O(\Delta t + h^2)
\end{align*}
so that
\[
    \mathcal{E}(U_h(t^{n+1})\vert \hat U(t^{n+1})) \le C(T)(\Delta t + h^2).
\]

This gives the desired convergence rates for $t = n\Delta t$, $n=0,1,\dots,N$. For other $t\in[0,T]$, using what we just showed and the Taylor expansion of $\hat{Q}$, we have for $n$ satisfying $n\Delta t \le t < (n+1)\Delta t$ that
\begin{align*}
    \|Q_h(t) - \hat{Q}(t)\|_{H^1} 
        &\le \|Q_h(t) - Q_h(t^n)\|_{H^1} + \|Q_h(t^n) - \hat{Q}(t^n)\|_{H^1} + \|\hat{Q}(t^n) - \hat{Q}(t)\|_{H^1}\\ 
        &= 0 + O(\Delta t^{1/2} + h) + O(\Delta t) = O(\Delta t^{1/2} + h),
\end{align*}
and similarly for the other error rates.
\end{proof}
\begin{remark}
Going through the proof above, we see that when  $\sigma=0$, we obtain a convergence rate of $O(\Delta t)$ with respect to time because all the lower order errors disappear. 
\end{remark}
\section{Convergence as $\sigma\to0$}\label{sec:zero_inertia}

Formally setting $\sigma=0$ in equation \eqref{eq:Qt}, we obtain the parabolic partial differential equation:
\begin{equation}
    \label{eq:Qt_no_inertia}\hat{Q}_t=L_1\Delta \hat{Q}+\frac{L_2+L_3}{2}\alpha(\hat{Q})-f(\hat{Q}).
\end{equation}
The goal of this section is to show that weak solutions of~\eqref{eq:Qt} converge to the strong solution of~\eqref{eq:Qt_no_inertia} as $\sigma\to 0$. To make clear the dependence of the solution on $\sigma$, we will use $\Qs$ to denote the solution of~\eqref{eq:qtensorflow} in this section. The main idea is to use the method of relative entropy that we used earlier in Theorem~\ref{thm:Existence_Strong_solution} to show weak-strong uniqueness. 
First, let us define strong solutions in the sense we will need in the following:
\begin{definition}
	\label{def:strongsolparabol}
	We call $\hQ$ a strong solution of~\eqref{eq:Qt_no_inertia} if the initial data $\hQ_0\in H^3(\dom)$ and $\hat{Q}$ satisfies~\eqref{eq:Qt_no_inertia} everywhere in $[0,T]\times \dom$ and 
	\begin{equation*}
	\hQ \in L^\infty(0,T;H^2(\dom))\cap L^2(0,T;H^3(\dom)),\quad \partial_t \hQ\in L^\infty(0,T;L^2(\dom))\cap L^2(0,T;H^2(\dom)),\quad \partial_t^2 \hQ\in L^2([0,T]\times\dom).
	\end{equation*}
\end{definition}
\begin{rem}[Regularity of $r$]
	\label{rem:regr}
	Using the equivalence of solutions of~\eqref{eq:Qt_no_inertia} and its reformulation~\eqref{eq:reformulation} (set $\sigma=0$), see~\cite[Lemma 5.2]{GWY2020}, we also obtain that for strong solutions, $\hr\in L^\infty([0,T]\times\dom)$. 
\end{rem}
	From~\cite{Yue2023} it follows that a strong solution of~\eqref{eq:Qt_no_inertia} exists with 	
	\begin{equation}\label{eq:aprioriregularity}
	\hQ \in L^\infty(0,T;H^2(\dom))\cap L^2(0,T;H^3(\dom)),\quad \partial_t \hQ\in L^\infty(0,T;L^2(\dom))\cap L^2(0,T;H^1(\dom)).
	\end{equation}
	as long as $\hQ_0\in H^2(\dom)$. The additional regularity can be proved by iterating the estimates for the derivatives of $\hQ$ as shown in the following lemma.

\begin{lemma}
\label{lem:Q_no_inertia_H4_estimate}
   Let $\hQ_0\in H^3(\Omega)$. Then, the strong solution $\hQ$ of~\eqref{eq:Qt_no_inertia} exists and satisfies $\partial_t^2\hQ\in L^2([0,T]\times\dom)$ and $\partial_t\hQ\in L^2(0,T;H^2(\dom))$.
\end{lemma}
\begin{proof}
	A sketch of the proof of this lemma can be found in Appendix~\ref{appen:lemmaproof}.
\end{proof}
We can now prove the main theorem of this section:
\begin{theorem}\label{thm:sigma_convergence}
    Let $\hQ_0,\Qs_0\in {H^3(\Omega)}$, $Q_{t,0}^\sigma\in L^2(\dom)$, where $Q_{t,0}^\sigma$ is the initial condition for the time derivative of $\Qs$.
    Then weak solutions $(\Qs,\rs)$ of \eqref{eq:Qt} satisfy 
    	\begin{multline*}
   \norm{\Grad\Qs(t)-\Grad\hQ(t)}_{L^2}^2+\norm{\rs(t)-\hr(t)}_{L^2}^2	\\
   \leq C\left(\norm{\Grad\Qs_0-\Grad\hQ_0}_{L^2}^2+\norm{\rs_0-\hr_0}_{L^2}^2+\frac{\sigma}{2}\norm{Q_{t,0}^\sigma-\partial_t \hQ(0)}_{L^2}^2+ C\sigma^2\right)\exp(C_T).
    \end{multline*}
    where $\hQ$ is the strong solution of~\eqref{eq:Qt_no_inertia} and $C_T$ a constant, depending on $T$ and the $L^\infty$-norm of $\partial_t\hQ$ and $\hr$. 
    Thus if 
    \begin{equation*}
    \norm{\Grad\Qs_0-\Grad\hQ_0}_{L^2}^2+\norm{\rs_0-\hr_0}_{L^2}^2\leq C \sigma^2,
    \end{equation*}
    and
    \begin{equation*}
    \norm{\Qs_{t,0}-\partial_t \hQ(0)}_{L^2}^2\leq C \sigma,
    \end{equation*}
    we obtain that weak solutions of \eqref{eq:Qt} converge to the strong solution of \eqref{eq:Qt_no_inertia} in the $H^1$-norm at a rate of $O({\sigma})$, i.e., for any $t\in [0,T]$,
    \begin{equation}
    \norm{\Qs(t)-\hQ(t)}_{H^1}\leq C {\sigma}.
    \end{equation}  
    Furthermore, we also have
    \begin{equation*}
    \norm{\partial_t\Qs-\partial_t\hQ}_{L^2([0,T]\times\dom)}\leq C{\sigma}.
    \end{equation*} 
\end{theorem}
\begin{proof}
	We use the method of relative entropy again. To do so, we denote $\Us=(\Qs,\rs)$ the solution of~\eqref{eq:Qt} with given $\sigma>0$ and $\hU=(\hQ,\hr)$ the solution of the (reformulation of) the parabolic system~\eqref{eq:Qt_no_inertia} and define the relative entropy
	\begin{equation}
	\label{eq:relentropysigma}
	\mathcal{E}(\Us|\hU)=E(\Us)-E(\hU)-dE(\hU)(\Us-\hU),
	\end{equation}
	where we recall the definition of the energy
	\begin{equation*}
	E(\Us)=\int_{\dom}\left[\frac{L_1}{2}|\Grad \Qs|^2+\frac{L_2+L_3}{2}|\Div \Qs|^2 +\frac{\sigma}{2}|\partial_t\Qs|^2+\frac12 (\rs)^2\right] dx. 
	\end{equation*}
	Therefore,
	\begin{equation*}
	\begin{split}
	E(\hU)&=\int_{\dom}\left[\frac{L_1}{2}|\Grad \hQ|^2+\frac{L_2+L_3}{2}|\Div \hQ|^2 +\frac12 \hr^2\right] dx,\quad \text{and}\\
	 dE(\hU)(\Us-\hU)&=\int_{\dom}\left[-\left(L_1\Delta\hQ+\frac{L_2+L_3}{2}\alpha(\hQ)\right):(\Qs-\hQ)+\hr(\rs-\hr)\right] dx,
	 \end{split}
	\end{equation*}
	and one can rewrite the relative entropy as
	\begin{equation*}
	\begin{split}
	\mathcal{E}(\Us|\hU)&=\int_{\dom}\left[\frac{L_1}{2}|\Grad (\Qs-\hQ)|^2+\frac{L_2+L_3}{2}|\Div (\Qs-\hQ)|^2 +\frac{\sigma}{2}|\partial_t\Qs|^2+\frac12 (\rs-\hr)^2\right] dx\\
	&=E(\Us)+E(\hU)-dE(\hU)\Us.
	\end{split}
	\end{equation*}
	Our goal is to get a bound on $\mathcal{E}(\Us|\hU)(t)$ in terms of the initial relative entropy $\mathcal{E}(\Us|\hU)(0)$ and terms that go to zero with $\sigma$. We first note that from the energy inequality, we have
	\begin{equation*}
	E(\Us(t))+E(\hU(t))\leq E(\Us_0)+E(\hU_0)-\int_0^t \left(\norm{\partial_t \hQ(s)}_{L^2}^2+\norm{\partial_t \Qs(s)}_{L^2}^2\right) ds.
	\end{equation*}
	Thus it remains to estimate $dE(\hU(t))\Us(t)$. Let $\theta_\epsilon=\mathbf{1}_{[\eta,t]}\star \omega_\epsilon$ be a regularized version of the indicator function of the interval $[\eta,t]$ where $\eta>\epsilon>0$ is a small number and $\omega_\epsilon(s)=\frac{1}{\epsilon}\omega(s/\epsilon)$ is a symmetric, nonnegative, smooth, compactly supported on $[-1,1]$ mollifier with $\int \omega(s)ds =1$. Then since $\Us$ is weakly continuous in time with values in $L^1(\dom)$ and $\hU$ is sufficiently smooth, we have
	\begin{equation*}
	dE(\hU(t))\Us(t) = dE(\hU(\eta))\Us(\eta)-\lim_{\epsilon\to 0}\int_0^T dE(\hU(s))\Us(s)\partial_s \theta_\epsilon(s) ds.
	\end{equation*}
	Therefore, we next compute an estimate for $\int_0^T dE(\hU(s))\Us(s)\partial_s \theta_\epsilon(s) ds$ and then pass $\epsilon\to 0$.
	Indeed, we have
	\begin{align*}
	&\int_0^T dE(\hU(s))\Us(s)\partial_s \theta_\epsilon(s) ds\\
	&\quad = \int_0^T\partial_s\theta_\epsilon\int_{\dom}\left(-L_1\Delta \hQ:\Qs -\frac{L_2+L_3}{2}\alpha(\hQ):\Qs+\hr\rs\right) dx ds\\
	&\quad = \int_0^T\theta_\epsilon\int_{\dom}\left(L_1\Delta \partial_s\hQ:\Qs +\frac{L_2+L_3}{2}\alpha(\partial_s\hQ):\Qs-\partial_s\hr\rs +L_1\Delta\hQ:\partial_s\Qs +\frac{L_2+L_3}{2}\alpha(\hQ):\partial_s\Qs-\hr\partial_s\rs\right) dx ds.
	\end{align*}
	The previous integration by parts in the time variable is well-defined since $\partial_t\Qs\in L^\infty(0,T;L^2(\dom))$ and $\partial_t\rs\in  L^2(0,T;L^1(\dom))$. We plug in the equations for $\partial_s\hr$, $\partial_s\rs$ and $L_1\Delta \hQ$ to obtain
	\begin{align*}
	&\int_0^T dE(\hU(s))\Us(s)\partial_s \theta_\epsilon(s) ds\\
	&\quad = \int_0^T\!\!\theta_\epsilon\!\int_{\dom}\Big(L_1\Delta \partial_s\hQ:\Qs +\frac{L_2+L_3}{2}\alpha(\partial_s\hQ):\Qs-P(\hQ):\partial_s\hQ\rs \\
 &\hphantom{\quad =\int_0^T\!\!\theta_\epsilon\!\int_{\dom}\Big(}\qquad+\left(\partial_s\hQ+P(\hQ)\hr\right):\partial_s\Qs -\hr P(\Qs):\partial_s \Qs\Big) dx ds\\
	&\quad = \int_0^T\!\!\theta_\epsilon\!\int_{\dom}\Big(L_1\Delta \partial_s\hQ:\Qs +\frac{L_2+L_3}{2}\alpha(\partial_s\hQ):\Qs+\partial_s\hQ:\partial_s\Qs \\
	&\hphantom{\quad = \int_0^T\!\!\theta_\epsilon\!\int_{\dom}\Big(}\qquad +P(\hQ):\left(\partial_s\Qs\hr-\partial_s\hQ \rs\right) -\hr P(\Qs):\partial_s \Qs\Big) dx ds\\
	&\quad = \int_0^T\!\!\theta_\epsilon\!\int_{\dom}\Big(-L_1\Grad \partial_s\hQ:\Grad\Qs -(L_2+L_3)\Div\partial_s\hQ \cdot\Div \Qs+\partial_s\hQ:\partial_s\Qs\\
 &\hphantom{\quad = \int_0^T\!\!\theta_\epsilon\!\int_{\dom}\Big(}\qquad
	+P(\hQ):\left(\partial_s\Qs\hr-\partial_s\hQ \rs\right) -\hr P(\Qs):\partial_s \Qs\Big) dx ds\\
&\quad = \int_0^T\!\!\int_{\dom}\Big(\theta_\epsilon \rs P(\Qs):\partial_s\hQ
+\theta_\epsilon \partial_s \hQ:\partial_s \Qs-\sigma\partial_s(\theta_\epsilon\partial_s \hQ):\partial_s\Qs \Big)dx ds\\
&\qquad+\int_0^T\!\!\theta_\epsilon\!\int_{\dom}\left(\partial_s\hQ:\partial_s\Qs
+P(\hQ):\left(\partial_s\Qs\hr-\partial_s\hQ \rs\right) -\hr P(\Qs):\partial_s \Qs\right) dx ds\\
&\quad = -\sigma\int_0^T\!\!\int_{\dom}\partial_s\theta_\epsilon\partial_s \hQ:\partial_s\Qs dx ds\\
&\qquad+\int_0^T\!\!\theta_\epsilon\!\int_{\dom}\left(2\partial_s\hQ:\partial_s\Qs
+(P(\hQ)- P(\Qs)):\left(\partial_s\Qs\hr-\partial_s\hQ \rs\right)-\sigma \partial_s^2 \hQ:\partial_s\Qs \right) dx ds.
\end{align*}
Letting $\epsilon\to 0$, the first term on the right hand side converges to	
\begin{equation*}
-\sigma\lim_{\epsilon\to 0}  \int_0^T\!\!\int_{\dom}\partial_s\theta_\epsilon\partial_s \hQ:\partial_s\Qs dx ds =- \sigma \int_{\dom} \partial_t\hQ(\eta):\partial_t \Qs(\eta)dx + \sigma \int_{\dom}\partial_t \hQ(t):\partial_t \Qs(t)dx,
\end{equation*}
since $\Qs_t$ is weakly continuous with values in $L^2$ and $\hQ_t$ is strongly continuous in time with values in $L^2(\dom)$. For the second term, we obtain since the integrand is uniformly integrable with respect to $\epsilon>0$,
\begin{equation*}
\begin{split}
&\lim_{\epsilon\to 0}\int_0^T\!\!\theta_\epsilon\!\int_{\dom}\left(2\partial_s\hQ:\partial_s\Qs
+(P(\hQ)- P(\Qs)):\left(\partial_s\Qs\hr-\partial_s\hQ \rs\right)-\sigma \partial_s^2 \hQ:\partial_s\Qs \right) dx ds\\
& = \int_\eta^t \int_{\dom}\left(2\partial_s\hQ:\partial_s\Qs
+(P(\hQ)- P(\Qs)):\left(\partial_s\Qs\hr-\partial_s\hQ \rs\right)-\sigma \partial_s^2 \hQ:\partial_s\Qs \right) dx ds.
\end{split}
\end{equation*}	
	Thus,
		\begin{align*}
	dE(\hU(t))\Us(t) &= dE(\hU(\eta))\Us(\eta)-\lim_{\epsilon\to 0}\int_0^T dE(\hU(s))\Us(s)\partial_s \theta_\epsilon(s) ds\\
	&=dE(\hU(\eta))\Us(\eta) +\sigma \int_{\dom} \partial_t\hQ(\eta):\partial_t \Qs(\eta)dx - \sigma \int_{\dom}\partial_t \hQ(t):\partial_t \Qs(t)dx\\
	&\quad - \int_\eta^t \int_{\dom}\left(2\partial_s\hQ:\partial_s\Qs
	+(P(\hQ)- P(\Qs)):\left(\partial_s\Qs\hr-\partial_s\hQ \rs\right)-\sigma \partial_s^2 \hQ:\partial_s\Qs \right) dx ds.
		\end{align*}
Again using the weak and strong continuity in time of $\partial_t\hQ$ and $\partial_t\Qs$, the weak continuity of $\Us$, and the strong continuity in time of $\hU$, we can send $\eta\to 0$ to obtain
		\begin{align*}
	dE(\hU(t))\Us(t) &= dE(\hU_0)\Us_0+\sigma \int_{\dom} \partial_t\hQ(0):\partial_t \Qs(0) dx - \sigma \int_{\dom}\partial_t \hQ(t):\partial_t \Qs(t)dx\\
	&\quad - \int_0^t \int_{\dom}\left(2\partial_s\hQ:\partial_s\Qs
	+(P(\hQ)- P(\Qs)):\left(\partial_s\Qs\hr-\partial_s\hQ \rs\right)-\sigma \partial_s^2 \hQ:\partial_s\Qs \right) dx ds.
	\end{align*}
	Therefore we get
	\begin{align*}
	\mathcal{E}(\Us|\hU)(t)&\leq \mathcal{E}(\Us_0|\hU_0)-\int_0^t\norm{\partial_s \Qs-\partial_s\hQ}_{L^2}^2 ds +\underbrace{\int_0^t\int_{\dom} (P(\Qs)- P(\hQ)):\left(\partial_s\Qs\hr-\partial_s\hQ \rs\right)dx ds}_{\text{I}}\\
	&-\underbrace{\sigma \int_{\dom} \partial_t\hQ(0):\partial_t \Qs(0) dx}_{\text{II}} + \underbrace{\sigma \int_{\dom}\partial_t \hQ(t):\partial_t \Qs(t)dx}_{\text{III}}-\underbrace{\sigma \int_0^t\int_{\dom}\partial_s^2\hQ:(\partial_s\Qs-\partial_s \hQ) dx ds}_{\text{IV}} \\
	&-\underbrace{\sigma \int_0^t\int_{\dom}\partial_s^2\hQ:\partial_s \hQ dx ds}_{\text{V}}.
	\end{align*}
	We continue to bound all the terms on the right hand side. For term I, using the Lipschitz continuity of $P$ (see \cite[Theorem 4.11]{GWY2020} or Lemma~\ref{lem:PLipschitz}), we obtain
	\begin{align*}
	\left|\text{I}\right|&\leq\int_0^t \int_{\dom} C |\Qs-\hQ|  \left|\hr\partial_s(\Qs-\hQ)+\partial_s \hQ(\hr-\rs)\right| dx ds\\
&\leq C \int_0^t\norm{\hr}_{L^\infty}\norm{\Qs-\hQ}_{L^2}\norm{\partial_s \Qs-\partial_s \hQ}_{L^2}ds +C\int_0^t\norm{\partial_s\hQ}_{L^\infty}\norm{\hr-\rs}_{L^2}\norm{\Qs-\hQ}_{L^2}ds\\
&\leq \frac14\int_0^t\norm{\partial_s \Qs-\partial_s\hQ}_{L^2}^2 ds + C\int_0^t \left(\norm{\hr}_{L^\infty}^2+\norm{\partial_s\hQ}_{L^\infty}\right)\left(\norm{\Qs-\hQ}_{L^2}^2+\norm{\hr-\rs}_{L^2}^2 \right)ds\\
&\leq \frac14\int_0^t\norm{\partial_s \Qs-\partial_s\hQ}_{L^2}^2 ds + C\int_0^t \left(\norm{\hr}_{L^\infty}^2+\norm{\partial_s\hQ}_{L^\infty}\right)\left(\norm{\Grad\Qs-\Grad\hQ}_{L^2}^2+\norm{\hr-\rs}_{L^2}^2 \right)ds,
	\end{align*}
where we used Poincar\'e inequality  for the last inequality.
For the second term, we have
\begin{equation*}
\text{II} = \frac{\sigma}{2}\int_{\dom}\left|\partial_t \hQ(0)-\partial_t \Qs(0)\right|^2dx -\frac{\sigma}{2}\int_{\dom}\left|\partial_t \hQ(0)\right|^2 dx-\frac{\sigma}{2}\int_{\dom}\left|\partial_t \Qs(0)\right|^2dx.
\end{equation*}
For term III, we have,
\begin{equation*}
\left|\text{III}\right|\leq \frac{\sigma}{2} \norm{\partial_t \Qs(t)}^2_{L^2}+\frac{\sigma}{2}\norm{\partial_t \hQ(t)}^2_{L^2}.
\end{equation*}
We combine this with term $\text{V}$, for which we have, using the chain rule and regularity of $\hQ$,
\begin{equation*}
\text{V} = - \frac{\sigma}{2}\left(\int_{\dom}|\partial_t \hQ(t)|^2dx -\int_{\dom}|\partial_t \hQ(0)|^2dx\right) =-\frac{\sigma}{2}\int_{\dom}|\partial_t \hQ(t)|^2dx +\frac{\sigma}{2}\int_{\dom}|\partial_t \hQ(0)|^2dx,
\end{equation*}
thus $|\text{III}|+\text{V}$ can be estimated as follows:
\begin{equation*}
|\text{III}|+\text{V} \leq \frac{\sigma}{2}\int_{\dom}|\partial_t \hQ(0)|^2dx+ \frac{\sigma}{2} \norm{\partial_t \Qs(t)}^2_{L^2}.
\end{equation*}
Now adding this to term $\text{II}$, we get
\begin{equation*}
\text{II} + |\text{III}|+\text{V} \leq  \frac{\sigma}{2} \norm{\partial_t \Qs(t)}^2_{L^2} + \frac{\sigma}{2}\norm{\partial_t \Qs(0)-\partial_t \hQ(0)}_{L^2}^2-\frac{\sigma}{2}\norm{\partial_t \Qs(0)}_{L^2}^2 .
\end{equation*}
Finally, we estimate term IV. Using Cauchy-Schwarz inequality, we have
\begin{equation*}
\left|\text{IV}\right|\leq \sigma^2\int_0^t\norm{\partial_s^2\hQ}_{L^2}^2+\frac14\int_0^t\norm{\partial_s\Qs-\partial_s\hQ}_{L^2}^2 ds\leq C\sigma^2+\frac14\int_0^t\norm{\partial_s\Qs-\partial_s\hQ}_{L^2}^2 ds,
\end{equation*}	
since $\partial_t^2\hQ\in L^2([0,T]\times\dom)$. 
In summary, we obtain
	
	\begin{align*}
	\mathcal{E}(\Us|\hU)(t)&\leq \mathcal{E}(\Us_0|\hU_0)-\frac12\int_0^t\norm{\partial_s \Qs-\partial_s\hQ}_{L^2}^2 ds+ C\sigma^2+ \frac{\sigma}{2} \norm{\partial_t\Qs(t)}^2_{L^2}+\frac{\sigma}{2}\norm{\partial_t \Qs(0)-\partial_t \hQ(0)}_{L^2}^2\\
	&\quad+C\int_0^t \left(\norm{\hr}_{L^\infty}^2+\norm{\partial_s\hQ}_{L^\infty}\right)\left(\norm{\Grad\Qs-\Grad\hQ}_{L^2}^2+\norm{\hr-\rs}_{L^2}^2 \right)ds-\frac{\sigma}{2}\norm{\partial_t \Qs(0)}_{L^2}^2 
	\end{align*}
	Defining
	\begin{equation*}
	\widetilde{\mathcal{E}}(\Us|\hU): = \int_{\dom}\left[\frac{L_1}{2}|\Grad (\Qs-\hQ)|^2+\frac{L_2+L_3}{2}|\Div (\Qs-\hQ)|^2 +\frac12 (\rs-\hr)^2\right] dx,
	\end{equation*}
	this implies
	\begin{align*}
	\widetilde{\mathcal{E}}(\Us|\hU)(t)&\leq \widetilde{\mathcal{E}}(\Us_0|\hU_0)-\frac12\int_0^t\norm{\partial_s \Qs-\partial_s\hQ}_{L^2}^2 ds+ C\sigma^2+\frac{\sigma}{2}\norm{\partial_t \Qs(0)-\partial_t \hQ(0)}_{L^2}^2\\
	&\quad+C\int_0^t \left(\norm{\hr}_{L^\infty}^2+\norm{\partial_s\hQ}_{L^\infty}\right) \widetilde{\mathcal{E}}(\Us|\hU)(s) ds.
	\end{align*}
	Using Gr\"onwall's inequality, we obtain
	\begin{align*}
		\widetilde{\mathcal{E}}(\Us|\hU)(t)&	\leq \left(\widetilde{\mathcal{E}}(\Us_0|\hU_0)+\frac{\sigma}{2}\norm{\partial_t \Qs(0)-\partial_t \hQ(0)}_{L^2}^2+ C\sigma^2\right)\exp\left( C\int_0^t \left(\norm{\hr}_{L^\infty}^2+\norm{\partial_s\hQ}_{L^\infty}\right)ds\right)\\
		&\leq \left(\widetilde{\mathcal{E}}(\Us_0|\hU_0)+\frac{\sigma}{2}\norm{\partial_t \Qs(0)-\partial_t \hQ(0)}_{L^2}^2+ C\sigma^2\right)\exp(C_T).
	\end{align*}
	We see that if $\tilde{\mathcal{E}}(\Us_0|\hU_0)\leq C\sigma^2$ and $\norm{\partial_t \Qs(0)-\partial_t \hQ(0)}_{L^2}^2\leq C \sigma$, then with Poincar\'e's inequality, we have that $\Qs(t)$ converges to $\hQ(t)$ in $H^1$ at rate $\sigma$ at best. We obtain the same rate for $\partial_t\Qs$ in $L^2([0,T]\times\dom)$. If one of the terms goes to zero at a slower rate, then the error is dominated by whichever of these terms is the largest. 
\end{proof}

\section{Numerical Results}\label{sec:num}
We now present some numerical results in 2D. The code to reproduce the results can be found on the Github repository \url{https://github.com/maxhirsch/Q-Tensor-Inertia}. In what follows, we define
\[
    E_1 = \|Q_0 - Q_0^\sigma\|_{H^1},\quad\text{and}\quad E_2 = \|Q_{t,0} - Q_{t,0}^\sigma\|_{H^1},
\]
the $H^1$ errors between the initial data and initial time derivative of solutions to the hyperbolic equation~\eqref{eq:Qt} and the parabolic equation~\eqref{eq:Qt_no_inertia}. We demonstrate numerical convergence rates in space and time then demonstrate the rate of convergence of the hyperbolic equation to the parabolic equation as we let $\sigma\to 0$.

As in \cite{GWY2020}, we consider the domain $\Omega = [0,2]^2$ and the parameters
\[
    L_1=0.001,\quad L_2 = L_3 = 0,\quad a = -0.2,\quad b = 1,\quad c = 1,\quad A_0 = 500,
\]
with the initial condition
\begin{equation}
    Q_0 = \boldsymbol{n}_0\boldsymbol{n}_0^\top - \frac{|\boldsymbol{n}_0|^2}{2}I_2,
\end{equation}
where
\begin{equation}
    \boldsymbol{n}_0 = \begin{bmatrix}
        x(2-x)y(2-y)\\
        \sin(\pi x)\sin(0.5\pi y)
    \end{bmatrix}.
\end{equation}
For the hyperbolic equation, we additionally take
\begin{equation}
    Q_{t, 0} = L_1\Delta Q(0) - r(0)P(Q(0)). 
\end{equation}

\subsection{Convergence test}

For the convergence in space and time, we take
\begin{equation}
    Q_0^\sigma = Q_0,\quad\text{and}\quad Q_{t,0}^\sigma = Q_{t,0}.
\end{equation}
We additionally set $\sigma = 0.025$.

\subsubsection{Refinement in space}

\begin{figure}[H]
    \centering
    \begin{tabular}{cc}
        \includegraphics[width=0.45\linewidth]{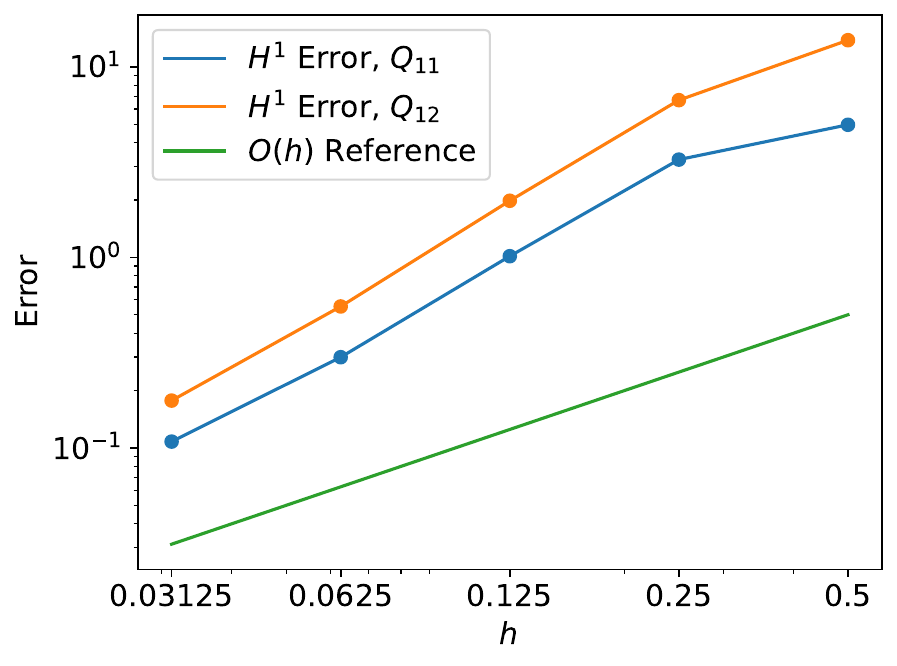} & \includegraphics[width=0.45\linewidth]{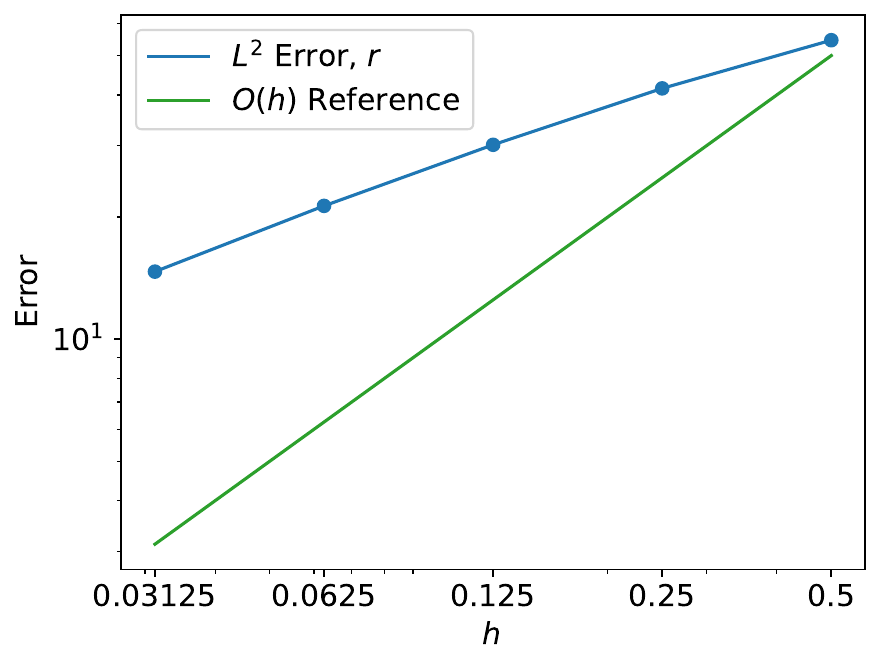}
    \end{tabular}
    \caption{Space refinement experiment $H^1$ errors for $(Q_h)_{11}$ and $(Q_h)_{12}$ and $L^2$ error for $r_h$. An $O(h)$ reference line is given in green without point markers.}
    \label{fig:space-refinement}
\end{figure}

For the convergence in space, we set $T=0.1$, $\Delta t = 1.25\times10^{-4}$, and the mesh size of the reference solution to be $h=2^{-8}$. We compute the $H^1$ errors of $(Q_h)_{11}$ and $(Q_h)_{12}$ from the reference solution for various $h$, as well as the $L^2$ error of $r_h$ from the reference solution. The results are given in Table~\ref{tab:space-refinement} and Figure~\ref{fig:space-refinement}.

\begin{table}[H]
\centering
    \begin{tabular}{|c|c|c|c|c|c|c|}
    \hline
     $h$ & $H^1$ error $Q_{11}$ & order for $Q_{11}$ & $H^1$ error $Q_{12}$ & order for $Q_{12}$ & $L^2$ error $r$ & order for $r$\\
     \hline
     0.5 & 4.97 & - & 13.80 & - & 54.56 & -\\
     \hline
     0.25 & 3.26 & 0.61 & 6.69 & 1.05 & 41.50 & 0.39\\
     \hline
     0.125 & 1.01 & 1.68 & 1.98 & 1.75 & 30.14 & 0.46\\
     \hline
     0.0625 & 0.30 & 1.76 & 0.55 & 1.84 & 21.31 & 0.50\\
     \hline
     0.03125 & 0.11 & 1.47 & 0.18 & 1.64 & 14.67 & 0.54\\
    \hline
    \end{tabular}
    \caption{Errors and rates for spatial refinement.}
    \label{tab:space-refinement}
\end{table}

Based on Theorem \ref{thm:convergence-rate}, we expect to see a convergence rate of $O(h)$, but instead, we see a faster rate for $Q_h$ and a slower rate for $r_h$. The average observed orders for $(Q_h)_{11}$ and $(Q_h)_{12}$ are $1.38$ and $1.57$, respectively, while the observed order for $r_h$ is $0.47$. The discrepancy with the theoretical estimates could be due to being in a pre-asymptotic regime.

\subsubsection{Refinement in time}

\begin{figure}[H]
    \centering
    \begin{tabular}{cc}
        \includegraphics[width=0.45\linewidth]{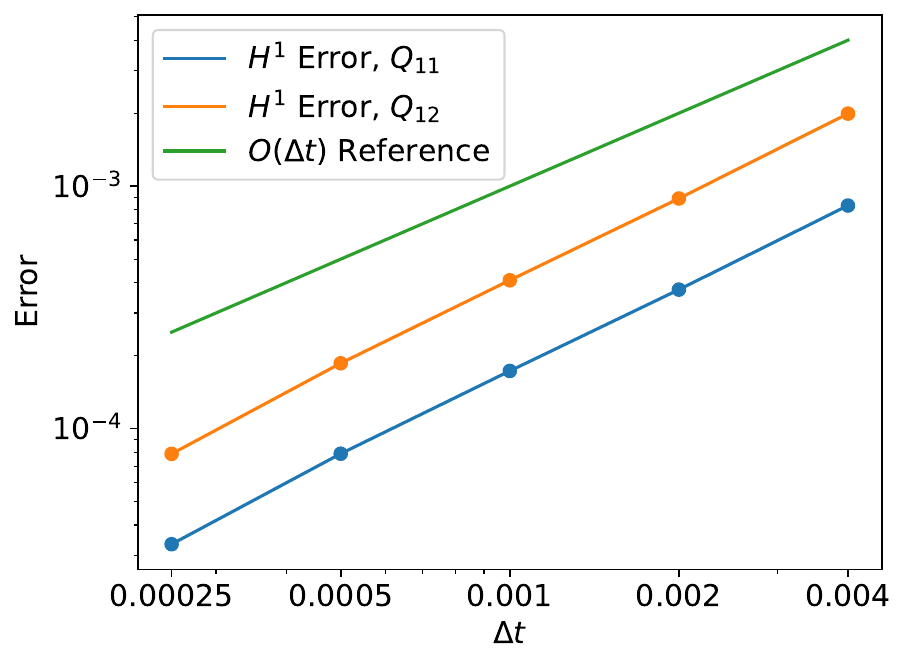} & \includegraphics[width=0.45\linewidth]{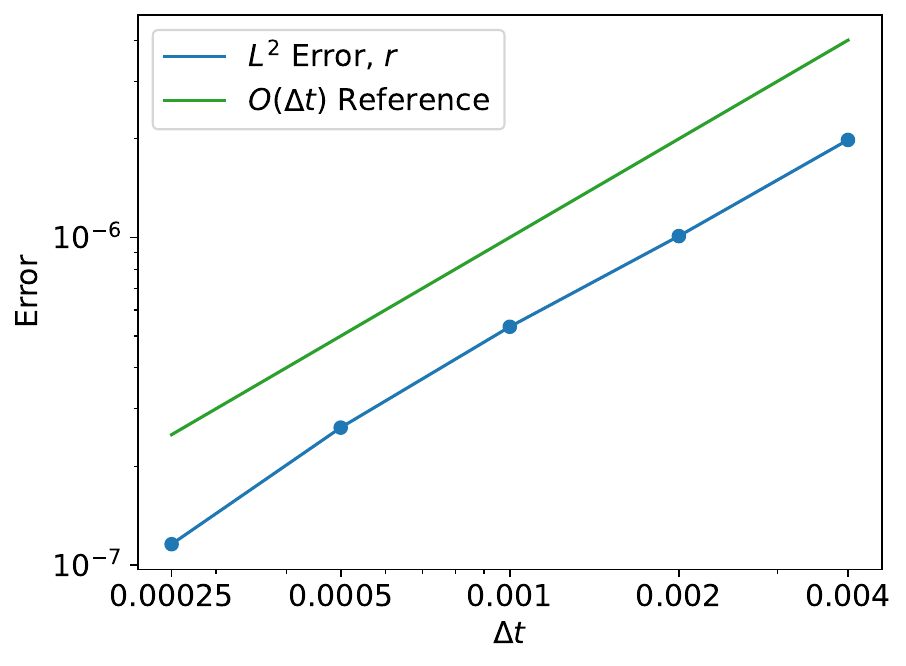}
    \end{tabular}
    \caption{Time refinement experiment $H^1$ errors for $(Q_h)_{11}$ and $(Q_h)_{12}$ and $L^2$ error for $r_h$. An $O(\Delta t)$ reference line is given in green without point markers.}
    \label{fig:time-refinement}
\end{figure}

For the convergence in time, we set $T = 0.1$, $h = 0.0625$, and the time step of the reference solution to be $\Delta t = 6.25\times10^{-5}$. We then compute the $H^1$ errors of $(Q_h)_{11}$ and $(Q_h)_{12}$ from the reference solution for various $\Delta t$, as well as the $L^2$ error of $r_h$ from the reference solution. The results are given in Table~\ref{tab:time-refinement} and Figure~\ref{fig:time-refinement}.

\begin{table}[H]
\centering
    \begin{tabular}{|c|c|c|c|c|c|c|}
    \hline
     $\Delta t$ & $H^1$ error $Q_{11}$ & order for $Q_{11}$ & $H^1$ error $Q_{12}$ & order for $Q_{12}$ & $L^2$ error $r$ & order for $r$\\
     \hline
     \hline
     $4\times10^{-3}$ & $8.32\times10^{-4}$ & - & $1.99\times10^{-3}$ & - & $1.98\times10^{-6}$ & -\\
     \hline
     $2\times10^{-3}$ & $3.74\times10^{-4}$ & 1.15 & $8.89\times10^{-4}$ & 1.16 & $1.01\times10^{-6}$ & 0.97\\
     \hline
     $1\times10^{-3}$ & $1.73\times10^{-4}$ & 1.11 & $4.09\times10^{-4}$ & 1.12 & $5.34\times10^{-7}$ & 0.92\\
     \hline
     $5\times10^{-4}$ & $7.88\times10^{-5}$ & 1.13 & $1.86\times10^{-4}$ & 1.14 & $2.63\times10^{-7}$ & 1.02\\
     \hline
     $2.5\times10^{-4}$ & $3.34\times10^{-5}$ & 1.24 & $7.87\times10^{-5}$ & 1.24 & $1.16\times10^{-7}$ & 1.18\\
    \hline
    \end{tabular}
    \caption{Errors and rates for time refinement.}
    \label{tab:time-refinement}
\end{table}

We see that we obtain a convergence rate of approximately $O(\Delta t)$, while Theorem~\ref{thm:convergence-rate} only guarantees a rate of $O(\Delta t^{1/2})$. We believe that Theorem~\ref{thm:convergence-rate} is not optimal as far as the established convergence rate in time is concerned, since the scheme is formally first order accurate and therefore such a rate should be expected for smooth solutions. However, we were not able to improve the rates in Theorem~\ref{thm:convergence-rate}.

\subsection{$\sigma$ convergence}

We now look at convergence rate of the hyperbolic equation~\eqref{eq:Qt} to the parabolic equation~\eqref{eq:Qt_no_inertia} as $\sigma\to 0$. We do this for various perturbations of the initial data and the initial time derivative. For the hyperbolic equation with $\sigma<1$, we take
\begin{equation}
    Q_0^\sigma = Q_0 + \frac{\sigma^{p_1}}{2}\begin{bmatrix}
        1 & 0\\ 0 & -1
    \end{bmatrix},
\end{equation}
with $p_1\in\{1/2,1,\infty\}$, and
\begin{equation}
    Q_{t,0}^\sigma = Q_{t,0} + \frac{\sigma^{p_2}}{2}\begin{bmatrix}
        1 & 0\\ 0 & -1
    \end{bmatrix},
\end{equation}
with $p_2\in\{1/2,\infty\}$. With $p_1 = \infty$ we have $Q_0^\sigma = Q_0$, and with $p_2 = \infty$ we have $Q_{t,0}^\sigma = Q_{t,0}$. Thus, the choices of $p_1 = 1/2, 1, \infty$ correspond to $E_1 = O(\sqrt{\sigma}), O(\sigma), 0$, respectively, and the choices of $p_2 = 1/2,\infty$ correspond to $E_2 = O(\sqrt{\sigma}), 0$, respectively. 

\begin{figure}[H]
\centering
\begin{tabular}{cc}
    \includegraphics[width=0.45\linewidth]{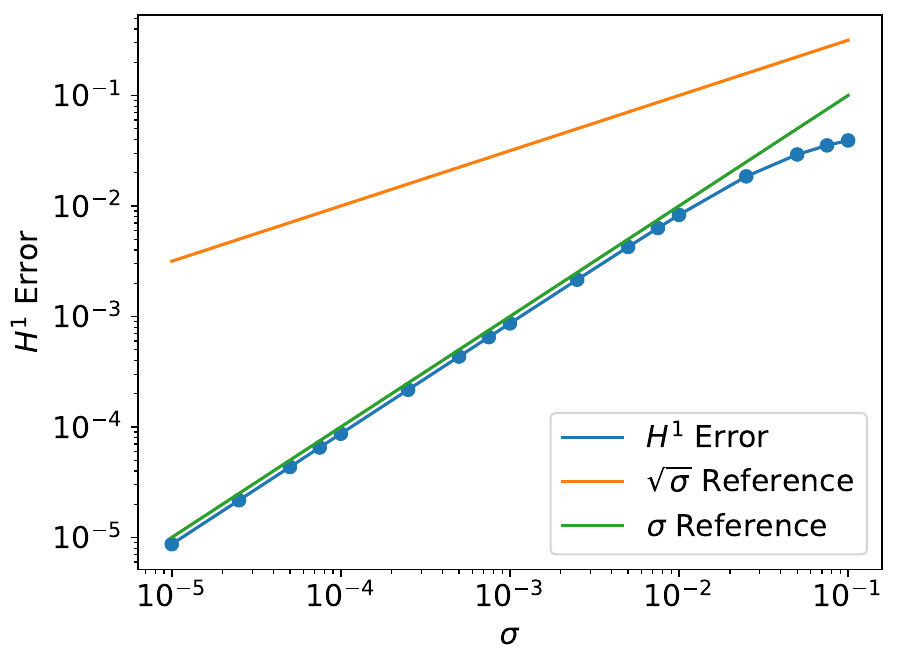}  & \includegraphics[width=0.45\linewidth]{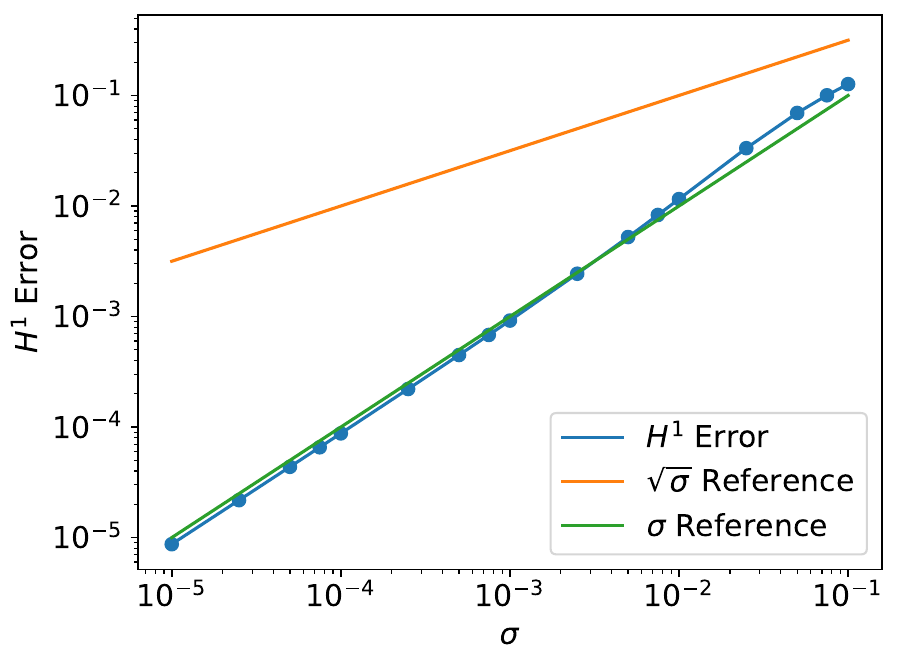}\\
    (a) $E_1 = 0$, $E_2 = 0$ & (b) $E_1 = 0$, $E_2 = O(\sqrt{\sigma})$\\
    \includegraphics[width=0.45\linewidth]{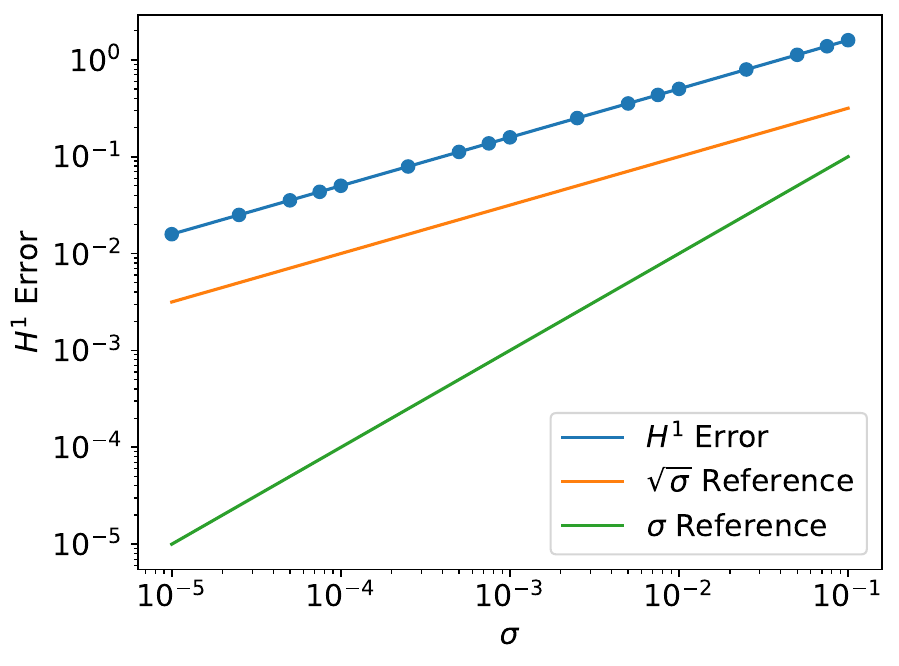}  & \includegraphics[width=0.45\linewidth]{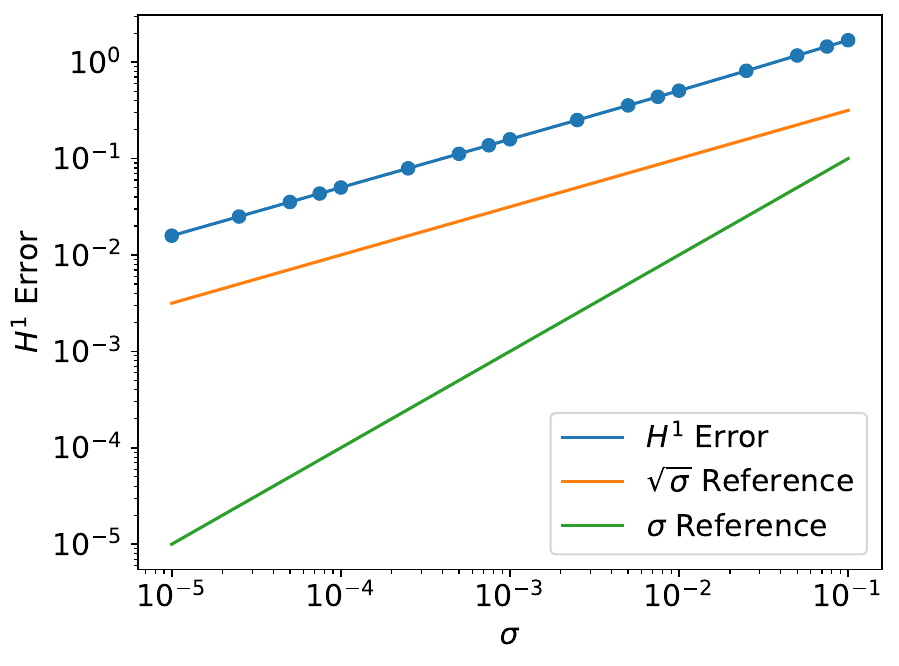}\\
    (c) $E_1 = O(\sqrt{\sigma})$, $E_2 = 0$ & (d) $E_1 = O(\sqrt{\sigma})$, $E_2 = O(\sqrt{\sigma})$\\
    \includegraphics[width=0.45\linewidth]{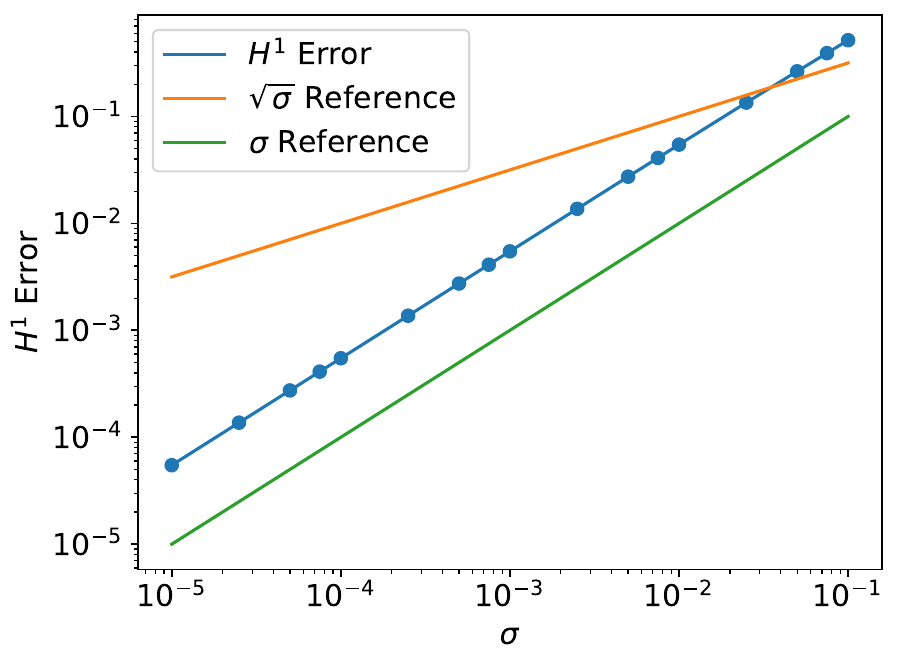}  & \includegraphics[width=0.45\linewidth]{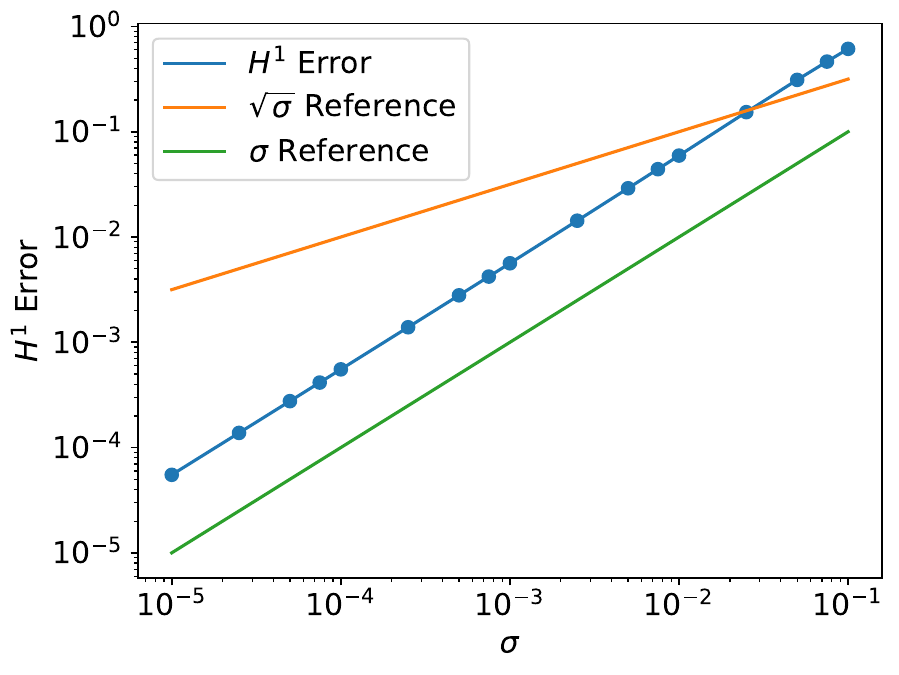}\\
    (e) $E_1 = O(\sigma)$, $E_2 = 0$ & (f) $E_1 = O(\sigma)$, $E_2 = O(\sqrt{\sigma})$\\
\end{tabular}
\caption{Convergence rates in $\sigma$ for various perturbations rates of the initial value and initial time derivative.}
\label{fig:sigma_experiment}
\end{figure}

We take $T=0.1$, $\Delta t = 10^{-5}$, and $h=0.125$, and for various values of $\sigma$, we compute the error
\[
    \|Q_h(T) - Q_h^\sigma(T)\|_{H^1},
\]
that is, the $H^1$ error at time $T$ between the discrete solution to the parabolic equation and the discrete solution to the hyperbolic equation for given $\sigma$. The results of the experiment are shown in Figure~\ref{fig:sigma_experiment}.

We see that the plots in Figure~\ref{fig:sigma_experiment} match the theoretical results in Theorem~\ref{thm:sigma_convergence}. Namely, since $E_2^2 \le C\sigma$, we have that the observed $H^1$ error between the two solutions will behave like $O(\max(\sigma, E_1))$.

\appendix
\section{Proof of Lemma \ref{lem:Q_no_inertia_H4_estimate}}\label{appen:lemmaproof}

\begin{lemma}
\label{appen:lem:Q_no_inertia_H4_estimate}
    Let $\hQ_0\in H^3(\Omega)$. Then the strong solution $\hQ$ of~\eqref{eq:Qt_no_inertia} exists and satisfies $\partial_t^2\hQ\in L^2([0,T]\times\dom)$ and $\partial_t\hQ\in L^2(0,T;H^2(\dom))$.
\end{lemma}
\begin{proof}
	We derive the necessary a priori estimates assuming that $\hQ$ is smooth, which allows us to do any formal manipulations. For a rigorous proof, one would have to construct a sequence of Galerkin approximations, prove the estimates at the level of the approximations and then pass to the limit similarly to what was done in Section~\ref{sec:strong_Solution}. Since the estimates are very similar, we skip that step.
    With this in mind, for our subsequent analysis, we will assume that $\hQ$ is a smooth function. Taking derivatives on both sides of \eqref{eq:Qt_no_inertia}, we get
    \begin{equation*}
        \hat{Q}_{tt}=L_1\Delta\hat{Q}_t+\frac{L_2+L_3}{2}\alpha(\hat{Q}_t)-f(\hat{Q})_t.
    \end{equation*}
    Taking the inner products with $-\hat{Q}_{tt}$ and integrating, we deduce,
    \begin{equation*}
        \begin{aligned}
           \frac{1}{2}\frac{d}{dt}\left[L_1\|\nabla \hat{Q}_t\|^2_{L^2} +(L_2+L_3) \|\Div \hat{Q}_t\|_{L^2}^2 \right]&= -\|\hat{Q}_{tt}\|^2_{L^2}-\int_\Omega\left[\frac{\delta f}{\delta Q}(\hat{Q})\hat{Q}_t\right]:\hat{Q}_{tt}\,dx\\
        &\leq  -\|\hat{Q}_{tt}\|_{L^2}^2+\left\|\frac{\delta f}{\delta Q}(\hat{Q})\right\|_{L^\infty}\|\hat{Q}_t\|_{L^2}\|\hat{Q}_{tt}\|_{L^2}.
        \end{aligned}
    \end{equation*}
    From~\cite{Yue2023}, we already get that 
   $\|\hat{Q}\|_{H^2}$ and $\|\hat{Q}_t\|_{L^2}$ are uniformly bounded, given $\hQ_0\in H^2(\Omega)$ (c.f.~\eqref{eq:aprioriregularity}). Since $\frac{\delta f}{\delta Q}(\hat{Q})$ is a quadratic polynomial of $\hat{Q}$, it can be controlled by a polynomial of $\|\hat{Q}\|_{L^\infty}$ and subsequently by $\|\hat{Q}\|_{H^2}$, by Sobolev embeddings. Together with the Cauchy-Schwarz inequality, we obtain
    \begin{equation*}
        \begin{aligned}
            \frac{1}{2}\frac{d}{dt}\left[L_1\|\nabla \hat{Q}_t\|^2_{L^2} +(L_2+L_3) \|\Div \hat{Q}_t\|_{L^2}^2 \right]&\leq -\|\hat{Q}_{tt}\|_{L^2}^2+C+\frac{1}{2}\|\hat{Q}_{tt}\|_{L^2}^2\leq C-\frac{1}{2}\|\hat{Q}_{tt}\|_{L^2}^2,
        \end{aligned}
    \end{equation*}
    for a constant $C$ depending on $\|Q_0\|_{H^2}$. Thus it follows by Gr\"onwall's inequality that
    \begin{equation}
        \|\nabla \hat{Q}_t(T)\|^2_{L^2}+\|\Div \hat{Q}_t(T)\|_{L^2}^2+\int_0^T\|\hat{Q}_{tt}(t)\|_{L^2}^2 dt\leq C(T+\|\nabla \hat{Q}_t(0)\|_{L^2}^2+\|\Div \hat{Q}_t(0)\|_{L^2}^2).
    \end{equation}
    Hence, if $Q_0\in H^3(\Omega)$, we obtain that $\int_0^T\|\hat{Q}_{tt}(t)\|^2_{L^2} dt$ remains bounded on any interval $[0,T]$.
    To prove the second estimate, we first note that the estimate of Lemma~\ref{lem:Q_H_control} can be shown for $\hQ_t$ instead of $\hQ$: We have
    \begin{equation*}
    L_1\Delta \hQ_t +\frac{L_2+L_3}{2}\alpha(\hQ_t)=f(\hQ)_t+H(\hQ)_t .
    \end{equation*}
    We take the inner product of this identity with $\Delta \hQ_t$ and integrate over the domain.
    \begin{equation*}
    L_1\norm{\Delta \hQ_t}_{L^2}^2 +\frac{L_2+L_3}{2}\int_\dom \alpha(\hQ_t):\Delta \hQ_t dx = \int_{\dom} f(\hQ)_t:\Delta\hQ_t +H(\hQ)_t:\Delta\hQ_t.
    \end{equation*}
    We estimate the second term on the left hand side
    \begin{align*}
    \left|\frac{L_2+L_3}{2}\int_{\dom}\alpha(\hQ)_t:\Delta\hQ_t dx \right|&\leq \frac{L_2+L_3}{2}\norm{\alpha(\hQ)_t}_{L^2}\norm{\Delta\hQ_t}_{L^2}\\
    &\leq C \frac{L_2+L_3}{2}\norm{\hQ_t}_{H^2}\norm{\Delta\hQ_t}_{L^2}\\
    &\leq C \frac{L_2+L_3}{2}\left(\norm{\hQ_t}_{L^2}+\norm{\Delta\hQ_t}_{L^2}\right)\norm{\Delta\hQ_t}_{L^2}\\
    &\leq \frac{L_1}{8}\norm{\Delta\hQ_t}_{L^2}^2+C \norm{\hQ_t}_{L^2}^2\\
    &\leq \frac{L_1}{8}\norm{\Delta\hQ_t}_{L^2}^2+C,
    \end{align*}
    where we used Lemma~\ref{lem:Q_H2_Laplacian} for $Q$ replaced by $\hQ_t$ and Assumption~\eqref{eq:L_1_large_assumption} and that we already know that $\hQ_t\in L^\infty(0,T;L^2(\dom))$. For the right hand side, we have
    \begin{equation*}
    \left|f(\hQ)_t:\Delta\hQ\right|\leq \frac{L_1}{4}\norm{\Delta\hQ}_{L^2}^2+C\norm{f(\hQ)_t}_{L^2}^2\leq \frac{L_1}{4}\norm{\Delta\hQ}_{L^2}^2+C\left(\norm{\hQ}_{L^\infty}^4+1\right)\norm{\hQ_t}_{L^2}^2\leq \frac{L_1}{4}\norm{\Delta\hQ}_{L^2}^2+C,
    \end{equation*}
    where we again used the already available bounds on $\hQ$. For the second term on the right hand side, we have
    \begin{equation*}
    \left|H(\hQ)_t:\Delta\hQ_t\right|\leq \frac{L_1}{2}\norm{\Delta\hQ_t}_{L^2}^2+C\norm{H(\hQ)_t}_{L^2}^2.
    \end{equation*}
    Combining these estimates, we obtain
    \begin{equation*}
    L_1\norm{\Delta\hQ_t}_{L^2}^2\leq \frac{7L_1}{8}\norm{\Delta\hQ_t}_{L^2}^2+C +C\norm{H(\hQ)_t}_{L^2}^2.
    \end{equation*}
    Thus,
    \begin{equation}
    \label{eq:DeltaQt}
    \norm{\Delta\partial_t\hQ}_{L^2}\leq C \left(1+\norm{\partial_t H(\hQ)}_{L^2}\right).
    \end{equation}
    Combining this with Lemma~\ref{lem:Q_H2_Laplacian} for $\hQ$ replaced by its time derivative, we obtain
    \begin{equation}
    \label{eq:H2Qtbound}
    \norm{\partial_t\hQ}_{H^2}\leq C\left(\norm{\partial_t\hQ}_{L^2}+\norm{\partial_t\Delta\hQ}_{L^2}\right)\leq C\left(1+\norm{\partial_t H(\hQ)}_{L^2}\right).
    \end{equation}
    The right hand side can be bounded by noticing that $\partial_t H(\hQ)=\partial_t^2\hQ$ which we have already shown is in $L^2([0,T]\times\dom)$, thus integrating~\eqref{eq:H2Qtbound} in time, we obtain the estimate on $\partial_t\hQ$ in $L^2(0,T;H^2(\dom))$.
     This concludes the proof.
        
\end{proof}

\begin{lemma}
\label{lem:thirdtimebound}
Let $Q_0\in H^3(\Omega)$ and $Q_{t,0} \in H^2(\Omega)$. Then we have for the solution of~\eqref{eq:qtensorflow} that
\[
    \|Q_{ttt}(t)\|_{L^2} \le C(\sigma)\quad \forall t\in[0,T],\quad\text{and}\quad \|Q_{ttt}\|_{L^2([0,T]\times\Omega)} \le C(\sigma).
\]
\end{lemma}
\begin{proof}
	We assume that the solution of~\eqref{eq:Qt} is smooth, otherwise, the same estimates can be derived at the level of a Galerkin approximation (cf. Section~\ref{subsection:galerkin}).
We differentiate \eqref{eq:Qt} twice and obtain
\begin{align*}
    Q_{ttt} &= -\sigma Q_{tttt} + \left(L_1\Delta Q_{tt} + \frac{L_2+L_3}{2}\alpha(Q)_{tt} - f(Q)_{tt}\right).
\end{align*}
Taking the inner product with $-Q_{ttt}$ on both sides and integrating then integrating by parts, we get
\[
    -\int_\Omega \lvert Q_{ttt}\rvert^2\, dx = \frac\sigma2 \frac{d}{dt} \int_\Omega \lvert Q_{ttt}\rvert^2\, dx + \frac{L_1}{2}\frac{d}{dt}\int_\Omega \lvert \nabla Q_{tt}\rvert^2\, dx + \frac{L_2 + L_3}{2}\frac{d}{dt}\int_\Omega \lvert\Div Q_{tt}\rvert^2\, dx + \int_\Omega f(Q)_{tt}:Q_{ttt}\, dx.
\]
We then have
\[
    \|f(Q)_{tt}\|_{L^2} \le C\left(\|Q_{tt}\|_{L^2} + \|Q\|_{L^\infty}\|Q_{tt}\|_{L^2} + \|Q_t\|_{L^4}^2 + \|Q\|_{L^\infty}^2\|Q_{tt}\|_{L^2} + \|Q\|_{L^\infty}\|Q_t\|_{L^4}^2\right) \le C(\sigma),
\]
where we used $\|Q\|_{L^\infty} \le C$ by Corollary \ref{cor:Q_H2_estimate} and the Sobolev embedding theorem, $\|Q_{tt}\|_{L^2} \le C(\sigma)$ by Corollary \ref{cor:Q_tt_estimate}, and $\|Q_t\|_{L^4} \le C$ again by Corollary \ref{cor:Q_H2_estimate} and the Sobolev embedding theorem.

Then using the Cauchy-Schwarz inequality and Young's inequality, we obtain
\begin{align*}
    &\frac12\frac{d}{dt} \int_\Omega \left[\sigma\lvert Q_{ttt}\rvert^2 + L_1\lvert \nabla Q_{tt}\rvert^2 + (L_2+L_3)\lvert \Div Q_{tt}\rvert^2\right]\, dx\\
    &= -\int_\Omega \lvert Q_{ttt}\rvert^2\, dx - \int_\Omega f(Q)_{tt}:Q_{ttt}\, dx\\
    &\le -\int_\Omega \lvert Q_{ttt}\rvert^2\, dx + \frac12 \int_\Omega \lvert Q_{ttt}\rvert^2\, dx + \frac12\int_\Omega \lvert f(Q)_{tt}\rvert^2\, dx\\
    &\le -\frac12\int_\Omega \lvert Q_{ttt}\rvert^2\, dx + C(\sigma).
\end{align*}
Integrating both sides with respect to time yields
\begin{align*}
    &\sigma\|Q_{ttt}(t)\|_{L^2}^2 + L_1\|\nabla Q_{tt}(t)\|_{L^2}^2 + (L_2+L_3)\|\Div Q_{tt}(t)\|_{L^2}^2 + \|Q_{ttt}\|_{L^2([0,t]\times\Omega)}^2\\
    &\le C(\sigma)T + \sigma\|Q_{ttt}(0)\|_{L^2}^2 + L_1\|\nabla Q_{tt}(0)\|_{L^2}^2 + (L_2+L_3)\|\Div Q_{tt}(0)\|_{L^2}^2.
\end{align*}
By using the equations
\begin{align*}
    \sigma Q_{ttt} &= -Q_{tt} + L_1\Delta Q_t + \frac{L_2+L_3}{2}\alpha(Q)_t - f(Q)_t,\\
    \sigma Q_{tt} &= -Q_t + L_1\Delta Q + \frac{L_2+L_3}{2}\alpha(Q) - f(Q),
\end{align*}
and using Sobolev inequalities, we see that the initial conditions give $Q_{ttt}(0) \in L^2(\Omega)$ and $Q_{tt}(0) \in H^1(\Omega)$. This along with the preceding inequality yields the desired result.
\end{proof}

\begin{lemma}
\label{lem:fourthtimebound}
Let $Q_0\in H^4(\Omega)$ and $Q_{t,0} \in H^3(\Omega)$. Then we have
\[
    \|Q_{tttt}(t)\|_{L^2} \le C(\sigma)\quad \forall t\in[0,T],\quad\text{and}\quad \|Q_{tttt}\|_{L^2([0,T]\times\Omega)} \le C(\sigma).
\]
\end{lemma}
\begin{proof}
Again, we assume the solution is smooth, if not, all calculations need to be done at the level of a Galerkin approximation. We differentiate \eqref{eq:Qt} three times and obtain
\begin{align*}
    Q_{tttt} &= -\sigma Q_{ttttt} + \left(L_1\Delta Q_{ttt} + \frac{L_2+L_3}{2}\alpha(Q)_{ttt} - f(Q)_{ttt}\right).
\end{align*}
Taking the inner product with $-Q_{tttt}$ on both sides and integrating then integrating by parts, we get
\[
    -\int_\Omega \lvert Q_{tttt}\rvert^2\, dx = \frac\sigma2 \frac{d}{dt} \int_\Omega \lvert Q_{tttt}\rvert^2\, dx + \frac{L_1}{2}\frac{d}{dt}\int_\Omega \lvert \nabla Q_{ttt}\rvert^2\, dx + \frac{L_2 + L_3}{2}\frac{d}{dt}\int_\Omega \lvert\Div Q_{ttt}\rvert^2\, dx + \int_\Omega f(Q)_{ttt}:Q_{tttt}\, dx.
\]
We then have
\begin{align*}
    \|f(Q)_{ttt}\|_{L^2} &\le C\Bigg(\|Q_{ttt}\|_{L^2} + \|Q\|_{L^\infty}\|Q_{ttt}\|_{L^2} + \|Q_t\|_{L^2}\|Q_{tt}\|_{L^2}\\ 
    &\hspace{10ex}+ \|Q\|_{L^\infty}^2\|Q_{tt}\|_{L^2} + \|Q\|_{L^\infty}\|Q_t\|_{L^2}\|Q_{tt}\|_{L^2} + \|Q_t\|_{L^6}^3\Bigg) \le C(\sigma),
\end{align*}
where we used Lemma \ref{lem:thirdtimebound}, that $\|Q\|_{L^\infty}, \|Q_t\|_{L^6} \le C$ by Corollary \ref{cor:Q_H2_estimate} and the Sobolev embedding theorem, and that $\|Q_{tt}\|_{L^2} \le C(\sigma)$ by Corollary \ref{cor:Q_tt_estimate}.

Then using the Cauchy-Schwarz inequality and Young's inequality, we obtain
\begin{align*}
    &\frac12\frac{d}{dt} \int_\Omega \left[\sigma\lvert Q_{tttt}\rvert^2 + L_1\lvert \nabla Q_{ttt}\rvert^2 + (L_2+L_3)\lvert \Div Q_{ttt}\rvert^2\right]\, dx\\
    &= -\int_\Omega \lvert Q_{tttt}\rvert^2\, dx - \int_\Omega f(Q)_{ttt}:Q_{tttt}\, dx\\
    &\le -\int_\Omega \lvert Q_{tttt}\rvert^2\, dx + \frac12 \int_\Omega \lvert Q_{tttt}\rvert^2\, dx + \frac12\int_\Omega \lvert f(Q)_{ttt}\rvert^2\, dx\\
    &\le -\frac12\int_\Omega \lvert Q_{tttt}\rvert^2\, dx + C(\sigma).
\end{align*}
Integrating both sides with respect to time yields
\begin{align*}
    &\sigma\|Q_{tttt}(t)\|_{L^2}^2 + L_1\|\nabla Q_{ttt}(t)\|_{L^2}^2 + (L_2+L_3)\|\Div Q_{ttt}(t)\|_{L^2}^2 + \|Q_{tttt}\|_{L^2([0,t]\times\Omega)}^2\\
    &\le C(\sigma)T + \sigma\|Q_{tttt}(0)\|_{L^2}^2 + L_1\|\nabla Q_{ttt}(0)\|_{L^2}^2 + (L_2+L_3)\|\Div Q_{ttt}(0)\|_{L^2}^2.
\end{align*}
By using the equations
\begin{align*}
    \sigma Q_{tttt} &= -Q_{ttt} + L_1\Delta Q_{tt} + \frac{L_2+L_3}{2}\alpha(Q)_{tt} - f(Q)_{tt},\\
    \sigma Q_{ttt} &= -Q_{tt} + L_1\Delta Q_{t} + \frac{L_2+L_3}{2}\alpha(Q)_{t} - f(Q)_{t},\\
    \sigma Q_{tt} &= -Q_t + L_1\Delta Q + \frac{L_2+L_3}{2}\alpha(Q) - f(Q),
\end{align*}
and using Sobolev inequalities, we see that the initial conditions give $Q_{tttt}(0) \in L^2(\Omega)$ and $Q_{ttt}(0) \in H^1(\Omega)$. This along with the preceding inequality yields the desired result.
\end{proof}

\vspace{2ex}
\noindent\textit{Data availability statement.} The code used to produce the numerical experiments in this article can be found on the Github repository \url{https://github.com/maxhirsch/Q-Tensor-Inertia} \cite{PaperCode}.

\bibliographystyle{abbrv}
\bibliography{related_literature}

\end{document}